\documentclass[12pt]{amsart}
\usepackage{etex}

\usepackage{amscd,latexsym,amsthm,amsfonts,amssymb,amsmath,amsxtra,mathdots,tikz}
\usepackage[mathscr]{eucal}
\usepackage{xcolor}
\usepackage[normalem]{ulem}
\usepackage[all,cmtip]{xy}
\usepackage{enumitem}
\usepackage{tabularx}
\usepackage{textcomp}
\usepackage{caption}

\renewcommand\theequation{\thesection.\arabic{equation}}

\newcommand{\BC}{{\mathbb {C}}}

\newcommand{\BR}{{\mathbb {R}}}

\newcommand{\BZ}{{\mathbb {Z}}}

\newcommand{\CO}{{\mathcal {O}}}

\newcommand{\CT}{{\mathcal {T}}}

\newcommand{\CX}{{\mathcal {X}}}

\newcommand{\Fe}{{\mathfrak {e}}}

\newcommand{\Fg}{{\mathfrak {g}}}

\newcommand{\Fl}{{\mathfrak {l}}}

\newcommand{\GL}{{\mathrm{GL}}}

\newcommand{\GSp}{{\mathrm{GSp}}}

\newcommand{\GSO}{{\mathrm{GSO}}}
\newcommand{\GSpin}{{\mathrm{GSpin}}}
\newcommand{\Spin}{{\mathrm{Spin}}}
\newcommand{\PGSO}{{\mathrm{PGSO}}}
\newcommand{\GHSpin}{{\mathrm{GHSpin}}}

\newcommand{\Hom}{{\mathrm{Hom}}}

\newcommand{\PGL}{{\mathrm{PGL}}}

\newcommand{\SL}{{\mathrm{SL}}}

\newcommand{\GU}{{\mathrm{GU}}}
\newcommand{\HSpin}{{\mathrm{HSpin}}}
\newcommand{\PGSp}{{\mathrm{PGSp}}}

\newcommand{\SO}{{\mathrm{SO}}}

\newcommand{\Sp}{{\mathrm{Sp}}}

\newcommand{\tr}{{\mathrm{tr}}}

\def\BR{{\mathbb R}}
\def\BC{{\mathbb C}}

\def\diag{{\rm diag}}

\def\eps{{\epsilon}}

\newtheorem{thm}{Theorem}[section]
\newtheorem{cor}[thm]{Corollary}

\newtheorem{prop}[thm]{Proposition}
\newtheorem {conj}[thm]{Conjecture}

\newtheorem {ques/conj}[thm]{Question/Conjecture}

\newtheorem{defn}[thm]{Definition}
\newtheorem{rmk}[thm]{Remark}

\newtheorem{lemma}[thm]{Lemma}

\makeatletter

\newcommand{\Rmnum}[1]{\expandafter\@slowromancap\romannumeral #1@}
\makeatother

\begin{document}
\renewcommand{\theequation}{\arabic{equation}}
\numberwithin{equation}{section}

\title{Multiplicities for Some Strongly Tempered Spherical Varieties}

\author{Chen Wan}
\address{Department of Mathematics \& Computer Science\\
Rutgers University – Newark\\
Newark, NJ 07102, USA}
\email{chen.wan@rutgers.edu}

\author{Lei Zhang}
\address{Department of Mathematics\\
National University of Singapore, Singapore}
\email{matzhlei@nus.edu.sg}

\date{}



\begin{abstract}
In this paper, we study the local multiplicity of 10 strongly tempered spherical varieties. We will formulate a uniform epsilon dichotomy conjecture for all these models regarding the unique distinguished element in tempered $L$-packets. Then we will prove this conjecture in many cases, including all the Archimedean cases.
\end{abstract}
\maketitle

% \tableofcontents

\section{Introduction and main results}

Let $F$ be a local field of characteristic 0, $G$ be a connected reductive group defined over $F$, $H$ be a connected closed subgroup of $G$, and $\chi$ be a unitary character of $H(F)$. Assume that $H$ is a spherical subgroup of $G$ (i.e. $H$ admitting an open orbit in the flag variety of
$G$). We say the spherical pair $(G,H)$ is reductive if $H$ is reductive. For every irreducible smooth 
\footnote{In the Archimedean case, smooth representations mean  Casselman--Wallach representations.
For the rest of the paper, all representations are assumed to be smooth.}
representation $\pi$ of $G(F)$, we define the multiplicity
$$m(\pi,\chi):=\dim(\Hom_{H(F)}(\pi,\chi)).$$
We say $\pi$ is $(H,\chi)$-distinguished (or just $H$-distinguished if the choice of $\chi$ is clear) if the multiplicity is nonzero. Also to simplify the notation we will use $m(\pi)$ instead of $m(\pi,\chi)$ to denote the multiplicity if the choice of $\chi$ is clear. One of the fundamental problems in the {\it Relative Langlands Program} is to study the multiplicity $m(\pi,\chi)$. In general, one expects the multiplicity $m(\pi,\chi)$ to be finite and to detect some functorial structures of $\pi$. We refer the reader to \cite{SV} for a detailed discussion of these kinds of problems.

Among all the spherical pairs, there is a special category called strongly tempered spherical pairs. More precisely, when $H$ is reductive, we say the pair $(G,H)$ is strongly tempered if all the matrix coefficients of tempered representations of $G(F)$ are integrable on $H(F)/Z_{G,H}(F)$ (here $Z_G$ is the center of $G$ and $Z_{G,H}=Z_G\cap H$). When $H$ is not reductive and if the model $(G,H)$ is the Whittaker induction (we refer the reader to Section 2.6 of \cite{Wan} for the definition of Whittaker induction) of a reductive spherical pair $(G_0,H_0)$, then we say the pair $(G,H)$ is strongly tempered if and only if  $(G_0,H_0)$ is strongly tempered. According to the general conjecture of Sakellaridis and Venkatesh in Conjecture 16.5.1 of \cite{SV}, for a strongly tempered spherical pair $(G,H)$, if we assume the spherical varieties $X=G/H$ does not have Type N spherical root (we refer the reader to Section 3.1 of \cite{SV} for the definition of spherical roots), then almost all the tempered local Vogan $L$-packets of $G(F)$ should contain at least one $(H,\chi)$-distinguished representation (i.e. almost all the tempered local Vogan $L$-packets are $(H,\chi)$-distinguished). Moreover, if the spherical variety only has one open Borel orbit over the local field $F$, then the general conjecture of Sakellaridis and Venkatesh predicts that almost all the tempered local Vogan $L$-packets of $G(F)$ should contain exactly one $(H,\chi)$-distinguished representation (this is usually called strong multiplicity one on  $L$-packets).

The most famous examples of strongly tempered spherical pairs without Type N root are the so call Gan--Gross--Prasad models $(\SO_{n+2k+1}\times \SO_n,\SO_n\ltimes U)$ and $(U_{n+2k+1}\times U_n,U_n\ltimes U)$. Here $U$ is some unipotent subgroup. For these cases, the local conjecture was formulated by Gan, Gross, and Prasad in Section 17 of \cite{GGP}. In it they not only conjectured the property of strong multiplicity one on generic $L$-packets, but they also conjectured about the unique distinguished element in those $L$-packets. More precisely, for each local $L$-packet $\Pi_\phi$, let $Z_\phi$ be the centralizer of the parameter and $S_\phi=Z_\phi/Z_{\phi}^{\circ}$ be its component group. The local Langlands conjecture states that there is a natural bijection between the $L$-packet   and the set of irreducible representations of $S_\phi$ (denoted by $\hat{S}_\phi$). In Section 17 of \cite{GGP}, they defined a quadratic character of $S_\phi$ using some local epsilon factor and conjectured that the unique distinguished element in a generic $L$-packet  is the one associated with this quadratic character. This is usually called the epsilon dichotomy conjecture.

In his pioneering works, Waldspurger developed a new method using local harmonic analysis to study the multiplicities. His idea is to first prove a local trace formula for the model which will imply a multiplicity formula $m(\pi,\chi)=m_{geom}(\pi,\chi)$. Here $m_{geom}(\pi,\chi)$ is defined via the Harish-Chandra character $\theta_\pi$ of $\pi$ and is called the geometric multiplicity. Then by using the multiplicity formula together with various relations of the Harish-Chandra characters of representations in a local $L$-packet, one can explicitly compute the multiplicity. In his works \cite{Wal1}, \cite{Wal2} and \cite{Wal3}, Waldspurger applied this idea to the orthogonal Gan--Gross--Prasad models and proved the local conjecture in the $p$-adic case. Later his method was adapted by Beuzart-Plessis (\cite{Beu1}, \cite{Beu2}, \cite{Beu3}, \cite{Beu4}), Wan (\cite{Wan15}, \cite{Wan16}, \cite{Wan17}), Beuzart-Plessis--Wan (\cite{BW}), Wan--Zhang (\cite{WZ1}, \cite{WZ2}), Luo (\cite{Luo}) for many other cases. Guided by all these works, in \cite{Wan}, the first author gave a uniform definition of the geometric multiplicity $m_{geom}(\pi,\chi)$ and proposed the conjectural multiplicity formula for all the spherical varieties.

In our previous paper \cite{WZ2}, we studied ten spherical pairs that are strongly tempered and without Type N spherical root. For each of the models, we computed its local relative character at unramified places. Our computation of the local relative characters shows that like the Gan--Gross--Prasad models case, the global period integrals of these models are related to the central values of certain automorphic $L$-functions $L(s,\pi,\rho_X)$ where $\rho_X$ is some finite dimensional representations of the $L$-group ${}^LG$ of symplectic type (see Table \ref{fig:1} for details). This allows us to formulate the Ichino-Ikeda type conjectures for these models. Locally, following the method of Waldspurger, we proved the multiplicity formulas in many cases. By using the multiplicity formulas, we proved the strong multiplicity one on  tempered $L$-packets (i.e. the summation of the multiplicities is equal to 1 over every tempered local $L$-packet) for these models. In particular, our results suggested that these models should have similar local and global behaviors as the Gan--Gross--Prasad models. In other words, many nice properties of the local multiplicities and global period integrals are not just enjoyed by the Gan--Gross--Prasad models, they can also be applied to general spherical varieties that are strongly tempered and without Type N spherical root.

Guided by this philosophy, it is natural to expect that for each of the models considered in our previous paper, like the Gan--Gross--Prasad models case, the unique distinguished element in the $L$-packet should be determined by the local epsilon factor $\epsilon(s,\pi,\rho_X)$. In this paper, we will formulate a uniform epsilon dichotomy conjecture for all these models. By studying the behaviors of the geometric multiplicities under parabolic induction and under endoscopy, we will prove the epsilon dichotomy conjecture in many cases including all the Archimedean cases.

\subsection{The conjectures and main results}

We recall the following table of spherical varieties from Section 1 of \cite{WZ2} (note that $\rho_X$ is a representation of the $L$-group ${}^L(G/Z_{G,H})$). Each model $(G,H)$ in the table is strongly tempered without Type N root and has a unique open Borel orbit. 

\begin{figure}[h!]
\begin{tabular}{| c | c | c |c| }
\hline
\textnumero & $G$ & $H$ &  $\rho_X$ \\
\hline
1 & $\GL_4\times \GL_2$ & $\GL_2\times \GL_2$ &  $(\wedge^2\otimes {\rm std}_{2})\oplus {\rm std}_{4} \oplus {\rm std}_{4}^\vee$ \\
\hline
2 &  $\GU_4\times \GU_2$ & $(\GU_2\times \GU_2)^0$ &  $(\wedge^2\otimes {\rm std}_{2})\oplus {\rm std}_{4} \oplus {\rm std}_{4}^\vee$  \\
\hline
3 &  $\GSp_6\times \GSp_4$ & $(\GSp_4\times \GSp_2)^0$ &  $\Spin_7 \otimes \Spin_5$ \\
\hline
4 &  $\GL_6$ & $\GL_2\ltimes U$ &  $\wedge^3$  \\
\hline
5 & $\GU_6$ & $\GU_2\ltimes U$ & $\wedge^3$  \\
\hline
6 & $\GSp_{10}$ & $\GL_2\ltimes U$ &  $\Spin_{11}$ \\
\hline
7 & $\GSp_{6}\times \GL_2$ & $\GL_2\ltimes U$ &  $\Spin_{7}\otimes {\rm std}_2$ \\
\hline
8 & $\GSO_8\times \GL_2$ & $\GL_2\ltimes U$ &  $\HSpin_8\otimes {\rm std}_{2}$ \\
\hline
9 & $\GSO_{12}$ & $\GL_2\ltimes U$ & $\HSpin_{12}$ \\
\hline
10 & $E_7$ & $\PGL_2\ltimes U$ & $\omega_7$  \\
\hline
\end{tabular}
\captionof{table}{}
\label{fig:1}
\end{figure}

Here ${\rm std}_n$ is the standard representation of $\GL_n(\BC)$ and ${\rm std}_n^\vee$ is its dual representation;
$\Spin_{2n+1}$ is the Spin representation of $\Spin_{2n+1}(\BC)$;
$\HSpin_{2n}$ is a half-Spin representation of $\Spin_{2n}(\BC)$;
$\omega_7$ is the 56 dimensional irreducible representation of $E_7(\BC)$. 
We refer readers to later sections for more details about $\rho_X$ in the unitary group cases (i.e. Models 2 and 5). We will also recall the definitions of all the models in later sections.

When $H$ is reductive, let $\chi=1$ be the trivial character of $H(F)$; when $H=H_0\ltimes U$ is not reductive, let $\chi=1\otimes \xi$ where $\xi$ is a generic character of $U(F)$ defined in our previous paper \cite{WZ2} (we will recall the definitions in later sections). Let $\pi$ be an irreducible representation of $G(F)$ whose central character is trivial on $Z_{G,H}(F)$, we want to study the multiplicity $m(\pi)=m(\pi,\chi)$. In our previous paper \cite{WZ2}, we have proved a multiplicity formula in the $p$-adic case and the complex case for all the models in the above table except the model associated to $E_7$. In the real case, we are only able to prove the multiplicity formula for the first four models. Moreover, we proved that if we assume the multiplicity formula and the local Langlands conjecture holds, then the summation of the multiplicities is equal to one over every tempered local Vogan $L$-packet (i.e. we have the strong multiplicity one on the $L$-packet). 

\begin{rmk}\label{distinguished is character}
In fact, as we explained in Section 9 of \cite{WZ2}, the multiplicity formula not only implies that the summation of the multiplicities is equal to one over every tempered local Vogan $L$-packet $\Pi_\phi$, it also implies that the unique distinguished element in the $L$-packet corresponds to a character of the component group $S_\phi$ (note that for Model 3 and Model 6-10 in Table \ref{fig:1}, the component group is not necessarily abelian). To be specific, combining the multiplicity formula and the character identity in the local Langlands correspondence, we have
$$\sum_{\pi} \dim(\chi_\pi)m(\pi)=1,$$
where $\pi$ runs over the representations in the packet and $\chi_\pi$ is the irreducible representation of the component group $S_\phi$ associated to $\pi$. Hence the irreducible representation of $S_\phi$ that corresponds to the unique distinguished element in the packet must be a character.
\end{rmk}

If $F=\BC$, the tempered $L$-packet only contains one element and its multiplicity is equal to 1. For the rest of this paper, we assume that $F\neq \BC$. To formulate the epsilon dichotomy conjecture, we need to define a character of the component group. Let $(G,H,\chi)$ be one of the models in the table above, and $\phi:W_F'\rightarrow {}^L G$ be a tempered Langlands parameter of $G$ (here ${}^LG$ is the $L$-group and $W_F'$ is the Weil-Deligne group) whose central character is trivial on $Z_{G,H}(F)$. This is equivalent to say that $\phi:W_F'\rightarrow {}^LG/Z_{G,H}$ is a tempered Langlands parameter of $G/Z_{G,H}$. Let $\Pi_\phi=\cup_{\alpha\in H^1(F,G/Z_{G,H})}\Pi_\phi(G_\alpha)$ be the associated tempered $L$-packet, $Z_\phi\subset \widehat{G/Z_{G,H}}$ be the centralizer of the parameter, and $S_\phi=Z_\phi/Z_{\phi}^{\circ}$ be the component group. The local Langlands conjecture states that we have a bijection between the $L$-packet $\Pi_\phi$ and the set $\hat{S}_\phi$ of irreducible representations of $S_{\phi}$. We refer the reader to Section \ref{sec:pre} for more details.

Next, we define a quadratic character of the component group $S_\phi$. We fix an additive character $\psi$ of $F$ and we use $V$ to denote the underlying space of the representation $\rho_X$ (i.e. $\rho_X:{}^LG/Z_{G,H}\rightarrow \GL(V)$). For $s\in S_\phi$, we will show in Lemma \ref{lem extended endoscopic triple} that there exists an elliptic extended endoscopic triple $(G',s',{}^L\eta)$ of $G$ such that the Langlands parameter $\phi$ factors through ${}^L\eta$ and $s'\in sZ_{\phi}^{\circ}$. For the model $(\GL_4\times \GL_2,\GL_2\times \GL_2)$, we require the lifting $s'$ to be of the form $\pm(I_4,I_2)$. Let $V_{s',-}$ be the $-1$ eigenspace of $V$ with respect to the operator $\rho_X(s')$. Since $s'$ commutes with $Im(\phi)$, the space $V_{s',-}$ is stable under $\rho_X(Im(\phi))$, this gives us a representation $\rho_{X,\phi,s'}$ of $W_F'$ on $V_{s',-}$, i.e. $\rho_{X,\phi,s'}:W_F'\rightarrow \GL(V_{s',-})$. If $(G,H)$ is not the two models associated to the unitary groups, we define
$$\omega_{\phi,H}(s')=\epsilon(\frac{1}{2},\rho_{X,\phi,s'},\psi).$$
In the two unitary group cases, we need to add some extra sign, we refer the reader to Section \ref{sec epsilon factor} for details. In Section \ref{sec epsilon factor}, we will show  that $\epsilon(\frac{1}{2},\rho_{X,\phi,s'},\psi)\in \{\pm 1\}$ and is independent of the choice of the additive character $\psi$ of $F$.

\begin{rmk}
For all the models in Table \ref{fig:1}, the Whittaker datum of $G$ is unique. So we don't need to discuss the choice of Whittaker datum.
\end{rmk}

\begin{rmk}
The extra sign in the two unitary group cases is an analogue of the extra sign in the Gan--Gross--Prasad model case (Section 6 of \cite{GGP}).
\end{rmk}

\begin{conj}[{\bf Epsilon Dichotomy Conjecture}] \label{main conj}
\begin{enumerate}
\item The function $\omega_{\phi,H}$ is well defined (i.e. it is independent of the choice of the elliptic extended endoscopic triple) and it is a quadratic character of $S_\phi$.
\item The unique $(H,\chi)$-distinguished element in the $L$-packet $\Pi_\phi$ is the one associated to the character $\omega_{\phi,H}$.
\end{enumerate}
\end{conj}

\begin{rmk}
Unlike the Gan--Gross--Prasad models case, the component group $S_\phi$ for most models in Table \ref{fig:1} is not necessarily a 2-group (not even necessarily abelian). But we still expect that the unique distinguished element in the $L$-packet corresponds to a quadratic character of $S_\phi$. 
\end{rmk}

There is also a weak form for Conjecture \ref{main conj}. For each model $(G,H)$ in Table \ref{fig:1} except the model $(\GU_4\times \GU_2, (\GU_2\times \GU_2)^0)$, the model has a unique pure inner form associated to the unique quaternion algebra $D$ defined over $F$. We will denote this model by $(G_D,H_D)$.

\begin{conj}\label{weak conjecture}
Let $(G,H)$ be a model in Table \ref{fig:1} that is not $(\GU_4\times \GU_2, (\GU_2\times \GU_2)^0)$ or $(\GU_6, \GU_2\ltimes U)$. The unique distinguished element in the packet $\Pi_\phi$ belongs to $\Pi_\phi(G)$ (resp. $\Pi_\phi(G_D)$) if and only if $\epsilon(\frac{1}{2},\Pi_\phi,\rho_X)=1$ (resp. $\epsilon(\frac{1}{2},\Pi_\phi,\rho_X)=-1$). 

For the model $(\GU_6, \GU_2\ltimes U)$, the unique distinguished element in the packet $\Pi_\phi$ belongs to $\Pi_\phi(G)$ (resp. $\Pi_\phi(G_D)$) if and only if $\eta_{E/F}(-1)\epsilon(\frac{1}{2},\Pi_\phi,\rho_X)=1$ (resp. $\eta_{E/F}(-1)\epsilon(\frac{1}{2},\Pi_\phi,\rho_X)=-1$) where $\eta_{E/F}$ is the quadratic character associated to the quadratic extension $E/F$ defining the unitary group.
\end{conj}

We refer the reader to Conjecture \ref{weak conjecture GU(4)xGU(2)} for the weak form of Conjecture \ref{main conj} for the model $(\GU_4\times \GU_2, (\GU_2\times \GU_2)^0)$. As in the Gan--Gross--Prasad model case, we also expect Conjecture \ref{main conj} and \ref{weak conjecture} to be true for all the generic local $L$-packets. 

\begin{prop}
Conjecture \ref{main conj} implies Conjecture \ref{weak conjecture} (or Conjecture \ref{weak conjecture GU(4)xGU(2)} if $(G,H)=(\GU_4\times \GU_2, (\GU_2\times \GU_2)^0)$). Moreover, if $\Pi_\phi$ is a discrete L-packet of $G$ with $|\Pi_\phi(G)|=1$, then Conjecture \ref{main conj} holds for the L-packet $\Pi_\phi$ if and only if Conjecture \ref{weak conjecture} (or Conjecture \ref{weak conjecture GU(4)xGU(2)} if $(G,H)=(\GU_4\times \GU_2, (\GU_2\times \GU_2)^0)$) holds for the L-packet $\Pi_\phi$.
\end{prop}

\begin{proof}
We will only consider the model $(E_7,\PGL_2\ltimes U)$. The argument for the other models is similar. In this case, $\hat{G}=E_{7,sc}(\BC)$ is the simply connected form of $E_7$ and its center $Z_{\hat{G}}$ is isomorphic to $\BZ/2\BZ$. Let $z$ be the nontrivial element in the center. 

We first prove the first part of the conjecture. Let $\Pi_\phi$ be a tempered L-packet of $G$ and assume that Conjecture \ref{main conj} holds for $\Pi_\phi$. We need to prove Conjecture \ref{weak conjecture} for $\Pi_\phi$. By our assumption we know that $\omega_{\phi,H}$ is a well defined quadratic character of $S_\phi$ and it corresponds to the unique distinguished element in the packet. Let $s_0\in S_\phi$ be the image of $z$. By our definition of $\omega_{\phi,H}$, we have (note that $\rho_X(z)=-\textbf{1}_{V}$)
$$\omega_{\phi,H}(s_0)=\omega_{\phi,H}(z)=\epsilon(\frac{1}{2},\Pi_\phi,\rho_X).$$
Then Conjecture \ref{weak conjecture} follows from the fact that $\omega_{\phi,H}$ corresponds to an element in $\Pi_\phi(G)$ (resp. $\Pi_\phi(G_D)$) if and only if $\omega_{\phi,H}(z)=1$ (resp. $\omega_{\phi,H}(z)=-1$). 

For the second part, if $\Pi_\phi$ is a discrete L-packet of $G$ with $|\Pi_\phi(G)|=1$, we have 
$$S_\phi=Z_\phi=Z_{\hat{G}}=\BZ/2\BZ,\; |\Pi_\phi(G_D)|=1.$$ 
In this case, it is clear that $\omega_{\phi,H}$ is well defined and is a quadratic character of $S_\phi$. Moreover, we have
$$\omega_{\phi,H}(1)=1,\;\omega_{\phi,H}(z)=\epsilon(\frac{1}{2},\Pi_\phi,\rho_X).$$
In particular, $\omega_{\phi,H}$ corresponds to the unique element in $\Pi_\phi(G)$ (resp. $\Pi_\phi(G_D)$) if and only if $\epsilon(\frac{1}{2},\Pi_\phi,\rho_X)=1$ (resp. $\epsilon(\frac{1}{2},\Pi_\phi,\rho_X)=-1$). This proves the proposition.
\end{proof}

In this paper, by using the multiplicity formulas and the character identities in the local Langlands conjecture, we will prove Conjecture \ref{main conj} in many cases, including all the archimedean cases. In order to state our result, we first need to define a partial order for the models in Table \ref{fig:1}.

\begin{defn}\label{smaller model}
Consider the following diagram of the models in Table \ref{fig:1}
$$(\GU_6,\GU_2\ltimes U)\rightarrow (\GU_4\times \GU_2, \GU_2\ltimes U)$$

\[
\xymatrix{
(\GSp_4\times \GL_2\times \GL_2,(\GL_2\times \GL_2)^0) & (\GSp_6\times \GSp_4,(\GSp_4\times \GSp_2)^0) \ar[l] \ar@{->}[d]\\
(\GSp_{10},\GL_2\ltimes U)\ar[r] \ar@{->}[d] & (\GSp_{6}\times \GL_2,\GL_2\ltimes U)    \\
(\GSO_{8}\times \GL_2,\GL_2\ltimes U)\ar@{->}[d] &(\GSO_{12},\GL_2\ltimes U) \ar[l] \ar@{->}[d] \\
(\GL_4\times \GL_2,\GL_2\times \GL_2) & (\GL_6,\GL_2\ltimes U) \ar[l] 
}
\] 
We say a model $(G,H)$ is smaller than another model $(G',H')$ if there is a line connecting these two models with the arrow pointing to $(G,H)$. 
\end{defn}
For example, there are two models smaller than $(\GSp_{10},\GL_2\ltimes U)$: $(\GSp_6\times \GL_2,\GL_2\ltimes U)$ and $(\GSO_8\times \GL_2, \GL_2\ltimes U)$.

\begin{rmk}
The only models in the above diagram that do not appear in Table \ref{fig:1} are the models $(\GU_4\times \GU_2, \GU_2\ltimes U)$ and $(\GSp_4\times \GL_2\times \GL_2,(\GL_2\times \GL_2)^0)$. The reason is that up to some finite isogeny, these two models are essentially the Gan-Gross-Prasad models $(\SO_6\times \SO_3,\SO_3\ltimes U)$  and $(\SO_5\times \SO_4,\SO_4)$. 
Although the epsilon dichotomy conjecture is known for Gan-Gross-Prasad model of special orthogonal groups and unitary groups, it is still open for these two models. We refer the reader to Sections \ref{sec GU} and \ref{sec GSp} for the details about these two models (note that each of these two models also has a unique pure inner form associated to the quaternion algebra $D$). The analogue of Conjecture \ref{weak conjecture} for these two models are stated in Conjecture \ref{weak conjecture GU(4)} and Conjecture \ref{weak conjecture GSp_4}.
\end{rmk}

\begin{thm}\label{main theorem}
Let $(G,H)$ be one of the models in Table \ref{fig:1} that is not $(E_7,\PGL_2\ltimes U)$. Assume that the multiplicity formula holds for the model $(G,H)$ and for all the models smaller than $(G,H)$. Also assume that the local Langlands conjecture (see Section \ref{sec L-packets}) holds for $G$. Let $\Pi_\phi=\cup_{\alpha\in H^1(F,G/Z_{G,H})}\Pi_\phi(G_\alpha)$ be a tempered $L$-packet.
\begin{enumerate}
\item If $(G,H)=(\GL_4\times \GL_2,\GL_2\times \GL_2)$ or $(\GU_4\times \GU_2, (\GU_2\times \GU_2)^0)$, there is no model smaller than $(G,H)$. Assume that the central character of $\Pi_\phi$ is trivial (not just trivial on $Z_{G,H}(F)$), or assume that $\Pi_\phi$ is not a discrete $L$-packet with $|\Pi_\phi(G)|=1$ (this is always the case when $F$ is Archimedean). Then Conjecture \ref{main conj} holds for the $L$-packet $\Pi_\phi$.
\item If $(G,H)$ is one of Models 3-9 of Table \ref{fig:1}, Assume that the weaker form of the conjecture (i.e. Conjecture \ref{weak conjecture}, Conjecture \ref{weak conjecture GU(4)} or Conjecture \ref{weak conjecture GSp_4}) holds for all the models that are smaller than $(G,H)$. Assume that $\Pi_\phi$ is not a discrete $L$-packet with $|\Pi_\phi(G)|=1$ (this is always the case when $F$ is Archimedean). Then Conjecture \ref{main conj} holds for the $L$-packet $\Pi_\phi$.
\end{enumerate}
\end{thm}

In particular, the above theorem proves Conjecture \ref{main conj} in the Archimedean case because by induction we can always assume that the weak form of the conjecture holds for all the models that are smaller than $(G,H)$.

\begin{cor}
Let $F=\BR$ and let $(G,H)$ be one of the models in Table \ref{fig:1} that is not $(E_7,\PGL_2\ltimes U)$. Assume that the multiplicity formula holds for the model $(G,H)$ and for all the models smaller than $(G,H)$. Then Conjecture \ref{main conj} and \ref{weak conjecture} (or Conjecture \ref{weak conjecture GU(4)xGU(2)} if $(G,H)=(\GU_4\times \GU_2, (\GU_2\times \GU_2)^0)$) hold for all tempered L-packets of $G.$
\end{cor}

\begin{rmk}
For Models 1--4, the multiplicity formula $m(\pi)=m_{geom}(\pi)$ has been proved for both the $p$-adic case and the real case (\cite{Wan15},  \cite{Wan16}, \cite{WZ2}, \cite{PWZ19}). For Models 5-9, the multiplicity formula $m(\pi,\chi)=m_{geom}(\pi,\chi)$ has been proved for the $p$-adic case (\cite{WZ1}, \cite{WZ2}). One can prove the multiplicity formula for the smaller models $(\GU_4\times \GU_2, \GU_2\ltimes U)$ and $(\GSp_4\times \GL_2\times \GL_2,(\GL_2\times \GL_2)^0)$ by a very similar argument. For the remaining cases, one needs to solve two technical issues in order to prove the multiplicity formulas (see the proof of Theorem 9.8 of \cite{WZ2} for details). We will recall the multiplicity formula for all the models in later sections.
\end{rmk}

When the packet $\Pi_\phi$ is discrete with $|\Pi_\phi(G)|=1$, since Conjecture \ref{main conj} is equivalent to Conjecture \ref{weak conjecture} (or Conjecture \ref{weak conjecture GU(4)xGU(2)}), the above theorem implies that if we assume that Conjecture \ref{weak conjecture} (or Conjecture \ref{weak conjecture GU(4)xGU(2)}) holds for all the models that are smaller than $(G,H)$ and for the model $(G,H)$, then Conjecture \ref{main conj} holds for $(G,H)$.

\begin{cor}\label{main cor}
Let $(G,H)$ be one of the models in Table \ref{fig:1} that is not $(E_7,\PGL_2\ltimes U)$. Assume that the multiplicity formula holds for the model $(G,H)$ and for all the models smaller than $(G,H)$. Also assume that the local Langlands conjecture holds for $G$. Moreover, assume that the weaker form of the conjecture (i.e. Conjecture \ref{weak conjecture}, Conjecture \ref{weak conjecture GU(4)xGU(2)}, Conjecture \ref{weak conjecture GU(4)} or Conjecture \ref{weak conjecture GSp_4}) holds for both the model $(G,H)$ and for all the models smaller than $(G,H)$. Then Conjecture \ref{main conj} holds for $(G,H)$.
\end{cor}

In fact, if we assume the weaker form of the conjecture holds for the model $(G,H)$ and the model is not $(\GSp_{10},\GL_2\ltimes U)$, we can prove the weaker form of the conjecture for all the models smaller than it.

\begin{thm}\label{thm weak conjecture smaller models}
Let $(G,H)$ be one of the models in Table \ref{fig:1} that is not $(\GSp_{10},\GL_2\ltimes U)$ or $(E_7,\PGL_2\ltimes U)$. Assume that the multiplicity formula holds for the model $(G,H)$ and for all the models smaller than $(G,H)$. Also assume that the local Langlands conjecture holds for $G$. Then the weaker form of the conjecture (i.e. Conjecture \ref{weak conjecture}, Conjecture \ref{weak conjecture GU(4)xGU(2)}) for the model $(G,H)$ implies the weaker form of the conjecture (i.e. Conjecture \ref{weak conjecture}, Conjecture \ref{weak conjecture GU(4)} or Conjecture \ref{weak conjecture GSp_4}) for all the models smaller than $(G,H)$.
\end{thm}

Combine the above theorem with Corollary \ref{main cor}, we get the following corollary.

\begin{cor}
Let $(G,H)$ be one of the models in Table \ref{fig:1} that is not $(\GSp_{10},\GL_2\ltimes U)$ or $(E_7,\PGL_2\ltimes U)$. Assume that the multiplicity formula holds for the model $(G,H)$ and for all the models smaller than $(G,H)$. Also assume that the local Langlands conjecture holds for $G$. Then Conjecture \ref{main conj} is equivalent to the weaker form of the conjecture (i.e. Conjecture \ref{weak conjecture}, Conjecture \ref{weak conjecture GU(4)xGU(2)}) for the model $(G,H)$.
\end{cor}

\begin{rmk}
For the model $(E_7,\PGL_2\ltimes U)$, by assuming the local Langlands conjecture and the multiplicity formula,  we can prove similar results as above. However, the smaller models are more complicated in this case. We will postpone the discussion of this model to Section \ref{sec E7} (see Theorem \ref{main theorem for E7} and \ref{thm weak conj for E7}).

For the model $(\GSp_{10},\GL_2\ltimes U)$, Conjecture \ref{weak conjecture} for this model will only imply Conjecture \ref{weak conjecture} for one of the smaller models $(\GSp_6\times \GL_2,\GL_2\ltimes U)$. It does not imply Conjecture \ref{weak conjecture} for the model $(\GSO_8\times \GL_2,\GL_2\ltimes U)$. We refer the reader to Remark \ref{GSp10 weak conjecture remark} for details.
\end{rmk}

\begin{rmk}
By Theorem 1.1 of \cite{GZ}, we also know that when $F$ is p-adic, Conjecture \ref{weak conjecture} for the model $(\GSO_{12},\GL_2\ltimes U)$ (resp. $(\GSO_8\times \GL_2,\GL_2\ltimes U)$) implies Conjecture \ref{weak conjecture} for the model $(\GSp_{10},\GL_2\ltimes U)$ (resp. $(\GSp_6\times \GL_2,\GL_2\ltimes U)$).
\end{rmk}

Now let's briefly explain the proof of Theorem \ref{main theorem} and \ref{thm weak conjecture smaller models}. Let $(G,H)$ be one of the models in Table \ref{fig:1}. Our assumption in Theorem \ref{main theorem} (i.e. the packet is not discrete with only one element) tells us the $L$-packet $\Pi_\phi$ is either the parabolic induction of some $L$-packet of a Levi subgroup $M$ of $G$, or the endoscopic transfer of some $L$-packet of an elliptic endoscopic group $G'$ of $G$. 
Then by studying the behaviors of the geometric multiplicities under parabolic induction and under endoscopic transfer, we can relate the multiplicity of the $L$-packet $\Pi_\phi$ to the multiplicities of certain models related to $M$ or $G'$. The only exception would be the model $(\GU_4\times \GU_2,(\GU_2\times \GU_2)^0)$ which requires some extra effort. This is because unlike the rest 9 models in the table, the model $(\GU_4\times \GU_2,(\GU_2\times \GU_2)^0)$ has more than one pure inner form. We refer the reader to Sections \ref{sec strategy} and \ref{sec GU} for details. There are two types of models we will get under this process,  either some model that has already been studied (e.g. Whittaker models, Gan--Gross--Prasad models), or  a model that is smaller than $(G,H)$ under Definition \ref{smaller model}. This is why we make the assumption on smaller models in Theorem \ref{main theorem}. This proves Conjecture \ref{weak conjecture} for the $L$-packet and also proves Theorem \ref{thm weak conjecture smaller models}. As a result, we get a formula of the epsilon factor $\epsilon(1/2,\Pi_\phi,\rho_X)$ in terms of the Harish-Chandra character of the $L$-packet. Combining the formula of epsilon factor with the formula of the geometric multiplicity under endoscopy and the definition of $\omega_{\phi,H}$, we can prove Theorem \ref{main theorem}.

On the other hand, for the discrete $L$-packets with only one element,
if one can prove the same formula of the epsilon factor $\epsilon(\frac{1}{2},\phi,\rho_X)$ in terms of the Harish-Chandra character of the $L$-packet, then
one can prove the conjecture for this case. 
% if one wants to prove the conjecture when the packet is discrete with only one element, one still hopes to prove a formula of the epsilon factor $\epsilon(\frac{1}{2},\phi,\rho_X)$ in term of the Harish-Chandra character of the $L$-packet. 
In the Gan--Gross--Prasad models case, such a formula was proved by Waldspurger and Beuzart-Plessis using the Rankin-Selberg integrals of the general linear groups and the theory of twisted endoscopy. However their method cannot be directly applied to our cases in Table \ref{fig:1} because the representations $\rho_X$ in Table \ref{fig:1} are more complicated than the Gan--Gross--Prasad models case (in particular, the Langlands functoriality between $G$ and $\GL_{\dim(\rho_X)}$ is not of twisted endoscopic type). In ongoing work, we are trying to completely prove Conjecture \ref{main conj} by studying the multiplicity of certain models related to the Rankin-Selberg integrals. 

Finally, we want to emphasize that the assumption that the spherical variety has a unique open Borel orbit is essential. Without this assumption, the multiplicity will no longer satisfy strong multiplicity one on the $L$-packet, and the summation of the multiplicities over the $L$-packet should be equal to the number of open Borel orbits (although we still expect the distinguished elements to be related to certain epsilon factors). In another ongoing work, we are trying to formulate the epsilon dichotomy conjecture for general strongly tempered spherical varieties.

\subsection{Organization of the paper}
In Section \ref{sec:pre}, we will give the basic notation of the paper and define the epsilon factors that appeared in our conjecture. We will also recall the structure of the $L$-packet under the local Langlands conjecture. In Section \ref{sec strategy}, we will explain the strategy of the proof. In Section \ref{Section GL}, we will consider the models $(\GL_4\times \GL_2,\GL_2\times \GL_2)$ and $(\GL_6,\GL_2\ltimes U)$. These are the easiest cases since the $L$-packet contains at most one element for each group (in particular Conjecture \ref{main conj} is equivalent to Conjecture \ref{weak conjecture}). In Section \ref{sec GU}, we will consider the models $(\GU_4\times \GU_2,(\GU_2\times \GU_2)^0)$ and $(\GU_6,\GU_2\ltimes U)$. The model $(\GU_4\times \GU_2,(\GU_2\times \GU_2)^0)$ is the most complicated case in this paper. In Section \ref{sec GSO and GSp}, we will consider the models $(\GSO_8\times \GL_2,\GL_2\ltimes U)$, $(\GSO_{12},\GL_2\ltimes U)$, $(\GSp_6\times \GL_2,\GL_2\ltimes U)$ and $(\GSp_{10},\GL_2\ltimes U)$. In Section \ref{sec GSp}, we will consider the model $(\GSp_6\times \GSp_4, (\GSp_4\times \GSp_2)^0)$. In Section \ref{sec E7}, we will consider the model $(E_7,\PGL_2\ltimes U)$.

\subsection{Acknowledgements}
 We would like to thank Wee Teck Gan, Tasho Kaletha, Diana Shelstad, Jun Yu for the helpful discussions.
The work of the first author is partially supported by the NSF grant DMS-2000192 and DMS-2103720. 
The work of the second author is partially supported by AcRF Tier 1 grants A-0004274-00-00 and A-0004279-00-00 of the National University of Singapore. We thank an anonymous referee for the helpful comments and corrections.

\section{Preliminary}\label{sec:pre}

\subsection{Notation}
Let $F$ be a local field of characteristic 0, and $\psi:F\rightarrow \BC^{\times}$ be a nontrivial additive character. Let $G$ be a connected reductive group defined over $F$, $\Fg$ be the Lie algebra of $G$, $Z_G$ be the center of $G$, and $A_G$ be the maximal split torus of $Z_G$. We use $G_{ss}$, $G_{reg}$ (resp. $\Fg_{ss}$, $\Fg_{reg}$) to denote the set of semisimple and regular semisimple elements of $G$ (resp. $\Fg$). For $x\in G_{ss}$ (resp. $X\in \Fg_{ss}$), let $Z_G(x)$ (resp. $Z_G(X)=G_X$) be the centralizer of $x$ (resp. $X$) in $G$ and let $G_x$ be the neutral component of $Z_G(x)$. Similarly, for any abelian subgroup $T$ of $G$, let $Z_G(T)$ be the centralizer of $T$ in $G$ and let $G_T$ be the neutral component of $Z_G(T)$. We say $x\in G_{ss}(F)$ is elliptic if $A_G=A_{G_x}$. We use $G(F)_{ell}$ (resp. $G(F)_{reg,ell}=G(F)_{ell}\cap G_{reg}(F)$) to denote the set of elliptic elements (resp. regular elliptic elements) of $G(F)$. For $x\in G_{ss}(F)$ (resp. $X\in \Fg_{ss}(F)$), let 
$$D^G(x)=|\det(1-Ad(x))_{|\Fg/\Fg_x}|_F \;(\text{resp.}\; D^G(X)=|\det(ad(X))_{|\Fg/\Fg_X}|_F)$$ 
be the Weyl determinant where $|\cdot|_F$ is the normalized absolute value on $F$. Finally, we use $\CT(G)$ (resp. $\CT_{ell}(G)$) to denote a set of representatives of maximal tori (resp. maximal elliptic tori) of $G(F)$. For $T\in \CT(G)$, we use $W(G,T)$ to denote the Weyl group.

For a quasi-character $\theta$ on $G(F)$ and $x\in G_{ss}(F)$, let $c_\theta(x)$ be the average of the regular germs of $\theta$ at $x$. For a regular nilpotent orbit $\CO$ of $\Fg_x(F)$, let $c_{\theta,\CO}(x)$ be the regular germ of $\theta$ at $x$ with respect to $\CO$. We refer the reader to Section 4 of \cite{Beu3} for the definition and basic properties of quasi-characters. If $\pi$ is a smooth finite length representation of $G(F)$, we use $\theta_\pi$ to denote its Harish-Chandra character (which is a quasi-character) and we use 
$$c_\pi(x)=c_{\theta_\pi}(x),\;c_{\pi,\CO}(x)=c_{\theta_\pi,\CO}(x)$$  
to denote the regular germs. If $M$ is a Levi subgroup of $G$ and $\theta^M$ is a quasi-character on $M(F)$, we use $i_{M}^{G}(\theta^M)$ to denote the parabolic induction of $\theta^M$ to $G(F)$. It is a quasi-character of $G(F)$. We refer the reader to Sections 3.4 and 4.7 of \cite{Beu3} for details. The following two propositions will be used in later sections.

\begin{prop}\label{regular germs} (Proposition 4.5.1 of \cite{Beu3})
Let $\theta$ be a quasi-character on $G(F)$ and $x\in G_{ss}(F)$. If $G_x$ is not quasi-split, then $c_\theta(x)=0$. If $G_x$ is quasi-split, let $B_x\subset G_x$ be a Borel subgroup and $T_{qs,x}\subset B_x$ be a maximal torus. Then we have
$$D^G(x)^{1/2}c_\theta(x)=|W(G_x,T_{qs,x})|^{-1} \lim_{x'\in T_{qs,x}(F)\rightarrow x} D^G(x')^{1/2}\theta(x'). $$
\end{prop}

\begin{prop}\label{germ parabolic induction} (Proposition 4.7.1 of \cite{Beu3})
Let $\theta=i_{M}^{G}(\theta^M)$ and $x\in G_{ss}(F)$. Let $\CX_M(x)$ be a set of representatives for the $M(F)$-conjugacy classes of elements in $M(F)$ that are $G(F)$-conjugated to $x$. Then $D^G(x)^{1/2}c_\theta(x)$ is equal to
$$|Z_G(x)(F):G_x(F)|  \sum_{y\in \CX_M(x)} |Z_M(y)(F):M_y(F)|^{-1} D^M(y)^{1/2} c_{\theta^M}(y).$$
In particular, $c_{\theta}(x)=0$ if the set $\CX_M(x)$ is empty.
\end{prop}

Lastly, we recall the endoscopic transfer of quasi-characters. Let $(G',s',{}^L\eta)$ be an extended endoscopic triple of $G$ (defined in Section 1.3 of \cite{K}), and let $\theta$ (resp. $\theta'$) be a quasi-character on $G(F)$ (resp. $G'(F)$). Assume that $\theta'$ is stable. We say $\theta$ is the endoscopic transfer of $\theta'$ if 
$$D^G(x)^{1/2}\theta(x)=\sum_y D^{G'}(y)^{1/2} \Delta(y,x) \theta'(y)$$
for all $x\in G_{reg}(F)$. Here $y$ runs over regular semisimple stable conjugacy classes of $G'(F)$,
and $\Delta(y,x)$ is the transfer factor defined in Section 2.3 of \cite{K} (the definition is the same as the one in \cite{LS} if $G$ is quasi-split). 
Note that for given $x$ there are only finitely many stable conjugacy classes $y$ such that the transfer factor is nonzero. 
In later sections, we will write down the explicit formula of the transfer factors in some special cases.

\begin{rmk}
For all the groups in Table \ref{fig:1}, the endoscopic group is always an L-group. Hence we only need to consider the extend endoscopic triple in this paper instead of the general endoscopic datum defined in \cite{LS}.
\end{rmk}

\subsection{The groups}
In this subsection, we will define all the reductive groups that will be used in later sections. Let $E=F(\sqrt{\eps})$ be a quadratic extension of $F$, $\eta_{E/F}$ be the quadratic character associated to $E$,  $N_{E/F}$ (resp. $\tr_{E/F}$) be the norm map (resp. trace map), and $x\rightarrow \bar{x}$ be the Galois action on $E$. Denote  $w_{n}$ to  be the longest Weyl element of $\GL_n$. Define the quasi-split even unitary similitude group $\GU_{n,n}(F)$ to be
\begin{equation}\label{eq:GU-n}
\GU_{n,n}=\{ g\in Res_{E/F}\GL_{2n}\mid  {}^{t}\bar{g}w_{2n}g=l(g) w_{2n}\}	
\end{equation}
where $l(g)\in F^{\times}$ is the similitude factor of $g$. If $F=\BR$, we can also define the groups
\begin{equation}\label{eq:GU-n}
\GU_{p,q}=\{ g\in Res_{\BC/\BR}\GL_{n}\mid  {}^{t}\bar{g}\cdot \diag(I_p,-I_q)g=l(g)\cdot \diag(I_p,-I_q)\}	
\end{equation}
for $p+q=n$ with $p\neq q$. 
To be compatible with the standard notation in the Archimedean case, when $F$ is $p$-adic, we will use $\GU_{n+1,n-1}=\GU_{n-1,n+1}$ to denote the non quasi-split inner form of $\GU_{n,n}$.

Let 
$$J_2'=\begin{pmatrix}0&-1\\1&0\end{pmatrix},\;J_{2n}'=\begin{pmatrix}0&J_{2n-2}'\\ J_2' &0 \end{pmatrix},\;L_{4}=\begin{pmatrix}0&J_2'\\-J_2'&0 \end{pmatrix}$$ and
$L_{4n}=\begin{pmatrix}0&0&J_2'\\0&L_{4n-4}&0\\-J_2'&0&0 \end{pmatrix}$. 
Define
\begin{align*}
&\GSO_{4n}=\{g\in \GL_{4n} \mid {}^tgL_{4n}g =l(g)L_{4n},\;\det(g)=l(g)^{2n}\},\\ 
&\GSO_{2n}(D)=\{g\in \GL_{2n}(D) \mid {}^t\bar{g} J_{2n}'g=l(g)J_{2n}'\}.
\end{align*}
Here $D/F$ is unique quaternion algebra over $F$. We also define
$$\GSp_{2n}=\{g\in \GL_{2n} \mid {}^t gJ_{2n}'g =l(g)J_{2n}'\},$$
$$ \GSp_{n}(D)=\{g\in \GL_n(D)\mid {}^t\bar{g}w_ng=l(g)w_n\}.
$$
We can also define $\PGSO_{4n}=\GSO_{4n}/\GL_1$, $\PGSO_{2n}(D)=\GSO_{2n}(D)/\GL_1$, $\PGSp_{2n}=\GSp_{2n}/\GL_1$ and $\PGSp_n(D)=\GSp_n(D)/\GL_1$. Also for any two similitude groups $GH_1$ and $GH_2$, we will use $G(H_1\times H_2)=(GH_1\times GH_2)^0$ to denote the subgroup
$$\{(g_1,g_2)\in GH_1\times GH_2 \mid l(g_1)=l(g_2)\}$$
of $GH_1\times GH_2$. And we use $(GH_1\times GH_2)^1$ to denote the subgroup
$$\{(g_1,g_2)\in GH_1\times GH_2 \mid l(g_1)=l(g_2)^{-1}\}$$
of $GH_1\times GH_2$.

Lastly, we use $E_7$ to denote the split adjoint reductive group of Type $E_7$ and we use $E_{7,sc}$ to denote the split simply connect reductive group of Type $E_7$.

\subsection{The local Langlands conjecture}\label{sec L-packets}
In this subsection we recall the local Langlands conjecture in Conjecture E of \cite{K}. Let $G$ be a quasi-split reductive group defined over $F$ and let $\{G_\alpha \mid \alpha\in H^1(F,G)\}$ be the set of pure inner forms of $G$. Let $\Pi_{irr,temp}(G_{\alpha})$ be the set of irreducible tempered representations of $G_{\alpha}(F)$. The local Langlands conjecture states that 
$$\cup_{\alpha\in H^1(F,G)}\Pi_{irr,temp}(G_\alpha)$$ 
is a disjoint union of finite sets (i.e. the local tempered Vogan $L$-packets)
$$\cup_{\phi} \Pi_{\phi}$$
where $\phi$ runs over all the tempered $L$-parameters of $G$ and $$\Pi_{\phi}=\cup_{\alpha\in H^1(F,G)} \Pi_{\phi}(G_\alpha)$$
consists of a finite number of tempered representations with $\Pi_{\phi}(G_\alpha)\subset \Pi_{irr,temp}(G_{\alpha})$ such that the following conditions hold.

\begin{itemize}
\item There is a unique generic element in $\Pi_\phi(G)$ with respect to any Whittaker datum of $G$.
\item For the given Whittaker datum, there is a bijection between $\hat{S_\phi}$, the set of irreducible representations of the component group $S_\phi=Z_\phi/Z_{\phi}^{\circ}$ of the Langlands parameter $\phi$, and $\Pi_\phi$ (denoted by $\pi\leftrightarrow \chi_\pi$) satisfies the following conditions. 
\begin{itemize}
\item The trivial character of $S_\phi$ corresponds to the unique generic element of $\Pi_\phi(G)$ with respect to the given Whittaker datum.
\item For $\alpha\in H^1(F,G)$, the distribution character $$\theta_{\Pi_{\phi}(G_\alpha)}=\sum_{\pi\in \Pi_\phi(G_{\alpha})}\dim(\chi_\pi)\theta_\pi$$ 
is stable. Moreover,  $\iota(G_{\alpha})\theta_{\Pi_\phi(G_\alpha)}$ is the transfer of $\theta_{\Pi_\phi(G)}$ where $\iota(G_\alpha)$ is the Kottwitz sign.
\item For any $\alpha\in H^1(F,G)$ and $\pi \in \Pi_\phi(G_\alpha)$, the restriction of the central character of $\chi_{\pi}$ to $Z(\hat{G})^{\Gamma_F}$ is equal to $\chi_\alpha$. Here $\chi_\alpha$ is the character of $Z(\hat{G})^{\Gamma_F}$ associated to $\alpha$ via the Kottwitz isomorphism. 
Note that the representation $\chi_\pi$ of the component group can be viewed as a representation of the centralizer $Z_\phi$ of the image of $\phi$, the group $Z(\hat{G})^{\Gamma_F}$ belongs to the center of $Z_\phi$ and hence it makes sense to talk about the restriction of the central character of $\chi_\pi$ to $Z(\hat{G})^{\Gamma_F}$.
\item For $s\in S_\phi$ and for an extended endoscopic triple $(G',s',{}^L\eta)$ of $G$ (defined in Section 1.3 of \cite{K}) such that $s'\in sZ_{\phi}^{\circ}$ and $\phi$ factors through ${}^L\eta$, let $\Pi_{\phi,s}(G')$ be the corresponding $L$-packet of $G'$ and let $\theta_{\Pi_{\phi,s}(G')}$ be the distribution character of that packet (which is a stable character on $G'(F)$). Then for $\alpha\in H^1(F,G)$, the character
$$\theta_{\Pi_\phi,\alpha,s}=\sum_{\pi\in \Pi_\phi(G_\alpha)} \tr(\chi_\pi(s))\theta_\pi$$
is the endoscopic transfer of $\iota(G_\alpha)\theta_{\Pi_{\phi,s}(G')}$.
\end{itemize}
\end{itemize}

In this paper, we will assume that the local Langlands conjecture holds for all the groups in Table \ref{fig:1}. To end this subsection, we will prove the existence of the lifting in our definition of the character of the component group. To be specific, let $(G,H)$ be one of the models in Table \ref{fig:1} and let $\phi:W_F'\rightarrow {}^L(G/Z_{G,H})$ be a tempered $L$-packet. Our goal is to prove the following lemma.

\begin{lemma}\label{lem extended endoscopic triple}
For $s\in S_\phi$, there exists an elliptic extended endoscopic triple $(G',s',{}^L\eta)$ of $G/Z_{G,H}$ such that $s'\in sZ_{\phi}^{\circ}$ and $\phi$ factors through ${}^L\eta$. 
\end{lemma}

\begin{proof}
We will only consider the $E_7$ case, the rest cases follow from a similar and easier argument. Let $G=E_7$ be the split adjoint reductive group of Type $E_7$ and its dual group $\hat{G}=E_{7,sc}(\BC)$ is simply connected (in particular $Z_{\hat{G}}(t)=\hat{G}_t$ for all $t\in \hat{G}_{ss}$). To prove the statement, it is enough to show that for any $s\in S_\phi$, there exists $s'\in sZ_{\phi}^\circ$ such that $s'$ is elliptic in $\hat{G}$. Here we say a semisimple element $t\in \hat{G}$ is elliptic if and only if the center of $\hat{G}_t$ is finite modulo the center of $\hat{G}$.

By induction we may assume that this is true for all the maximal Levi subgroups $M$ of $G$, i.e. for any tempered $L$-parameter $\phi_M$ of $M$ and for any $s_M\in S_{\phi_M}$, there exists $s_M'\in s_MZ_{\phi_M}^{\circ}$ that is elliptic in $\hat{M}$. 

Now we are ready to prove the statement. Let $s'$ be any element in $sZ_{\phi}^{\circ}$. If the center of $\hat{G}_{s'}$ is finite then we are done. If not, let $A$ be a split torus of the center of $\hat{G}_{s'}$. The centralizer of $A$ in $\hat{G}$ is a Levi subgroup of $\hat{G}$. Both $s'$ and $Im(\phi)$ belong  to the centralizer of $A$ in $\hat{G}$. In particular, there exists a maximal Levi subgroup $M$ of $G$ such that $s',Im(\phi)\subset \hat{M}$. Let $Z_{\phi,M}$ be the centralizer of $Im(\phi)$ in $\hat{M}$. We have $Z_{\phi,M}\subset Z_\phi$ and $Z_{\hat{M}}^{\circ}\subset Z_{\phi,M}^{\circ}\subset Z_{\phi}^{\circ}$.

By induction, we know that we may choose $s'$ so that it is elliptic in $\hat{M}$. If $\hat{M}$ is one of the following four maximal Levi subgroups of $\hat{G}$:
$$\SL_2(\BC)\times \SL_3(\BC)\times \SL_4(\BC)\times \GL_1(\BC)/(\BZ/12\BZ),$$
$$\SL_3(\BC)\times \SL_5(\BC)\times \GL_1(\BC)/(\BZ/15\BZ),$$
$$\SL_2(\BC)\times \SL_6(\BC)\times \GL_1(\BC)/(\BZ/6\BZ),\;\SL_7(\BC)\times \GL_1(\BC)/(\BZ/7\BZ),$$
the only elliptic elements in $\hat{M}$ are the center $Z_{\hat{M}}$. Moreover, each connected component of the center $Z_{\hat{M}}$ contains an element of $Z_{\hat{G}}$. Hence we may choose $s'$ so that it belongs to the center $Z_{\hat{G}}$ (in particular it is elliptic in $\hat{G}$).

If $\hat{M}=E_{6,sc}(\BC)\times \GL_1(\BC)/(\BZ/3\BZ)$, we have three cases. If $s'$ belongs to the center of $\hat{M}$, as in the previous case, we may choose $s'$ so that it belongs to the center $Z_{\hat{G}}$. If the centralizer of $s'$ in $\hat{M}$ is of Type $A_1\times A_5$, then it is easy to see that there exists an element in $s'Z_{\hat{M}_{s'}}^{\circ}\subset s'Z_{\phi,M}^{\circ}\subset s'Z_{\phi}^{\circ}$ that is elliptic in $\hat{G}$ and whose centralizer in $\hat{G}$ is isomorphic to $\Spin_{12}(\BC)\times \SL_2(\BC)/(\BZ/2\BZ)$. If the centralizer of $s'$ in $\hat{M}$ is of Type $A_2\times A_2\times A_2$, then it is easy to see that there exists an element in $s'Z_{\hat{M}_{s'}}^{\circ}\subset s'Z_{\phi,M}^{\circ}\subset s'Z_{\phi}^{\circ}$ that is elliptic in $\hat{G}$ and whose centralizer in $\hat{G}$ is isomorphic to $\SL_6(\BC)\times\SL_3(\BC)/(\BZ/3\BZ)$. This proves the case when $\hat{M}=E_{6,sc}(\BC)\times \GL_1(\BC)/(\BZ/3\BZ)$.

If $\hat{M}=\Spin_{12}(\BC)\times \GL_1(\BC)/(\BZ/2\BZ)$, we have three cases. If $s'$ belongs to the center of $\hat{M}$, as in the previous case, we may choose $s'$ so that it belongs to the center $Z_{\hat{G}}$. If the centralizer of $s'$ in $\hat{M}$ is of Type $D_4\times A_1\times A_1$, then it is easy to see that there exists an element in $s'Z_{\hat{M}_{s'}}^{\circ}\subset s'Z_{\phi,M}^{\circ}\subset s'Z_{\phi}^{\circ}$ that is elliptic in $\hat{G}$ and whose centralizer in $\hat{G}$ is isomorphic to $\Spin_{12}(\BC)\times \SL_2(\BC)/(\BZ/2\BZ)$. If the centralizer of $s'$ in $\hat{M}$ is of Type $A_3\times A_3$, then it is easy to see that there exists an element in $s'Z_{\hat{M}_{s'}}^{\circ}\subset s'Z_{\phi,M}^{\circ}\subset s'Z_{\phi}^{\circ}$ that is elliptic in $\hat{G}$ and whose centralizer in $\hat{G}$ is isomorphic to $\SL_8(\BC)/(\BZ/2\BZ)$. This proves the case when $\hat{M}=\Spin_{12}(\BC)\times \GL_1(\BC)/(\BZ/2\BZ)$.

If $\hat{M}=\Spin_{10}(\BC)\times \SL_2(\BC)\times \GL_1(\BC)/(\BZ/4\BZ)$, we have two cases. If $s'$ belongs to the center of $\hat{M}$, as in the previous case, we may choose $s'$ so that it belongs to the center $Z_{\hat{G}}$. If the centralizer of $s'$ in $\hat{M}$ is of Type $A_3\times A_1\times A_1\times A_1$, then it is easy to see that there exists an element in $s'Z_{\hat{M}_{s'}}^{\circ}\subset s'Z_{\phi,M}^{\circ}\subset s'Z_{\phi}^{\circ}$ that is elliptic in $\hat{G}$ and whose centralizer in $\hat{G}$ is isomorphic to $\SL_4(\BC)\times\SL_4(\BC)\times\SL_2(\BC)/(\BZ/4\BZ)$. This finishes the proof.
\end{proof}

\subsection{Transfer factors for $\GU_{2n},\;\GSp_{2n}$ and $\GSO_{2n}$}\label{section transfer factor}
In this subsection, we will discuss the transfer factors for $\GU_{2n},\;\GSp_{2n}$ and $\GSO_{2n}$, which are defined by the same formula as the classical groups case in \cite{Wal}. The only difference is that the similitude groups have fewer conjugacy classes compared with the classical groups. This will be used in later sections when we study the behavior of the geometric multiplicity under endoscopy.

We first discuss the semisimple conjugacy classes for these groups. For $\GSp_{2n}(F)$, the conjugacy classes of $\Sp_{2n}(F)$ are given by (Section 1.3 of \cite{Wal})
$$(F_i,F_{\pm i}, x_i,c_i)_{i\in I}$$
where
\begin{itemize}
\item $F_{\pm i}$ is a finite extension of degree $d_i$ with $\sum_{i\in I} d_i=n$, and $F_{i}$ is either a quadratic extension of $F_{\pm i}$ or $F_i=F_{\pm i}\oplus F_{\pm i}$. 
\item $x_i\in ker(N_{F_i/F_{\pm i}})$ and 
$$c_i\in (ker(\tr_{F_i/F_{\pm i}}) \cap F_{i}^{\times})/Im(N_{F_{i}/F_{\pm i}})$$ 
where $N_{F_i/F_{\pm i}}$ (resp. $\tr_{F_i/F_{\pm i}}$) is the norm map (resp. trace map).
\end{itemize}
The conjugacy classes of $\GSp_{2n}(F)$ are very similar to $\Sp_{2n}(F)$, 
and the only difference is that $(c_i)_{i\in I}$ needs to belong to a quotient of 
$$\Pi_{i\in I} (ker(\tr_{F_i/F_{\pm i}}) \cap F_{i}^{\times})/Im(N_{F_{i}/F_{\pm i}}).$$ 
To be specific, we say two elements $(c_i)_{i\in I}$ and $(c_i')_{i\in I}$ in 
$$\Pi_{i\in I} (ker(\tr_{F_i/F_{\pm i}}) \cap F_{i}^{\times})/Im(N_{F_{i}/F_{\pm i}})$$ 
are equivalent if they are differed by an element of $F^{\times}$, i.e. there exists $c\in F^{\times}$ such that $cc_i=c_i'$ for all $i\in I$. We use 
$$\Pi_{i\in I} (ker(\tr_{F_i/F_{\pm i}}) \cap F_{i}^{\times})/Im(N_{F_{i}/F_{\pm i}})/\sim$$ 
to denote the quotient under this equivalence. 
Then the conjugacy classes of $\GSp_{2n}(F)$ are given by 
$$(F_i,F_{\pm i}, x_i,c_i)_{i\in I}$$
where
\begin{itemize}
\item $F_{\pm i}$ is a finite extension of degree $d_i$ with $\sum_{i\in I} d_i=n$, and $F_{i}$ is either a quadratic extension of $F_{\pm i}$ or $F_i=F_{\pm i}\oplus F_{\pm i}$. 
\item $x_i\in F_{i}^{\times}$ and 
$$(c_i)_{i\in I}\in \Pi_{i\in I} (ker(\tr_{F_i/F_{\pm i}}) \cap F_{i}^{\times})/Im(N_{F_{i}/F_{\pm i}})/\sim$$ 
such that $N_{F_i/F_{\pm i}}(x_i)=N_{F_j/F_{\pm j}}(x_j)\in F^\times$ for all $i,j\in I$.
\end{itemize}

Next we consider $\GSO_{2n}(F)$. For our application, we only need to consider the split case. The conjugacy classes of $\SO_{2n}(F)$ are described by (Section 1.3 of \cite{Wal} and Section 1.4 of \cite{Wal3})
$$(F_i,F_{\pm i}, x_i,c_i)_{i\in I}$$
where
\begin{itemize}
\item $F_{\pm i}$ is a finite extension of degree $d_i$ with $\sum_{i\in I} d_i=n$, and $F_{i}$ is either a quadratic extension of $F_{\pm i}$ or $F_i=F_{\pm i}\oplus F_{\pm i}$.
\item $c_i\in F_{\pm i}^{\times}/Im(N_{F_{i}/F_{\pm i}})$ and $x_i\in ker(N_{F_i/F_{\pm i}})$.
\item The quadratic form associated to $(F_i,F_{\pm i},c_i)_{i\in I}$ (defined in Section 1.3 of \cite{Wal}) defines the split even special orthogonal group.  
\end{itemize}

Unlike the symplectic case, each $(F_i,F_{\pm i}, x_i,c_i)_{i\in I}$ determines two conjugacy classes in $\SO_{2n}(F)$ differed by the outer automorphism. 
The conjugacy classes of $\GSO_{2n}(F)$ are very similar to $\SO_{2n}(F)$, the only difference is that $(c_i)_{i\in I}$ needs to belong to the quotient (the equivalence relation is defined in the same way as in the symplectic case)
$$\Pi_{i\in I} F_{\pm i}^{\times}/Im(N_{F_{i}/F_{\pm i}})/\sim.$$  The conjugacy classes of $\GSO_{2n}(F)$ are described by 
$$(F_i,F_{\pm i}, x_i,c_i)_{i\in I}$$
where
\begin{itemize}
\item $F_{\pm i}$ is a finite extension of degree $d_i$ with $\sum_{i\in I} d_i=n$, and $F_{i}$ is either a quadratic extension of $F_{\pm i}$ or $F_i=F_{\pm i}\oplus F_{\pm i}$. 
\item $x_i\in F_{i}^{\times}$ and 
$$(c_i)_{i\in I}\in \Pi_{i\in I} F_{\pm i}^{\times}/Im(N_{F_{i}/F_{\pm i}})/\sim$$ 
such that $N_{F_i/F_{\pm i}}(x_i)=N_{F_j/F_{\pm j}}(x_j)\in F^\times$ for all $i,j\in I$.
\item The quadratic form associated to $(F_i,F_{\pm i},c_i)_{i\in I}$ defines the split even special orthogonal group.
\end{itemize}
Again each $(F_i,F_{\pm i}, x_i,c_i)_{i\in I}$ determines two conjugacy classes in $\GSO_{2n}(F)$ differed by the outer automorphism.

For $\GU_{2n}$, let $E/F$ be a quadratic field extension and we consider the group $\GU_{p,q}(F)$ with $p+q=2n$ (if $F$ is $p$-adic we require $p\in \{n,n\pm 1\}$). The conjugacy class of $U_{p,q}(F)$ is given by (Section 1.3 of \cite{Wal})
$$(F_i,F_{\pm i}, x_i,c_i)_{i\in I}$$
where
\begin{itemize}
\item $F_{\pm i}$ is a finite extension of degree $d_i$ with $\sum_{i\in I} d_i=n$, and $F_{i}=F_{\pm i}\otimes_F E$.
\item $c_i\in (ker(\tr_{F_i/F_{\pm i}})\cap F_{ i}^{\times})/Im(N_{F_{i}/F_{\pm i}})$ and $x_i\in ker(N_{F_i/F_{\pm i}})$.
\item The Hermitian form associated to $(F_i,F_{\pm i},c_i)_{i\in I}$ (defined in Section 1.3 of \cite{Wal}) defines the unitary group $U_{p,q}$. 
\end{itemize}
The conjugacy classes of $\GU_{p,q}(F)$ are very similar to $U_{p,q}(F)$, the only difference is that $(c_i)_{i\in I}$ needs to belong to a quotient of 
$$\Pi_{i\in I} (ker(\tr_{F_i/F_{\pm i}})\cap F_{ i}^{\times})/Im(N_{F_{i}/F_{\pm i}}).$$ 
To be specific, we say two elements $(c_i)_{i\in I}$ and $(c_i')_{i\in I}$ in 
$$\Pi_{i\in I} (ker(\tr_{F_i/F_{\pm i}})\cap F_{ i}^{\times})/Im(N_{F_{i}/F_{\pm i}})$$ 
are equivalent if they are differed by an element of $F^{\times}$. We use 
$$\Pi_{i\in I} (ker(\tr_{F_i/F_{\pm i}})\cap F_{ i}^{\times})/Im(N_{F_{i}/F_{\pm i}})/\sim$$ 
to denote the quotient under this equivalence. 
Then the conjugacy classes of $\GU_{p,q}(F)$ are given by 
$$(F_i,F_{\pm i}, x_i,c_i)_{i\in I}$$
where
\begin{itemize}
\item $F_{\pm i}$ is a finite extension of degree $d_i$ with $\sum_{i\in I} d_i=n$, and $F_{i}=F_{\pm i}\otimes_F E$. 
\item $x_i\in F_{i}^{\times}$ and 
$$(c_i)_{i\in I}\in \Pi_{i\in I} (ker(\tr_{F_i/F_{\pm i}})\cap F_{ i}^{\times})/Im(N_{F_{i}/F_{\pm i}})/\sim$$ 
such that $N_{F_i/F_{\pm i}}(x_i)=N_{F_j/F_{\pm j}}(x_j)\in F^\times$ for all $i,j\in I$.
\item The Hermitian form associated to $(F_i,F_{\pm i},c_i)_{i\in I}$ defines the unitary group $U_{p,q}$. 
\end{itemize}

\begin{rmk}
The stable semisimple conjugacy classes for all the cases above are given by $(F_i,F_{\pm i},x_i)$, i.e. the only difference between rational conjugacy classes and stable conjugacy classes is the extra $c_i$ for rational conjugacy classes.
\end{rmk}

Next we discuss the extended elliptic endoscopic triple $(G',s',{}^L\eta)$ for these groups. For $\GSp_{2n}$, the group $G'$ is of the form $G(\Sp_{2n_1}\times \SO_{2n_2})$ with $n=n_1+n_2$, $n_2\neq 1$, and $\SO_{2n_2}$ is the split even special orthogonal group. The projection of the element $s'\in \GSpin_{2n+1}(\BC)$ to $\SO_{2n+1}(\BC)$ is conjugate to the matrix $\diag(I_{2n_1+1},-I_{2n_2})$. 

For  $\GSO_{2n}$, the group $G'$ is of the form $G(\SO_{2n_1}\times \SO_{2n_2})$ with $n=n_1+n_2$, $n_1,n_2\neq 1$, and $\SO_{2n_1}$ (resp. $\SO_{2n_2}$) is the split even special orthogonal group. The projection of the element $s'\in \GSpin_{2n}(\BC)$ to $\SO_{2n}(\BC)$ is conjugate  to the matrix $\diag(I_{2n_1},-I_{2n_2})$. 

For $\GU_{p,q}$, the group $G'$ is of the form $G(U_{n_1,n_1}\times U_{n_2,n_2})$ with $n=n_1+n_2$ and the projection of the element $s'\in \GL_{2n}(\BC)\times \GL_1(\BC)$ to $\GL_{2n}(\BC)$ is conjugate to the matrix $\diag(I_{2n_1},-I_{2n_2})$. In all these cases ${}^L\eta$ is the natural embedding from ${}^LG'$ into ${}^LG$.

Finally we recall the definition of the transfer factor. Let $(G',s',{}^L\eta)$ be an elliptic extended endoscopic triple for 
$G=\GSp_{2n}\; \text{(resp.}\; \GSO_{2n},$ $ \GU_{p,q})$
with 
$$G'=G(\Sp_{2n_1}\times \SO_{2n_2})\; \text{(resp.}\; G(\SO_{2n_1}\times \SO_{2n_2}),\;G(U_{n_1,n_1}\times U_{n_2,n_2})).$$ 
Let $y=(y^+,y^-)$ be a stable conjugacy class of $G'(F)$ corresponding to
$$(I^+,(F_{ i}',F_{\pm i}',y_i)_{i\in I^+}),\;(I^-,(F_{i}',F_{\pm i}',y_i)_{i\in I^-})$$
where $y^+$ is a stable conjugacy class  of $\GSp_{2n_1}$ (resp. $\GSO_{2n_1},\;\GU_{n_1,n_1}$) correspondes to $(I^+,(F_{ i}',F_{\pm i}',y_i)_{i\in I^+})$ and $y^-$ is a stable conjugacy class  of $\GSO_{2n_2}$ (resp. $\GSO_{2n_2},\;\GU_{n_2,n_2}$) correspondes to $(I^-,(F_{ i}',F_{\pm i}',y_i)_{i\in I^-})$. Let $x$ be a conjugacy class of $G(F)$ corresponding to $(F_i,F_{\pm i}, x_i,c_i)_{i\in I}$. If $G$ is not the even special orthogonal group, we say $x$ and $y$ are associated to each other if 
$$I=I^+\cup I^-,\;(F_i,F_{\pm i},x_i)=(F_i',F_{\pm i}',y_i),\qquad\forall i\in I.$$
If $G$ is the even special orthogonal group, then we still need the above relation holds. In addition, there are 2 conjugacy classes of $G$ associated to $(F_i,F_{\pm i},x_i,c_i)_{i\in I}$ and there are four stable conjugacy classes of $G'$ associated to 
$$(I^+,(F_{ i}',F_{\pm i}',y_i)_{i\in I^+})\cup (I^-,(F_{i}',F_{\pm i}',y_i)_{i\in I^-})$$ 
Each conjugacy class of $G$ corresponds to two stable conjugacy classes of $G'$.

The transfer factor $\Delta(y,x)$ is nonzero only if $x$ and $y$ are associated to each other. If this is the case, then the transfer factor is given by (Section 1.10 of \cite{Wal})
$$\Delta_{G',G}(y,x)=\Pi_{i\in I^-}\eta_{F_i/F_{\pm i}}(c_i \cdot \ast)$$
where $\ast$ only depends on the stable conjugacy class of $x$. In this paper, we do not need the explicit formula of $\ast$. Instead, we only need to know whether the transfer factors are trivial or non-trivial. Hence we will not recall the definition of $\ast$ here and we refer the reader to Section 1.10 of \cite{Wal} for details.

\subsection{The epsilon factor}\label{sec epsilon factor}
Let $(G,H,\chi)$ be one of the models in Table \ref{fig:1}, $\phi:W_F'\rightarrow {}^L G/Z_{G,H}$ be a tempered Langlands parameter of $G$, $\Pi_\phi$ be the associated tempered $L$-packet, $Z_\phi$ be the centralizer of the parameter and $S_\phi$ be the component group. We fix an additive character $\psi$ of $F$ and we use $V$ to denote the underlying space of the representation $\rho_X$ (i.e. $\rho_X:{}^LG\rightarrow \GL(V)$). For $s\in S_\phi$, by Lemma \ref{lem extended endoscopic triple}, we can find an elliptic extended endoscopic triple $(G',s',{}^L\eta)$ such that $s'\in sZ_{\phi}^{\circ}$ and $\phi$ factors through ${}^L\eta$. For the model $(\GL_4\times \GL_2,\GL_2\times \GL_2)$, we require the lifting $s'$ to be of the form $\pm(I_4,I_2)$. Let $V_{s',-}$ be the $-1$ eigenspace of $V$ with respect to the operator $\rho_X(s')$. Since $s'$ commutes with $Im(\phi)$, the space $V_{s',-}$ is stable under $\rho_X(Im(\phi))$, this gives us a representation $\rho_{X,\phi,s'}$ of $W_F'$ on $V_{s',-}$, i.e. $\rho_{X,\phi,s'}:W_F'\rightarrow \GL(V_{s',-})$. If the model is not the two models in Table \ref{fig:1} related to unitary groups, we define
$$\omega_{\phi,H}(s)=\epsilon(\frac{1}{2},\rho_{X,\phi,s'},\psi).$$
For the two models related to unitary groups, we refer the reader to Sections \ref{sec:pre-GU4} and \ref{sec:pre-GU6} for its definition. Since $\rho_X$ is a sympletic representation, $\rho_{X,\phi,s'}$ is also symplectic and hence we have the following proposition.

\begin{prop}\label{prop epsilon}
The function $\omega_{\phi,H}$ is independent of the choice of the character $\psi$ and takes values in $\{\pm 1\}$. 
\end{prop}

For the rest of this paper, we will skip $\psi$ in the expression of the epsilon factor. In later sections we will also show that under certain assumption, the definition of $\omega_{\phi,H}(s)$ is independent of the choice of the elliptic extended endoscopic triple $(G',s',{}^L\eta)$ and $\omega_{\phi,H}$ is a quadratic character of $S_\phi$.

The goal of this subsection is to explicitly describe the representation $\rho_{X,\phi,s'}$ for each case. 

\subsubsection{The general linear group case} \label{sec:pre-GL}

We first consider the two models for general linear groups. For the model $(\GL_4\times \GL_2,\GL_2\times \GL_2)$, the dual group of $G/Z_{G,H}$ is 
\begin{eqnarray*}
(\GL_4(\BC)\times \GL_2(\BC))^1&:=&\{(g_1,g_2)\in \GL_4(\BC)\times \GL_2(\BC) \mid \\
&&\det(g_1)=\det(g_2)^{-1}\}.
\end{eqnarray*}
In this case, the component group $S_\phi$ is either the trivial group or $\BZ/2\BZ$. If $S_\phi=\BZ/2\BZ$ then $-(I_4,I_2)$ belongs to the non-neutral component component of $Z_\phi$. Recall that in this case we require the lifting of the element in $S_\phi$ to $Z_\phi$ to be of the form $\pm(I_4,I_2)$. We have $\omega_{\phi,H}(s)=1$ if $s'=(I_4,I_2)$  and $$\omega_{\phi,H}(s)=\epsilon(1/2,\Pi_\phi,\rho_X)\in \{\pm 1\}$$ 
if $s'=-(I_4,I_2)$.

For the model $(\GL_6,\GL_2\ltimes U)$, the dual group of $G/Z_{G,H}$ is $\SL_6(\BC)$. The component group $S_\phi$ is $\BZ/6\BZ$, $\BZ/3\BZ$, $\BZ/2\BZ$ or the trivial group. In this case, the lifting $s'\in Z_\phi$ is of the form $s'=aI_6$ with $a^6=1$. We have $\omega_{\phi,H}(s)=1$ if $a$ has order $1$ or $3$, and 
$$\omega_{\phi,H}(s)=\epsilon(1/2,\Pi_\phi,\rho_X)\in \{\pm 1\}$$ 
if $a$ has order $2$ or $6$.

\subsubsection{The model $(\GU_4\times \GU_2,(\GU_2\times \GU_2)^0)$} \label{sec:pre-GU4}

Next we consider the model $(\GU_4\times \GU_2,(\GU_2\times \GU_2)^0)$. In this case the dual group of $G/Z_{G,H}$ is $(\GL_4(\BC)\times \GL_2(\BC))^1 \times \GL_1(\BC)$ and the $L$-group is 
$$((\GL_4(\BC)\times \GL_2(\BC))^1 \times \GL_1(\BC))\rtimes \{1,\sigma\}$$ 
where $\sigma$ acts by the involution
$$(g,h,a)\in (\GL_4(\BC)\times \GL_2(\BC))^1 \times \GL_1(\BC)$$
$$\mapsto (J_{4}{}^tg^{-1}J_{4}^{-1},J_2{}^th^{-1}J_{2}^{-1},a\det(g)).$$

The representation $\rho_X=\rho_1\oplus \rho_2$ where $\rho_1$ is the tensor product of the exterior square representation of $\GL_4(\BC)$ with the standard representation of $\GL_2(\BC)$ and $\rho_2$ is the standard representation of $\GL_4(\BC)$ plus its dual. In particular, we have $\dim(\rho_1)=12$ and $\dim(\rho_2)=8$. We refer the reader to Section 6.4 of \cite{WZ2} and Section 3.1 of \cite{Z} for the $\sigma$-action on these spaces.

In this case the lifting $s'$ is of the form $s'=(s_1,s_2,a)$ with $s_1\in \GL_4(\BC)$ either equal to $\pm I_4$ or conjugate to $\diag(I_2,-I_2)$, $s_2=\pm I_2$, and $a\in \BC^\times$.

If $s_1=\pm I_4$, then
$$\epsilon(\frac{1}{2},\rho_{X,\phi,s'})=1\;\text{if}\;s_1=I_4,s_2=I_2,$$
$$\epsilon(\frac{1}{2},\rho_{X,\phi,s'})=\epsilon(\frac{1}{2},\Pi_\phi,\rho_1)\in \{\pm 1\}\;\text{if}\;s_1=I_4,s_2=-I_2,$$ $$\epsilon(\frac{1}{2},\rho_{X,\phi,s'})=\epsilon(\frac{1}{2},\Pi_\phi,\rho_2)\in \{\pm 1\}\;\text{if}\;s_1=-I_4,s_2=I_2,$$
$$\epsilon(\frac{1}{2},\rho_{X,\phi,s'})=\epsilon(\frac{1}{2},\Pi_\phi,\rho_X)\in \{\pm 1\}\;\text{if}\;s_1=-I_4,s_2=-I_2.$$

Let $\eta_{E/F}$ be the quadratic character of $F^\times$ and let $\chi_{\phi}$ be the central character of the $L$-packet of $\GU_4$ induced by $\Pi_\phi$. We define
$$\omega_{\phi,H}(s)=1\;\text{if}\;s_1=I_4,s_2=I_2,$$
$$\omega_{\phi,H}(s)=\eta_{E/F}(-1)\chi_\phi(-1)\epsilon(\frac{1}{2},\Pi_\phi,\rho_1)\in \{\pm 1\}\;\text{if}\;s_1=I_4,s_2=-I_2,$$ $$\omega_{\phi,H}(s)=\chi_\phi(-1)\epsilon(\frac{1}{2},\Pi_\phi,\rho_2)\in \{\pm 1\}\;\text{if}\;s_1=-I_4,s_2=I_2,$$
$$\omega_{\phi,H}(s)=\eta_{E/F}(-1)\epsilon(\frac{1}{2},\Pi_\phi,\rho_X)\in \{\pm 1\}\;\text{if}\;s_1=-I_4,s_2=-I_2.$$

If $s_1$ does not belong to the center, let $W$ be the 4 dimensional standard representation of $\GL_4(\BC)$ and we can decompose $W$ as
$$W=W_{s_1,+}\oplus W_{s_1,-}$$ 
where $W_{s_1,+}$ (resp. $W_{s_1,-}$) is the $+1$ (resp. $-1$) eigenspace of $s_1$ and $\dim(W_{s_1,+})=\dim(W_{s_1,-})=2$. We can also decompose $\rho_1\circ \phi$ as 
$$\rho_1\circ \phi=\rho_{1,s',\phi}\oplus \rho_{1,s',\phi}'$$ 
where the underlying vector space of $\rho_{1,s',\phi}$ is the tensor product representation of $\GL_2(\BC)\times \GL_2(\BC)\times \GL_2(\BC)$ (the first and second copy of $\GL_2(\BC)$ comes from the decomposition $W=W_{s_1,+}\oplus W_{s_1,-}$) and the underlying vector space of $\rho_{1,s',\phi}'$ is the tensor product of the standard representation of the third $\GL_2(\BC)$ copy with the direct product of the  exterior square representation of the first two $\GL_2(\BC)$ copy. 

We can also decompose $\rho_2\circ \phi$ as 
$$\rho_2\circ \phi=\rho_{2,s',\phi,+}\oplus \rho_{2,s',\phi,-}$$ 
where the underlying vector space of $\rho_{2,s',\phi,+}$ (resp. $\rho_{2,s',\phi,-}$) is the standard representation of $\GL_2(\BC)$ associated to $W_{s_1,+}$ (resp. $W_{s_1,-}$) plus its dual. All of these four representations are self-dual of symplectic type. We have
$$\epsilon(\frac{1}{2},\rho_{X,\phi,s'})=\epsilon(\frac{1}{2},\rho_{1,s',\phi}\oplus \rho_{2,s',\phi,-})\in \{\pm 1\}\;\text{if}\;s_2=I_2,$$
$$\epsilon(\frac{1}{2},\rho_{X,\phi,s'})=\epsilon(\frac{1}{2},\rho_{1,s',\phi}'\oplus \rho_{2,s',\phi,-})\in \{\pm 1\}\;\text{if}\;s_2=-I_2.$$
In this case, the parameter $\phi$ factors through the $L$-group of $G'=G(U_2\times U_2)\times \GU_2$ where the first two copies of $U_2$ correspond to the decomposition $W=W_{s_1,+}\oplus W_{s_1,-}$. We let $\Pi_\phi(G')$ be the associated $L$-packet and let $\chi_{\phi,s',1}$ (resp. $\chi_{\phi,s',2},\chi_{\phi,s',3}$) be the central character of the $L$-packet of $U_2$ obtained by restricting the $L$-packet $\Pi_\phi(G')$ to the first (resp. second, third) $U_2$ copy in $G'$. Then we have 
$$\chi_{\phi,s',1}\chi_{\phi,s',2}\chi_{\phi,s',3}=1,\;\chi_{\phi,s',1}(-1)\chi_{\phi,s',2}(-1)=\chi_{\phi}(-1).$$
We define
$$\omega_{\phi,H}(s)=\eta_{E/F}(-1)\chi_{\phi,s',2}(-1)\epsilon(\frac{1}{2},\rho_{1,s',\phi}\oplus \rho_{2,s',\phi,-})\in \{\pm 1\}\;\text{if}\;s_2=I_2,$$
$$\omega_{\phi,H}(s)=\chi_{\phi,s',1}(-1)\epsilon(\frac{1}{2},\rho_{1,s',\phi}'\oplus \rho_{2,s',\phi,-})\in \{\pm 1\}\;\text{if}\;s_2=-I_2.$$

\subsubsection{The model $(\GU_6,\GU_2\ltimes U)$}  \label{sec:pre-GU6}
Next we consider the model $(G,H)=(\GU_6,\GU_2\ltimes U)$. In this case the dual group of $G/Z_{G,H}$ is $\SL_6(\BC)$, and the $L$-group is $\SL_6(\BC)\rtimes \{1,\sigma\}$ where $\sigma$ acts by the involution 
$$g\in \SL_6(\BC)\mapsto J_6{}^tg^{-1}J_{6}^{-1}.$$ 
The representation $\rho_X$ is the exterior cube representation of $\SL_6(\BC)$ and we refer the reader to Section 3.1 of \cite{Z} for the $\sigma$-action on this space. Let $W$ be the 6-dimensional standard representation of $\SL_6(\BC)$. In this case the lifting $s'$ is conjugate to one of the following 4 matrices: $\diag(\pm I_2,\pm I_4)$.

If $s'=\pm I_6$, then $\epsilon(\frac{1}{2},\rho_{X,\phi,s'})=1$ if $s'=I_6$ and 
$$\epsilon(\frac{1}{2},\rho_{X,\phi,s'})=\epsilon(1/2,\Pi_\phi,\rho_X,\psi)\in \{\pm 1\}$$ 
if $s'=-I_6$. In this case we define $\omega_{\phi,H}(s)=1$ if $s'=I_6$ and 
$$\omega_{\phi,H}(s)=\eta_{E/F}(-1)\epsilon(1/2,\Pi_\phi,\rho_X,\psi)\in \{\pm 1\}$$ 
if $s'=-I_6$.

If $s'$ does not belong to the center of the dual group, we can decompose $W$ as
$$W=W_{s',+}\oplus W_{s',-}$$ 
where 
$W_{s',+}$ (resp. $W_{s',-}$) is the $+1$ (resp. $-1$) eigenspace of $s'$ and $\dim(W_{s',+})\in \{2,4\}$. In this case we can also decompose $\rho_X\circ \phi$ as 
$$\rho_X\circ \phi=\rho_{1,\phi,s'}\oplus \rho_{2,\phi,s'}$$ 
where the underlying vector space of $\rho_{1,\phi,s'}$ is the tensor product of the exterior square representation of $\GL_4(\BC)$ with the standard representation of $\GL_2(\BC)$, and the underlying vector space of $\rho_{2,\phi,s'}$ is the direct sum of the exterior cube representation of $\GL_4(\BC)$ and the tensor product of the standard representation of $\GL_4(\BC)$ with the exterior square representation of $\GL_2(\BC)$. If $\dim(W_{s',+})=2$ (resp. $\dim(W_{s',+})=4$), we have 
$$\epsilon(\frac{1}{2},\rho_{X,\phi,s'})=\epsilon(1/2,\rho_{2,\phi,s'})\in \{\pm 1\}$$ 

$$(\text{resp.}\; \epsilon(\frac{1}{2},\rho_{X,\phi,s'})=\epsilon(1/2,\rho_{1,\phi,s'})\in \{\pm 1\}).$$ 

In this case the parameter $\phi$ factors through the $L$-group of $G'=G(U_4\times U_2)$. Let $\Pi_\phi(G')$ be the associated $L$-packet and let $\chi_{\phi,s'}$ be the central character of the $L$-packet of $U_4$ induced by $\Pi_\phi(G')$. We define $$\omega_{\phi,H}(s)=\chi_{\phi,s'}(-1)\epsilon(1/2,\rho_{2,\phi,s'})\in \{\pm 1\}$$ 
$$\text{(resp.}\; \omega_{\phi,H}(s)=\eta_{E/F}(-1)\chi_{\phi,s'}(-1)\epsilon(1/2,\rho_{1,\phi,s'})\in \{\pm 1\})$$ 
if $\dim(W_{s',+})=2$ (resp. $\dim(W_{s',+})=4$).

\subsubsection{The model $(\GSp_6\times \GL_2,\GL_2\ltimes U)$}

Next we consider the model $(\GSp_6\times \GL_2,\GL_2\ltimes U)$. The dual group of $G/Z_{G,H}$ is 
\begin{eqnarray*}
(\GSpin_{7}(\BC)\times \GL_2(\BC))^1&=&\{(g_1,g_2)\in \GSpin_{7}(\BC)\times \GL_2(\BC) \mid \\
&&l(g_1)^{-1}=\det(g_2)\}.
\end{eqnarray*} 
The center of the dual group is isomorphic to $\GL_1(\BC)\times \BZ/2\BZ$. Let $W$ be the 7 dimensional standard representation of $\GSpin_{7}(\BC)$. In this case, up to multiplying the lifting by some element belonging to the neutral component of the center, the lifting $s'$ is of the form $s'=(s_1,I_2)$ with $s_1\in \Spin_7(\BC)$ such that $s_1$ induces a decomposition of $W$
$$W=W_{s_1,+}\oplus W_{s_1,-}$$ 
where 
$W_{s_1,+}$ (resp. $W_{s_1,-}$) is the $+1$ (resp. $-1$) eigenspace of $s_1$ and $\dim(W_{s_1,-})\in \{0,4,6\}$. 

If $s'$ belongs to the center of the dual group ($\iff \dim(W_{s_1,-})=0$), then $\omega_{\phi,H}(s)=1$ if $s'$ belongs to the neutral component of the center of of the dual group and 
$$\omega_{\phi,H}(s)=\epsilon(1/2,\Pi_\phi,\rho_X)\in \{\pm 1\}$$ 
if $s'$ does not belong to the neutral component of the center of of the dual group. 

If $s'$ does not belong to the center of the dual group and if the order of $s'$ is 4, then $\dim(W_{s_1,-})=6$ and it is easy to see that the space $V_{s',-}$ is zero and we have $\omega_{\phi,H}(s)=1$. In this case we can decompose $\rho_X\circ \phi$ as 
$$\rho_X\circ \phi=\rho_\phi\oplus \rho_{\phi}^\vee$$ 
where the underlying vector space of $\rho_\phi$ (resp. $\rho_{\phi}^\vee$) is the tensor product of the standard representation of $\GL_2(\BC)$ with a Half-Spin representation of the even Spin group associated to $W_{s_1,-}$. This implies that $\epsilon(\frac{1}{2},\Pi_\phi,\rho_X)=1$.

If $s'$ does not belong to the center of the dual group and the order of $s'$ is 2, then $\dim(W_{s_1,-})=4$.  Moreover, in this case, we can decompose the representation $\rho_X\circ \phi$ as 
$$\rho_X\circ \phi=\rho_{s',\phi,+}\oplus \rho_{s',\phi,-}$$ 
where the underlying vector space of $\rho_{s',\phi,+}$ (resp. $\rho_{s',\phi,-}$) is the tensor product of the Spin representation of the odd Spin group associated to $W_{s_1,+}$ with a Half-Spin representation of the even Spin group associated to $W_{s_1,-}$ and the standard representation of $\GL_2(\BC)$, and it is the $+1$ (resp. $-1$) eigenspace of $\rho_X(s')$. Both of them are self-dual of symplectic type and we have 
$$\omega_{\phi,H}(s)=\epsilon(\frac{1}{2},\rho_{s',\phi,-})\in \{\pm 1\}.$$

\begin{rmk}
For each decomposition $W=W_{+}\oplus W_{-}$ of the space $W$ with $\dim(W_{-})=4$, there are exactly two elements 
$$s'=(s_1,I_2),s''=(s_1',I_2)\in (\GSpin_{7}(\BC)\times \GL_2(\BC))^0$$ differed by the nontrivial element in center of $\Spin_{7}(\BC)$ such that $$W_+=W_{s_1,+}=W_{s_1',+},\;W_{-}=W_{s_1,-}=W_{s_1',-}.$$ 
The $+1$ (resp. $-1$) eigenspace of $\rho_X(s')$ is equal to the $-1$ (resp. $+1$) eigenspace of $\rho_X(s'')$. A similar version of this remark also applies to all the other models related to $\GSp$ and $\GSO$.
\end{rmk}

\subsubsection{The model $(\GSp_{10},\GL_2\ltimes U)$}

Next we consider the model $(G,H)$ $=(\GSp_{10},\GL_2\ltimes U)$. The dual group of $G/Z_{G,H}$ is $\Spin_{11}(\BC)$. Let $W$ be the 11-dimensional standard representation of $\Spin_{11}(\BC)$. In this case, the lifting $s'$ induces a decomposition of $W$ 
$$W=W_{s',+}\oplus W_{s',-}$$ 
where 
$W_{s',+}$ (resp. $W_{s',-}$) is the $+1$ (resp. $-1$) eigenspace of $s'$ and $\dim(W_{s',-})\in \{0,4,6,8,10\}$. 

If $s'$ belongs to the center of the dual group ($\iff \dim(W_{s',-})=0$), then $\omega_{\phi,H}(s)=1$ if $s'=1$ and 
$$\omega_{\phi,H}(s)=\epsilon(1/2,\Pi_\phi,\rho_X)\in \{\pm 1\}$$ 
if $s'\neq 1$. 

If $s'$ does not belong to the center of the dual group and if the order of $s'$ is 4, then $\dim(W_{s',-})\in \{6,10\}$ and it is easy to see that the space $V_{s',-}$ is zero and we have $\omega_{\phi,H}(s)=1$. In this case, we can decompose $\rho_X\circ \phi$ as 
$$\rho_X\circ \phi=\rho_\phi\oplus \rho_{\phi}^\vee$$ 
where the underlying vector space of $\rho_\phi$ (resp. $\rho_{\phi}^\vee$) is the tensor product of the Spin representation of the odd Spin group associated to $W_{s',+}$ with a Half-Spin representation of the even Spin group associated to $W_{s',-}$. This implies that $\epsilon(\frac{1}{2},\Pi_\phi,\rho_X)=1$.

If $s'$ does not belong to the center of the dual group and the order of $s$ is 2, then $\dim(W_{s',-})\in \{4,8\}$.  Moreover, in this case, we can decompose the representation $\rho_X\circ \phi$ as 
$$\rho_X\circ \phi=\rho_{s',\phi,+}\oplus \rho_{s',\phi,-}$$ 
where the underlying vector space of $\rho_{s',\phi,+}$ (resp. $\rho_{s',\phi,-}$) is the tensor product of the Spin representation of the odd Spin group associated to $W_{s',+}$ with a Half-Spin representation of the even Spin group associated to $W_{s',-}$ and it is the $+1$ (resp. $-1$) eigenspace of $\rho_X(s')$. Both of them are self-dual of symplectic type and we have $\omega_{\phi,H}(s)=\epsilon(\frac{1}{2},\rho_{s',\phi,-})\in \{\pm 1\}$.

\subsubsection{The model $(\GSp_6\times \GSp_4,G(\Sp_4\times \Sp_2))$}

Next we consider the model $(\GSp_6\times \GSp_4,G(\Sp_4\times \Sp_2))$. The dual group of $G/Z_{G,H}$ is 
\begin{eqnarray*}
(\GSpin_7(\BC)\times \GSpin_5(\BC))^1&=&\{(g_1,g_2)\in \GSpin_7(\BC)\times \GSpin_5(\BC) \mid \\
&&l(g_1)l(g_2)=1\}.
\end{eqnarray*}
The center of the dual group is isomorphic to $\GL_1(\BC)\times \BZ/2\BZ$. Let $V_1$ (resp. $V_2$) be the 7 (resp. 5) dimensional standard representation of $\GSpin_7(\BC)$ (resp. $\GSpin_5(\BC)$). In this case, up to multiplying the lifting by some element belonging to the neutral component of the center, the lifting $s'$ is of the form $s'=(s_1,s_2)$ with $(s_1,s_2)\in \Spin_7(\BC)\times \Spin_5(\BC)$. It induces the decomposition 
$$V_1=V_{1,s_1,+}\oplus V_{1,s_1,-},\; V_2=V_{2,s_2,+}\oplus V_{2,s_2,-}$$ where 
\begin{itemize}
\item $V_{1,s_1,+}$ (resp. $V_{1,s_1,-}$) is the $+1$ (resp. $-1$) eigenspace of $s_1$ and $\dim(V_{1,s_1,-})\in \{0,4,6\}$. 
\item $V_{2,s_2,+}$ (resp. $V_{2,s_2,-}$) is the $+1$ (resp. $-1$) eigenspace of $s_2$ and $\dim(V_{2,s_2,-})\in \{0,4\}$.
\end{itemize}

If $s'$ belongs to the center of the dual group ($\iff \dim(V_{1,s_1,-})=\dim(V_{2,s_2,-})=0$), then $\omega_{\phi,H}(s)=1$ if $s'$ belongs to the neutral component of the center of of the dual group and 
$$\omega_{\phi,H}(s)=\epsilon(1/2,\Pi_\phi,\rho_X)\in \{\pm 1\}$$ 
if $s'$ does not belong to the neutral component of the center of the dual group. 

Then we consider the case when $s$ does not belong to the center of the dual group. If the order of $s'$ is 4 ($\iff \dim(V_{1,s_1,-})=6$), it is easy to see that $\rho_X(s')$ does not have $-1$ eigenspace and we have $\omega_{\phi,H}(s)=1$. Moreover, we can decompose $\rho_X\circ \phi$ as $$\rho_X\circ \phi=\rho_\phi\oplus \rho_{\phi}^{\vee}$$ 
where the underlying vector space of $\rho_\phi$ (resp. $\rho_{\phi}^{\vee}$) is the tensor product of a Half Spin representation of the even Spin group associated to $V_{1,s_1,-}$ with the Spin representation of the Spin group associated to $V_2$. This implies that $\epsilon(\frac{1}{2},\Pi_\phi,\rho_X)=1$.

If the order of $s'$ is 2, there are three cases. The first case is when $\dim(V_{1,s_1,-})=4$ and $\dim(V_{2,s_2,-})=0$ (in particular, $s_2$ belongs to the center of $\GSpin_5(\BC)$). We can decompose the representation $\rho_X\circ \phi$ as 
$$\rho_X\circ \phi=\rho_{s',\phi,+}\oplus \rho_{s',\phi,-}$$  
where the underlying vector space of $\rho_{s',\phi,+}$ (resp. $\rho_{s',\phi,-}$) is the $+1$ (resp. $-1$) eigenspace of $\rho_X(s')$, and it is the tensor product of the Spin representation of $\GSpin_5(\BC)$ with the Spin representation of the odd Spin group associated to $V_{1,s_1,+}$ and a Half-Spin representation of the even Spin group associated to $V_{1,s_1,-}$. Both representations are self-dual of symplectic type. We have $\omega_{\phi,H}(s)=\epsilon(\frac{1}{2},\rho_{s',\phi,-})\in \{\pm 1\}$.

The second case is when $\dim(V_{1,s_1,-})=0$ and $\dim(V_{2,s_2,-})=4$ (in particular, $s_1$ belongs to the center of $\GSpin_7(\BC)$). We can decompose the representation $\rho_X\circ \phi$ as 
$$\rho_X\circ \phi=\rho_{s',\phi,+}\oplus \rho_{s',\phi,-}$$  
where the underlying vector space of $\rho_{s',\phi,+}$ (resp. $\rho_{s',\phi,-}$) is the $+1$ (resp. $-1$) eigenspace of $\rho_X(s')$, and it is the tensor product of the Spin representation of $\GSpin_7(\BC)$ with a Half-Spin representation of the even Spin group associated to $V_{2,s_2,-}$. Both representations are self-dual of symplectic type. We have $\omega_{\phi,H}(s)=\epsilon(\frac{1}{2},\rho_{s',\phi,-})\in \{\pm 1\}$.  

The last case is when $\dim(V_{1,s_1,-})=\dim(V_{2,s_2,-})=4$. We can decompose the representation $\rho_X\circ \phi$ as
$$\rho_X\circ \phi=\rho_{s',\phi,++}\oplus \rho_{s',\phi,+-}\oplus \rho_{s',\phi,-+}\oplus \rho_{s',\phi,--}$$
where the underlying vector space of each of the four representations is the tensor product of a Half-Spin representation of the even Spin group associated to $V_{2,s_2,-}$ with the Spin representation of the odd Spin group associated to $V_{1,s,+}$ and a Half-Spin representation of the even Spin group associated to $V_{1,s,-}$. Moreover $\rho_{s',\phi,++}\oplus \rho_{s',\phi,--}$ (resp. $\rho_{s',\phi,+-}\oplus \rho_{s',\phi,-+}$) is the $+1$ (resp. $-1$) eigenspace of $\rho_X(s')$. All of the four representations in the decomposition are self-dual of symplectic type. In this case, $\omega_{\phi,H}(s)$ is equal to
$$\epsilon(\frac{1}{2},\rho_{s',\phi,+-}\oplus \rho_{s',\phi,-+})\in \{\pm 1\}.$$

\subsubsection{The model $(\GSO_8\times \GL_2,\GL_2\ltimes U)$}
Next we consider the model $(\GSO_8\times \GL_2,\GL_2\ltimes U)$. The dual group of $G/Z_{G,H}$ is 
\begin{eqnarray*}
(\GSpin_{8}(\BC)\times \GL_2(\BC))^1&=&\{(g_1,g_2)\in \GSpin_{8}(\BC)\times \GL_2(\BC) \mid \\
&& l(g_1)^{-1}=\det(g_2)\}.
\end{eqnarray*}
The center of the dual group is isomorphic to 
$$\BC^{\times}\times \BZ/2\BZ\times \BZ/2\BZ$$ 
where $\BC^{\times}$ and the first copy of $\BZ/2\BZ$ act trivially under the Half-Spin representation $\rho_X$, and the second copy of $\BZ/2\BZ$ acts via the sign character.

Let $W$ be the 8-dimensional standard representation of $\GSpin_{8}(\BC)$. In this case, up to multiplying the lifting by an element belonging to the neutral component of the center, the lifting $s'$ is of the form $s'=(s_1,I_2)$ with $s_1\in \Spin_8(\BC)$ and $s_1$ induces a decomposition $$W=W_{s_1,+}\oplus W_{s_1,-}$$ 
where 
$W_{s_1,+}$ (resp. $W_{s_1,-}$) is the $+1$ (resp. $-1$) eigenspace of $s_1$ and $\dim(W_{s_1,-})\in \{0,4,8\}$.

If $s'$ belongs to the center of the dual group ($\iff \dim(W_{s_1,-})\in \{0,8\}$), then $\omega_{\phi,H}(s)=1$ if 
$$s'\in \BC^{\times}\times \BZ/2\BZ\times \{1\}$$ 
and $\omega_{\phi,H}(s)=\epsilon(1/2,\Pi_\phi,\rho_X)\in \{\pm 1\}$ if 
$$s'\in \BC^{\times}\times \BZ/2\BZ\times \{-1\}.$$

If $s'$ does not belong to the center of the dual group, then $\dim(W_{s_1,-})=4$.  We can decompose the representation $\rho_X\circ \phi$ as 
$$\rho_X\circ \phi=\rho_{s',\phi,+}\oplus \rho_{s',\phi,-}$$ 
where the underlying vector space of $\rho_{s',\phi,+}$ (resp. $\rho_{s',\phi,-}$) is the tensor product of a Half-Spin representation of the even Spin group associated to $W_{s_1,+}$ with a Half-Spin representation of the even Spin group associated to $W_{s_1,-}$ and the standard representation of $\GL_2(\BC)$, and it is the $+1$ (resp. $-1$) eigenspace of $\rho_X(s')$. Both representations are self-dual of symplectic type and we have $\omega_{\phi,H}(s)=\epsilon(\frac{1}{2},\rho_{s',\phi,-})\in \{\pm 1\}$.

\subsubsection{The model $(\GSO_{12},\GL_2\ltimes U)$}
Next we consider the model $(G,H)$ $=(\GSO_{12},\GL_2\ltimes U)$. The dual group of $G/Z_{G,H}$ is $\Spin_{12}(\BC)$. The center of the dual group is $(\BZ/2\BZ)^2$ and we will denote it by 
$$\{1,z,z',zz'\}$$ 
where $z$ belongs to the kernel of the map $\Spin_{12}(\BC)\rightarrow \SO_{12}(\BC)$ and $z'$ is the unique nontrivial element in the center that acts trivially on the Half-Spin representation $\rho_X$.

Let $W$ be the 12 dimensional standard representation of $\Spin_{12}(\BC)$. The lifting $s'$ induces a decomposition 
$$W=W_{s',+}\oplus W_{s',-}$$ 
where 
$W_{s',+}$ (resp. $W_{s',-}$) is the $+1$ (resp. $-1$) eigenspace of $s'$ and $\dim(W_{s',-})\in \{0,4,6,8,12\}$.

If $s'$ belongs to the center of the dual group ($\iff \dim(W_{s',-})\in \{0,12\}$), then $\omega_{\phi,H}(s)=1$ if $s'\in \{1,z'\}$ and $$\omega_{\phi,H}(s)=\epsilon(1/2,\Pi_\phi,\rho_X)\in \{\pm 1\}$$ 
if $s'\in \{z,zz'\}$. 

If $s'$ does not belong to the center of the dual group and if the order of $s'$ is 4, then $\dim(W_{s',-})=6$ and it is easy to see that the space $V_{s',-}$ is zero. This implies that $\omega_{\phi,H}(s)=1$. In this case we can decompose $\rho_X\circ \phi$ as 
$$\rho_X\circ \phi=\rho_\phi\oplus \rho_{\phi}^\vee$$ 
where the underlying vector space of $\rho_\phi$ (resp. $\rho_{\phi}^\vee$) is the tensor product of a Half-Spin representation of the even Spin group associated to $W_{s',+}$ with a Half-Spin representation of the even Spin group associated to $W_{s',-}$. This implies that $\epsilon(\frac{1}{2},\Pi_\phi,\rho_X)=1$.

If $s'$ does not belong to the center of the dual group and the order of $s'$ is 2, then $\dim(W_{s',-})\in \{4,8\}$.  Moreover, in this case, we can decompose the representation $\rho_X\circ \phi$ as 
$$\rho_X\circ \phi=\rho_{s',\phi,+}\oplus \rho_{s',\phi,-}$$ 
where the underlying vector space of $\rho_{s',\phi,+}$ (resp. $\rho_{s',\phi,-}$) is the tensor product of a Half-Spin representation of the even Spin group associated to $W_{s',+}$ with a Half-Spin representation of the even Spin group associated to $W_{s',-}$ and it is the $+1$ (resp. $-1$) eigenspace of $\rho_X(s')$. Both of them are self-dual of symplectic type and hence we have $\omega_{\phi,H}(s)=\epsilon(\frac{1}{2},\rho_{s',\phi,-})\in \{\pm 1\}$.

\subsubsection{The $E_7$ case}
The last case is the model $(E_7,\PGL_2\ltimes U)$. The dual group is the simply connected form $E_{7,sc}(\BC)$. We have five cases for the lifting $s'$. 

The first case is when $s'$ belongs to the center of the dual group (which is isomorphic to $\BZ/2\BZ$).  
In this case, we have $\omega_{\phi,H}(s)=1$ if $s'=1$ and 
$$\omega_{\phi,H}(s)=\epsilon(1/2,\Pi_\phi,\rho_X)\in \{\pm 1\}$$ 
if $s'\neq 1$. 

The second case is when $\hat{G}_{s'}\simeq \Spin_{12}(\BC)\times \SL_2(\BC)/(\BZ/2\BZ)$. Up to conjugation there are two such elements differed by the nontrivial element in the center, and both of them has order 2. In this case the restriction of the representation $\rho_X$ to the centralizer $\hat{G}_{s'}$ decomposes as $\rho_1\oplus \rho_2$ where $\rho_1$ is a Half-Spin representation of $\Spin_{12}(\BC)$ and $\rho_2$ is the tensor product of the standard representation of $\Spin_{12}(\BC)$ with the standard representation of $\SL_2(\BC)$. Both representations are self-dual of symplectic type. We can decompose the representation $\rho_X\circ \phi$ as 
$$\rho_X\circ \phi=\rho_{s',\phi,+}\oplus \rho_{s',\phi,-}$$ 
where the underlying vector space of $\rho_{s',\phi,+}$ (resp. $\rho_{s',\phi,-}$) is the $+1$ (resp. $-1$) eigenspace of $\rho_X(s')$. Moreover, we know that the underlying vector space of $\rho_{s',\phi,+}$ (resp. $\rho_{s',\phi,-}$) is either the space of $\rho_1$ (resp. $\rho_2$) or the space of $\rho_2$ (resp. $\rho_1$) depending on whether $\rho_X(s')$ acts identically on the space of $\rho_1$ or on the space of $\rho_2$. In both cases we have $\omega_{\phi,H}(s)=\epsilon(\frac{1}{2},\rho_{s',\phi,-})\in \{\pm 1\}$.

The third case is when $\hat{G}_{s'}\simeq \SL_6(\BC)\times\SL_3(\BC)/(\BZ/3\BZ)$. Up to conjugation there are two such elements differed by the nontrivial element in the center, one has order 6 and the other one has order 3. In this case the restriction of the representation $\rho_X$ to the centralizer $\hat{G}_{s'}$ decomposes as $$\rho_X=\rho_1\oplus \rho_2\oplus (\rho_2)^\vee$$ 
where $\rho_1$ is the exterior cube representation of $\SL_6(\BC)$ and $\rho_2$ is the tensor product of the standard representation of $\SL_6(\BC)$ with the standard representation of $\SL_3(\BC)$. We can decompose $\rho_X\circ \phi$ as 
$$\rho_X\circ \phi=\rho_{s',\phi}\oplus \rho_{s',\phi}'$$ 
where the underlying vector space of $\rho_{s',\phi}$ (resp. $\rho_{s',\phi}'$) is $\rho_1$ (resp. $\rho_2\oplus (\rho_2)^\vee$). In particular, we have
$$\epsilon(1/2,\Pi_\phi,\rho_X)=\epsilon(1/2,\rho_{s',\phi}).$$

If the order of $s'$ is 3, then $\rho_X(s')$ does not have $-1$ eigenspace and we have $\omega_{\phi,H}(s)=1$. If the order of $s'$ is 6, then the $-1$ eigenspace of $\rho_X(s')$ is the space of $\rho_1$ and we have
$\omega_{\phi,H}(s)=\epsilon(\frac{1}{2},\rho_{s',\phi})\in \{\pm 1\}$.

The fourth case is when $\hat{G}_{s'}\simeq \SL_4(\BC)\times\SL_4(\BC)\times \SL_2(\BC)/(\BZ/4\BZ)$. Such an element is unique up to conjugation. In this case, the restriction of the representation $\rho_X$ to the centralizer $\hat{G}_{s'}$ decomposes as $$\rho_X=\rho_1\oplus \rho_2\oplus \rho_3\oplus (\rho_3)^\vee$$ 
where $\rho_1$ (resp. $\rho_2$) is the tensor product of the exterior square representation of the first (resp. second) $\SL_4(\BC)$ copy with the standard representation of $\SL_2(\BC)$, and $\rho_3$ is the tensor product of the standard representation of the first $\SL_4(\BC)$ copy with the dual of the standard representation of the second $\SL_4(\BC)$ copy. We can decompose $\rho_X\circ \phi$ as 
$$\rho_X\circ \phi=\rho_{s',\phi,1}\oplus \rho_{s',\phi,2}\oplus \rho_{s',\phi}'$$ 
where the underlying vector space of $\rho_{s',\phi,1}$ (resp. $\rho_{s',\phi,2}$, $\rho_{s',\phi}'$) is the space of $\rho_1$ (resp. $\rho_2$, $\rho_3\oplus (\rho_3)^\vee$). In particular, we have
$$\epsilon(1/2,\Pi_\phi,\rho_X)=\epsilon(1/2,\rho_{s',\phi,1}\oplus \rho_{s',\phi,2}).$$

If $s'$ is of the form $(I_4,\pm iI_4,I_2)$ (resp. $(I_4,\pm iI_4,-I_2)$), the $-1$ eigenspace of $\rho_X(s')$ is the space of $\rho_2$ (resp. $\rho_1$), and we have $\omega_{\phi,H}(s)=\epsilon(\frac{1}{2},\rho_{s',\phi,2})\in \{\pm 1\}$ (resp. $\omega_{\phi,H}(s)=\epsilon(\frac{1}{2},\rho_{s',\phi,1})\in \{\pm 1\}$).

The last case is when $\hat{G}_{s'}=\SL_8(\BC)/\BZ_2$. Such an element is unique up to conjugation. In this case the restriction of the representation $\rho_X$ to the centralizer $\hat{G}_{s'}$ decomposes as $\rho_1\oplus (\rho_1)^\vee$ where $\rho_1$ is the exterior square representation of $\SL_8(\BC)$. It is easy to see that the $-1$ eigenspace of $\rho_X(s')$ is zero and hence $\omega_{\phi,H}(s)=1$. Moreover,  We can decompose $\rho_X\circ \phi$ as 
$$\rho_X\circ \phi=\rho_{s',\phi}\oplus \rho_{s',\phi}^{\vee}$$ 
where the underlying vector space of $\rho_{s',\phi}$ (resp. $\rho_{s',\phi}^{\vee}$) is the space of $\rho_1$ (resp. $(\rho_{1})^{\vee}$) and we have
$$\epsilon(1/2,\Pi_\phi,\rho_X)=1.$$

This finishes the description of $\rho_{X,\phi,s'}$ for all the models in Table \ref{fig:1}.

\section{The strategy of the proof}\label{sec strategy}
In this section, we will explain the strategy of the proof of Theorem \ref{main theorem}. Roughly speaking, the idea is to use the multiplicity formula $m(\pi)=m_{geom}(\pi)$ of the model $(G,H)$ to study the behaviors of the multiplicity under parabolic induction and endoscopic transfer (note that our assumption of the $L$-packet in Theorem \ref{main theorem} implies that the $L$-packet is either of endoscopic type or the parabolic induction of an $L$-packet of some Levi subgroup). Then we can reduce the problem to some models that are smaller than $(G,H)$, for which we can use the assumption in Theorem \ref{main theorem}.

We first explain the strategy for all the models $(G,H)$ in Table \ref{fig:1} except the model $(\GU_4\times \GU_2,(\GU_2\times \GU_2)^0)$. The model $(G,H)$ has a unique pure inner form $(G_D,H_D)$ corresponding to the unique quaternion algebra $D$ over $F$. We have the multiplicity formulas
$$m(\pi)=m_{geom}(\pi),\;m(\pi_D)=m_{geom}(\pi_D)$$
for these models. We will recall the definition of the geometric multiplicities $m_{geom}(\pi)$ and $m_{geom}(\pi_D)$ in later sections. Now let $\Pi_\phi=\Pi_\phi(G)\cup \Pi_{\phi}(G_D)$ be a tempered $L$-packet. We first consider the case when $\Pi_\phi$ is not discrete with $|\Pi_\phi(G)|=1$. In this case, one of the two statements is correct.
\begin{enumerate}
\item The $L$-packet $\Pi_\phi$ is the parabolic induction of an $L$-packet $\Pi_{\phi,M}$ of a proper parabolic subgroup $M$ of $G$.
\item There exists  a proper elliptic extended endoscopic triple $(G',s',{}^L\eta)$ of $G$ such that the Langlands parameter $\phi$ factors through the $L$-group of $G'$, i.e. the $L$-packet is of endoscopic type.
\end{enumerate}

If the $L$-packet is the parabolic induction of an $L$-packet of a proper Levi subgroup (we may assume that the Levi subgroup is a maximal Levi subgroup), by studying the behavior of the geometric multiplicity $m_{geom}(\pi)$ under parabolic induction, we only need to consider certain model associated to the Levi subgroup (denoted by $(M,M_H)$). If there is an analogue of the Levi subgroup $M$ in the pure inner form $G_D$ (denoted by $M_D$), then we can also study the behavior of $m_{geom}(\pi_D)$ under parabolic induction and hence we only need to consider certain model associated to the Levi subgroup (denoted by $(M_D,M_{H,D})$). For all the cases in Table 1, the models $(M,M_H)$ and $(M_D,M_{H,D})$ are either models that are smaller than $(G,H)$, or they have been considered in previous works. On the other hand, if  there is no analogue of the Levi subgroup $M$ in the pure inner form $G_D$, we will show in later sections that the model $(M,M_H)$ is just the Whittaker model which is well understood. Combining with our assumption in Theorem \ref{main theorem}, we can prove the weak conjecture for this $L$-packet.

For example, for the model $(G,H)=(\GSp_{6}\times \GSp_4,(\GSp_4\times \GSp_2)^0)$, if the Levi subgroup is 
$$\GSp_6\times (\GL_2\times \GL_1) \;\text{(resp.}\; (\GSp_2\times \GL_2\times \GL_1)\times \GSp_4),$$ 
we can reduce to the model 
$$(M,M_H)=(\GSp_6\times \GL_2,\GL_2\ltimes U)$$ 
$$\text{(resp.} \; (M,M_H)=(\GSp_4\times \GL_2\times \GL_2,(\GL_2\times \GL_2)^0))$$ 
which is smaller than $(G,H)$. All the other maximal Levi subgroups do not have an analogue in the pure inner form $G_D(F)=\GSp_3(D)\times \GSp_2(D)$ and the model $(M,M_H)$ is just the Whittaker model.

If the $L$-packet is of endoscopic type, by studying the behaviors of the geometric multiplicity $m_{geom}(\pi)$ (resp. $m_{geom}(\pi_D)$) under the endoscopic transfer from $G'$ to $G$ (resp. $G_D$), we only need to consider certain model associated to $G'$. For all the cases in Table 1, the model associated to $G'$ is either smaller than $(G,H)$ or  has been considered in previous works. Combining with our assumption in Theorem \ref{main theorem}, we can prove the weak conjecture for this $L$-packet. This also proves Theorem \ref{thm weak conjecture smaller models}.

For example, for the model $(G,H)=(\GSp_{6}\times \GSp_4,(\GSp_4\times \GSp_2)^0)$, if the endoscopic group is 
$$\GSp_6\times \GSO_4 \; \text{(resp.}\; G(\Sp_2\times \SO_4)\times \GSp_4),$$ 
then we can reduce to the model  
$$(\GSp_6\times \GL_2,\GL_2\ltimes U) \; \text{(resp.} \;(\GSp_4\times \GL_2\times \GL_2,(\GL_2\times \GL_2)^0))$$ 
which is smaller than $(G,H)$. If the endoscopic group is $\GSO_6\times \GSp_4$, we get the Whittaker model.

In both cases above, once we have proved the weak conjecture for the $L$-packet, we get a formula of the epsilon factor $\epsilon(1/2,\Pi_\phi,\rho_X)$ in terms of the Harish-Chandra character of the $L$-packet. We also have analogies of the formulas for smaller models by our assumption in Theorem \ref{main theorem}. Combining the formulas of epsilon factors with the formula of the geometric multiplicity under endoscopy and the definition of $\omega_{\phi,H}$, we can prove Theorem \ref{main theorem}. 

To be specific, by Remark \ref{distinguished is character}, the unique distinguished element in the packet corresponds to a character of the component group. We use $\omega_\phi$ to denote this character and we can view it as a character of the centralizer $Z_\phi$. For $s\in S_\phi$, let $(G',s',{}^L\eta)$ be an elliptic extended endoscopic triple such that $\phi$ factors through ${}^L\eta$ and $s'\in SZ_{\phi}^{\circ}$ (its existence was proved in Lemma \ref{lem extended endoscopic triple}). We just need to show that $\omega_{\phi,H}(s)=\omega_\phi(s')$. In fact, this will imply that our definition of $\omega_{\phi,H}$ is independent of the choice of the lifting $s'$ (i.e. it is a well-defined function on $S_\phi$) and it is a character of $S_\phi$. Combining with the fact that $Im(\omega_{\phi,H})\in \{\pm 1\}$, we know that $\omega_{\phi}=\omega_{\phi,H}$ is a quadratic character of $S_\phi$.

To prove the identity $\omega_{\phi,H}(s)=\omega_\phi(s')$, there are two cases. If $s'$ belongs to the center of the dual group, then the identity follows from the definition of $\omega_{\phi,H}$, the multiplicity formula, and the formula of the epsilon factor $\epsilon(1/2,\Pi_\phi,\rho_X)$. If $s'$ does not belong to the center of the dual group, the identity follows from the behavior of the geometric multiplicity under endoscopy, the definition of $\omega_{\phi,H}$, and the formula of some epsilon factor associated to $G'$ (obtained from the weak conjecture).

\begin{rmk}
In general it is not enough to compute the multiplicities by studying the endoscopic relations because one also needs to compute the summations of the multiplicities 
$$\sum_{\pi\in \Pi_\phi(G)} m(\pi), \; \;\sum_{\pi_D\in \Pi_\phi(G_D)}m(\pi_D).$$ 
But for all the cases in Table \ref{fig:1}, we have already proved the multiplicity one on the $L$-packet:
$$\sum_{\pi\in \Pi_\phi(G)} m(\pi)+\sum_{\pi_D\in \Pi_\phi(G_D)}m(\pi_D)=1.$$
The endoscopic relations will tell us the parity of the summations $$\sum_{\pi\in \Pi_\phi(G)} m(\pi), \; \;\sum_{\pi_D\in \Pi_\phi(G_D)}m(\pi_D),$$  
which allows us to compute these two summations (i.e. if the summation is odd then it must be 1, and if it is even then it must be 0).
\end{rmk}

The remaining case is when the packet is discrete with only one element. In this case, the packet is not of the endoscopic type or parabolic type. Moreover, the epsilon dichotomy conjecture and the weak conjecture are equivalent in this case. If we want to prove Conjecture \ref{main conj}, we still need to prove a formula of the epsilon factor $\epsilon(1/2,\Pi_\phi,\rho_X)$ in terms of the Harish-Chandra character of the $L$-packet. However, the epsilon factors $\epsilon(1/2,\Pi_\phi,\rho_X)$ in Table \ref{fig:1} are more complicated than the Gan--Gross--Prasad models case and we are not able to prove such a formula at this moment (one of the key obstacles is that the Langlands functoriality from $G/Z_{G,H}$ to $\GL_{\dim(\rho_X)}$ induced by $\rho_X$ is not known and this functoriality is not of twisted endoscopic type as in the Gan--Gross--Prasad model cases). This is why we exclude this case for Models 3--10 of Table \ref{fig:1} in our main theorem. 

For the remaining two models $(\GL_4\times \GL_2,\GL_2\times \GL_2)$ and $(\GU_4\times \GU_2,(\GU_2\times \GU_2)^0)$, we can prove the conjecture when the central character is trivial. Under this assumption, we can replace the groups $\GL_4,\GL_2,\GU_4,\GU_2$ by $\PGL_4,\PGL_2,\mathrm{PGU}_4,\mathrm{PGU}_2$. Then using some lower rank isomorphisms together with the multiplicity formulas, we can reduce to the cases of the Gan--Gross--Prasad models 
$$(\SO_6\times \SO_3,\SO_3\ltimes U),\;(U_4\times U_1,U_1\ltimes U)$$ which have been studied in \cite{Wal3} and \cite{Beu1}. 

In ongoing work, we are trying to prove this formula of epsilon factor and hence completely prove Conjecture \ref{main conj} by studying the multiplicity of certain models related to the Rankin-Selberg integrals.

Finally, we consider the model $(G,H)=(\GU_4\times \GU_2,(\GU_2\times \GU_2)^0)$, which is the most difficult  model in Table \ref{fig:1}. The reason is that it has more than one pure inner form.  
In the $p$-adic case (resp. real case), it has 3 (resp. 4) pure inner forms which will be denoted by $(G_i,H_i)$ for $1\leq i\leq 3$ (resp. $1\leq i\leq 4$). We refer the reader to Section \ref{sec GU} for details. If the packet is the parabolic induction of some packet of a Levi subgroup, we can prove the conjecture by the same argument as in all the other cases. On the other hand, if the packet is of endoscopic type, this will be more difficult than all the other cases.

To be specific, the group $G$ has a unique elliptic endoscopic group that is $G'=G(U_2\times U_2)\times \GU_2$. Like in all the other cases, we want to study the behavior of the geometric multiplicities $m_{geom}(\pi)$ and $m_{geom}(\pi_i)$ (where $1\leq i\leq 3$ in the $p$-adic case and $1\leq i\leq 4$ in the Archimedean case) under the endoscopic transfer. Here $\pi$ (resp. $\pi_i$) is a tempered representation of $G(F)$ (resp. $G_i(F)$). However, unlike all the other cases, these will not be related to the multiplicities of some models associated to $G'$. The reason is that some terms in the geometric multiplicities can not be eliminated under the endoscopic transfer (to be specific, the terms correspond to the regular elliptic conjugacy classes of $H(F)$ and $H_i(F)$), which is largely due to the fact that there are more than one pure inner forms. 
To solve this issue, instead of considering the behavior of each geometric multiplicity under the endoscopic transfer, we will consider some combinations of them. More specifically, in the p-adic case, we will consider
$$m_{geom}(\pi)-m_{geom}(\pi_1), \;m_{geom}(\pi_2)-m_{geom}(\pi_3).$$
In the real case, we will consider
$$m_{geom}(\pi)-m_{geom}(\pi_1)-m_{geom}(\pi_4), \;m_{geom}(\pi_2)-m_{geom}(\pi_3).$$
By considering these combinations, we can eliminate the terms corresponding to  the regular elliptic conjugacy classes of $H(F)$ and $H_i(F)$.
Then the endoscopic transfer of these combinations can be related to some models of $G'$. We refer the reader to Section \ref{sec GU} for details.

\section{The models $(\GL_4\times \GL_2,\GL_2\times \GL_2)$ and $(\GL_6,\GL_2\ltimes U)$}\label{Section GL}
\subsection{The models and the conjectures}
In this subsection, we recall the definitions of the models $(\GL_4\times \GL_2,\GL_2\times \GL_2)$ and $(\GL_6,\GL_2\ltimes U)$. We will also state Conjecture \ref{main conj} more explicitly for these models (note that for these two models Conjecture \ref{main conj} and Conjecture \ref{weak conjecture} are equivalent because the packets $\Pi_\phi(G)$ and $\Pi_\phi(G_D)$ contain at most one element).

We start with the model $(\GL_4\times \GL_2,\GL_2\times \GL_2)$. Let $G=\GL_4\times \GL_2$ and $H=\GL_2\times \GL_2$ embedded into $G$ via the map 
$$(a,b)\mapsto (\diag(a,b),b).$$ 
For the pure inner form, let $D/F$ be the quaternion algebra, $G_D=\GL_2(D)\times \GL_1(D)$ and $H_D=\GL_1(D)\times \GL_1(D)$ embedded into $G_D$ via the map 
$$(a,b)\mapsto (\diag(a,b),b).$$ 
In this case, $\rho_X=\rho_1\oplus \rho_2$ with $\rho_1=\wedge^2\otimes std_2$ and $\rho_2=std_4\oplus std_{4}^{\vee}$. Moreover, the local $L$-packet contains at most one element for each group and we have 
$$\epsilon(1/2,\pi,\rho_X)=\epsilon(1/2,\pi,\rho_1)\cdot \epsilon(1/2,\pi,\rho_2)=\epsilon(1/2,\pi,\rho_1)\cdot \omega_{\pi_1}(-1)$$
for all irreducible representations $\pi=\pi_1\otimes \pi_2$ of $G(F)$. Here $\omega_{\pi_i}$ is the central character of $\pi_i$ and we have $\omega_{\pi_1}(-1)=\omega_{\pi_2}(-1)$. As a result, Conjecture \ref{main conj} becomes the following conjecture.

\begin{conj}\label{conj GL(4)xGL(2)}
Let $\pi=\pi_1\otimes \pi_2$ be an irreducible tempered representation of $G(F)$ whose central character is trivial on $Z_{G,H}(F)$. Let $\pi_D$ be the Jacquet-Langlands correspondence of $\pi$ from $G(F)$ to $G_D(F)$ if it exists; otherwise let $\pi_D=0$. Then (note that $\chi=1$ in this case)
$$m(\pi)=1\iff \epsilon(1/2,\pi,\rho_X)=1,$$
$$m(\pi_D)=1\iff \epsilon(1/2,\pi,\rho_X)=-1.$$
\end{conj}

For the model $(\GL_6,\GL_2\ltimes U)$, let $G=\GL_6$, $H=H_0\ltimes U$ with 
$$H_0=\{diag(h,h,h) \mid h\in \GL_2\},$$
$$U=\{u(X,Y,Z)=\begin{pmatrix} I_2&X&Z\\0&I_2&Y\\0&0&I_2 \end{pmatrix} \mid X,Y,Z\in Mat_{2\times 2}\}.$$
The generic character on $U(F)$ is defined to be $$\xi(u(X,Y,Z))=\psi(\tr(X)+\tr(Y)).$$ 
This extends to a character $\chi$ of $H(F)$ that is trivial on $H_0(F)$. Similarly, we can define its pure inner form $(G_D,H_D,\chi_D)$ with  $G_D=\GL_3(D)$ and $H_{0,D}=\GL_1(D)$.

In this case, $\rho_X=\wedge^3$ and the local $L$-packet contains at most one element for each group. Conjecture \ref{main conj} becomes the following conjecture.

\begin{conj}\label{conj GL(6)}
Let $\pi$ be an irreducible tempered representation of $G(F)$ with trivial central character. Let $\pi_D$ be the Jacquet-Langlands correspondence of $\pi$ from $G(F)$ to $G_D(F)$ if it exists; otherwise let $\pi_D=0$. Then
$$m(\pi)=1\iff \epsilon(1/2,\pi,\rho_X)=1,$$
$$m(\pi_D)=1\iff \epsilon(1/2,\pi,\rho_X)=-1.$$
\end{conj}

To end this subsection, we discuss the multiplicity formula of these two models. For each of these two models, there is a canonical embedding from $\GL_2$ into $H_0$ and hence into $G$ (in the first case, this is the diagonally embedding), denoted by $\nu$. Similarly, there is a canonical embedding from $\GL_1(D)$ to $H_{0,D}(F)$ and hence into $G_D(F)$, denoted by $\nu_D$. 
Let $\theta$ (resp. $\theta_D$) be a quasi character of $G(F)$ (resp. $G_D(F)$), and define the geometric multiplicity to be
\begin{eqnarray*}
m_{geom}(\theta)&=&c_\theta(1)+\sum_{T\in \CT_{ell}(\GL_2)}|W(\GL_2,T)|^{-1}\\
&&\cdot \int_{T(F)/Z_{\GL_2}(F)} D^H(\nu(t)) c_\theta(\nu(t))dt,
\end{eqnarray*}
\begin{eqnarray*}
m_{geom}(\theta_D)&=&\sum_{T\in \CT_{ell}(\GL_1(D))}|W(\GL_1(D),T)|^{-1}\\ &&\cdot \int_{T(F)/Z_{\GL_1(D)}(F)} D^{H_D}(\nu_D(t)) c_{\theta_D}(\nu_D(t))dt.
\end{eqnarray*}
Recall that $\CT_{ell}(\GL_2)$ (resp. $\CT_{ell}(\GL_1(D))$) is a set of representatives of maximal elliptic tori of $\GL_2(F)$ (resp. $\GL_1(D)$). We have the multiplicity formulas (\cite{Wan15}, \cite{Wan16}, \cite{PWZ19})
$$m(\pi)=m_{geom}(\theta_\pi),\;m(\pi_D)=m_{geom}(\theta_{\pi_D})$$
for all tempered representations $\pi$ (resp. $\pi_D$) of $G(F)$ (resp. $G_D(F)$) with trivial central character.

As we explained in our previous papers (\cite{Wan15}, \cite{Wan16}, \cite{PWZ19}), the multiplicity formula implies the strong multiplicity one on the L-packet, i.e. 
$$m(\pi)+m(\pi_D)=m_{geom}(\theta_\pi)+m_{geom}(\theta_{\pi_D})=c_{\pi}(1)=1$$ 
where $\pi$ and $\pi_D$ are as in Conjecture \ref{conj GL(4)xGL(2)} and Conjecture \ref{conj GL(6)}.

By the multiplicity formula, we know that for these two models, Conjecture \ref{main conj} is equivalent to the following conjecture which expresses the epsilon factor in terms of the Harish-Chandra character.

\begin{conj}\label{epsilon formula GL}
Let $\pi$ be an irreducible tempered representation of $G(F)$ whose central character is trivial on $Z_{G,H}(F)$. Then
$$m_{geom}(\theta_\pi)=\frac{\epsilon(1/2,\pi,\rho_X)+1}{2}.$$
\end{conj}

\subsection{The proof of Theorem \ref{main theorem} and \ref{thm weak conjecture smaller models} for $(\GL_4\times \GL_2,\GL_2\times \GL_2)$ and $(\GL_6,\GL_2\ltimes U)$}
For the model $(\GL_6,\GL_2\ltimes U)$, when $F$ is Archimedean,  or when $\pi$ is not a discrete series of $\GL_6(F)$ or the parabolic induction of a discrete series of $\GL_4(F)\times \GL_2(F)$, Conjecture \ref{conj GL(6)} has already been proved in Theorem 1.4 of \cite{Wan16}. The same argument can be applied to the $(\GL_4\times \GL_2,\GL_2\times \GL_2)$ model case to prove Conjecture \ref{conj GL(4)xGL(2)} when $\pi$ is not a discrete series.

To prove the remaining parts of the theorem, we consider another model $(\GL_4\times \GL_2,\GL_2\ltimes U)$ defined in Appendix A.3 of \cite{Wan16}. This model is essentially the Gan--Gross--Prasad model $(\SO_6\times \SO_3,\SO_3\ltimes U)$ because of the isomorphisms $\PGL_4\simeq \PGSO_6$ and $\PGL_2\simeq \SO_3$. The geometric multiplicities for the models $(\GL_4\times \GL_2,\GL_2\times \GL_2)$ and $(\GL_4\times \GL_2,\GL_2\ltimes U)$ are the same (Appendix A.3 of \cite{Wan16} and Section 9.5 of \cite{PWZ19}). This implies that these two models have the same multiplicity for all tempered representations. 

We first consider the model $(\GL_4\times \GL_2,\GL_2\times \GL_2)$. The only case remaining is when $\pi$ has trivial central character. If this is the case, $\pi$ can be identified with a tempered representation of $\GSO_6\times \SO_3$ and we have $\epsilon(1/2,\pi,\rho_X)=\epsilon(1/2,\pi,\rho_1)$. By restriction we get a tempered $L$-packet of $\SO_6\times \SO_3$ (denoted by $\Pi=\Pi_1\otimes \Pi_2$). Moreover, the geometric multiplicity of the model $(\GL_4\times \GL_2,\GL_2\ltimes U)$ in Appendix A.3 of \cite{Wan16} is the same as the geometric multiplicity of the Gan--Gross--Prasad model $(\SO_6\times \SO_3,\SO_3\ltimes U)$ in Section 13.1 of \cite{Wal1}. Hence the multiplicity formulas imply that $m(\pi)$ is equal to the multiplicity of the $L$-packet $\Pi$ for the Gan--Gross--Prasad model $(\SO_6\times \SO_3,\SO_3\ltimes U)$. Combining with Theorem 4.3 of \cite{Wal3}, we know that 
$$m(\pi)=1\iff \epsilon(1/2,\Pi_1\times \Pi_2)=1,$$
$$m(\pi_D)=1\iff \epsilon(1/2,\Pi_1\times \Pi_2)=-1.$$
Then Conjecture \ref{conj GL(4)xGL(2)} follows from the fact that $\epsilon(s,\pi,\rho_1)=\epsilon(s,\Pi_1\times \Pi_2)$. This finishes the proof of Theorem \ref{main theorem} for the model $(\GL_4\times \GL_2,\GL_2\times \GL_2)$. 

For the model $(\GL_6,\GL_2\ltimes U)$, it remains to prove the case when $\pi$ is the parabolic induction of a discrete series $\pi'=\pi_1\otimes \pi_2$ of $\GL_4(F)\times \GL_2(F)$ under the assumption that Conjecture \ref{weak conjecture} holds for the model $(\GL_4\times \GL_2,\GL_2\times \GL_2)$ ($\iff$ Conjecture \ref{conj GL(4)xGL(2)} holds). In this case, Corollary 5.15 of \cite{Wan16} implies that the multiplicity of $\pi$ is equal to the multiplicity of $\pi'$ for the model $(\GL_4\times \GL_2,\GL_2\ltimes U)$, which is equal to the multiplicity of $\pi'$ for the model $(\GL_4\times \GL_2,\GL_2\times \GL_2)$ by the above discussion. Together with the assumption that Conjecture \ref{conj GL(4)xGL(2)} holds, we have
$$m(\pi)=1\iff \epsilon(1/2,\pi',\rho_1\oplus \rho_2)=1,$$
$$m(\pi_D)=1\iff \epsilon(1/2,\pi',\rho_1\oplus \rho_2)=-1.$$
Then Conjecture \ref{conj GL(4)xGL(2)} follows from the fact that 
\begin{align*}
\epsilon(1/2,\pi,\rho_X)=&\epsilon(1/2,\pi',\rho_1)\epsilon(1/2,\pi_{1}^\vee \otimes \omega_{\pi_1})\epsilon(1/2,\pi_1\otimes \omega_{\pi_1}^{-1})\\ 
=&\omega_{\pi_1}(-1)\cdot \epsilon(1/2,\pi',\rho_1)=\epsilon(1/2,\pi',\rho_1\oplus \rho_2),
\end{align*} 
where  $\omega_{\pi_1}$ is the central character of $\pi_1$. 
This finishes the proof of Theorem \ref{main theorem} for the model $(\GL_6,\GL_2\ltimes U)$. 

\begin{rmk}\label{GL(6) implies GL(4)}
The discussion above also shows that if Conjecture \ref{conj GL(6)} holds for all irreducible representations $\pi$ induced from $\GL_4(F)\times \GL_2(F)$, then Conjecture \ref{conj GL(4)xGL(2)} holds. This proves Theorem \ref{thm weak conjecture smaller models} for the model $(\GL_6,\GL_2\ltimes U)$.
\end{rmk}

We can also slightly generalize the conjecture for the model $(\GL_6,\GL_2\ltimes U)$. To be specific, let $\chi'$ be a character of $H(F)$ defined by 
$$\chi'(\diag(h,h,h)u(X,Y,Z))=\alpha(\det(h))\xi(u(X,Y,Z))$$ 
where $\alpha$ is any character of $F^\times$ (in particular if $\alpha=1$ then we recover the character $\chi$). Let $\pi$ be an irreducible representation of $G(F)$ with central character $\alpha^2$ and we can define the multiplicity $m(\pi,\chi')$. Similarly we can also define the character $\chi_D'$ of $H_D(F)$ and the multiplicity $m(\pi_D,\chi_D')$. Similar to the case when $\alpha=1$, we can define the geometric multiplicity
\begin{eqnarray*}
m_{geom}(\theta,\chi')&=&c_\theta(1)+\sum_{T\in \CT_{ell}(\GL_2)}|W(\GL_2,T)|^{-1} \\
&&\cdot \int_{T(F)/Z_{\GL_2}(F)} D^H(\nu(t)) c_\theta(\nu(t))\chi'(\nu(t))^{-1}dt,
\end{eqnarray*}
\begin{eqnarray*}
m_{geom}(\theta_D,\chi_D')&=&\sum_{T\in \CT_{ell}(\GL_1(D))}|W(\GL_1(D),T)|^{-1} \int_{T(F)/Z_{\GL_1(D)}(F)}\\
&& D^{H_D}(\nu_D(t)) c_{\theta_D}(\nu_D(t))\chi_D'(\nu_D(t))^{-1}dt.
\end{eqnarray*}
We have the multiplicity formulas (\cite{Wan15}, \cite{Wan16})
$$m(\pi,\chi')=m_{geom}(\theta_\pi,\chi'),\;m(\pi_D,\chi_D')=m_{geom}(\theta_{\pi_D},\chi_D')$$
for all tempered representations $\pi$ (resp. $\pi_D$) of $G(F)$ (resp. $G_D(F)$) with central character $\alpha^2$. Again the multiplicity formula implies the strong multiplicity one on the L-packet.

For the epsilon dichotomy conjecture, let $\phi_\alpha$ be the character of $W_F'$ corresponding to $\alpha$. Then the following conjecture is a generalization of Conjecture \ref{conj GL(6)}.

\begin{conj}\label{conj GL(6) general}
Let $\pi$ be an irreducible tempered representation of $G(F)$ with central character $\alpha^2$ and let $\phi_\pi$ be the Langlands parameter of $\pi$. Let $\pi_D$ be the Jacquet-Langlands correspondence of $\pi$ from $G(F)$ to $G_D(F)$ if it exists; otherwise let $\pi_D=0$. Then
$$m(\pi,\chi')=1\iff \epsilon(1/2,(\rho_X\circ\phi_\pi)\otimes \phi_{\alpha}^{-1})=1,$$
$$m(\pi_D,\chi_D')=1\iff \epsilon(1/2,(\rho_X\circ\phi_\pi)\otimes \phi_{\alpha}^{-1})=-1.$$
\end{conj}

By the same argument as in the case when $\alpha=1$, we can prove the above conjecture when $\pi$ is not a discrete series by assuming Conjecture \ref{conj GL(4)xGL(2)} holds.

\begin{rmk}
For the rest models in Table \ref{fig:1}, we can also put some nontrivial algebraic characters on the reductive part of $H(F)$ and we can still formulate the epsilon dichotomy conjecture. However, by twisting the representation of $G(F)$ by some suitable characters we can easily reduce the conjecture to the case when the character is trivial on the reductive part of $H(F)$. For example, for the model $(\GSO_{12},\GL_2\ltimes U)$, when we put a character $\alpha\circ \det$ on the $\GL_2(F)$-part, by twisting the representation of $\GSO_{12}(F)$ by the character $\alpha^{-1}\circ l$ we can reduce to the case when $\alpha$ is trivial. Hence for the rest models we will only consider the case when the character is trivial on the reductive part of $H(F)$.
\end{rmk}

\section{The models $(\GU_4\times \GU_2,(\GU_2\times \GU_2)^0)$ and $(\GU_6,\GU_2\ltimes U)$}\label{sec GU}
In this section we will consider the models 
$$(\GU_4\times \GU_2,(\GU_2\times \GU_2)^0),\;(\GU_6,\GU_2\ltimes U).$$ 
In Section \ref{sec GU(4) 1}, we define the model $(\GU_4\times \GU_2,(\GU_2\times \GU_2)^0)$ and study the behaviors of the geometric multiplicities under parabolic induction and endoscopic transfer. This is the most complicated case of this paper. In Section \ref{sec GU(4) preparation}, we will define and study the smaller model $(\GU_4,\GU_2\ltimes U)$, which will be used in our proof. In Section \ref{sec GU(4) 2}, we will prove Theorem \ref{main theorem} for $(\GU_4\times \GU_2,(\GU_2\times \GU_2)^0)$. Finally, in Section \ref{sec GU(6)}, we will prove Theorem \ref{main theorem} for $(\GU_6,\GU_2\ltimes U)$. Throughout this section, let $E=F(\sqrt{\eps})$ be a quadratic extension of $F$, $\eta_{E/F}$ be the quadratic character associated to $E$,  $N_{E/F}$ (resp. $\tr_{E/F}$) be the norm map (resp. trace map), and $x\rightarrow \bar{x}$ be the Galois action on $E$.

\subsection{The model $(\GU_4\times \GU_2,(\GU_2\times \GU_2)^0)$}\label{sec GU(4) 1}

For the model $(\GU_4\times \GU_2,(\GU_2\times \GU_2)^0)$,
let $G=\GU_{2,2}\times \GU_{1,1}$ and 
$$H=(\GU_{1,1}\times \GU_{1,1})^0=\{(h_1,h_2)\in \GU_{1,1}\times \GU_{1,1}\mid l(h_1)=l(h_2)\}.$$ 
We can embed  H into $G$ via the map
$$(h_1,h_2)\in H\mapsto (\begin{pmatrix}a&0&b\\0&h_1&0\\c&0&d \end{pmatrix} ,h_1)\in G,\;h_2=\begin{pmatrix}a&b\\c&d\end{pmatrix}.$$
The pure inner forms are 
$$(G_1,H_1)=(\GU_{2,2}\times \GU_{2,0},(\GU_{2,0}\times \GU_{0,2})^0),$$
$$(G_2,H_2)=(\GU_{3,1}\times \GU_{1,1},(\GU_{1,1}\times \GU_{2,0})^0),$$
$$(G_3,H_3)=(\GU_{3,1}\times \GU_{2,0},(\GU_{2,0}\times \GU_{1,1})^0),$$
$$(G_4,H_4)=(\GU_{4,0}\times \GU_{2,0},(\GU_{2,0}\times \GU_{2,0})^0)$$
where the last pair $(G_4,H_4)$ only appears in the Archimedean case.

Now we formulate the analogy of Conjecture \ref{weak conjecture} for this case. Let $\Pi_\phi=\Pi_\phi(G)\cup (\cup_i \Pi_\phi(G_i))$ be a tempered $L$-packet whose central character is trivial on $Z_{G,H}(F)$. Recall that we have defined the epsilon factors $\epsilon(\frac{1}{2},\Pi_\phi,\rho_1)$, $\epsilon(\frac{1}{2},\Pi_\phi,\rho_2)$ and the central character $\chi_\phi$ in Section \ref{sec:pre}.

\begin{conj}\label{weak conjecture GU(4)xGU(2)}
The unique distinguished element belongs to $\Pi_\phi(G)$ if and only if $$\chi_\phi(-1)\eta_{E/F}(-1)\epsilon(\frac{1}{2},\Pi_\phi,\rho_1)=\chi_\phi(-1)\epsilon(\frac{1}{2},\Pi_\phi,\rho_2)=1.$$ 
The unique distinguished element belongs to $\Pi_\phi(G_1)\cup \Pi_\phi(G_4)$ if and only if $$-\chi_\phi(-1)\eta_{E/F}(-1)\epsilon(\frac{1}{2},\Pi_\phi,\rho_1)=\chi_\phi(-1)\epsilon(\frac{1}{2},\Pi_\phi,\rho_2)=1.$$ 
The unique distinguished element belongs to $\Pi_\phi(G_2)$ if and only if 
$$\chi_\phi(-1)\eta_{E/F}(-1)\epsilon(\frac{1}{2},\Pi_\phi,\rho_1)=-\chi_\phi(-1)\epsilon(\frac{1}{2},\Pi_\phi,\rho_2)=1.$$ 
The unique distinguished element belongs to $\Pi_\phi(G_3)$ if and only if 
$$\chi_\phi(-1)\eta_{E/F}(-1)\epsilon(\frac{1}{2},\Pi_\phi,\rho_1)=\chi_\phi(-1)\epsilon(\frac{1}{2},\Pi_\phi,\rho_2)=-1.$$
\end{conj}

Next we recall the definition of the geometric multiplicities from Section 9 of \cite{WZ2}. Let $T_0$ be the unique element in $\CT_{ell}(\GU_{1,1})=\CT_{ell}(\GU_{2,0})$ that is isomorphic to 
$$E^{2,0}:=\{(a,b)\in E^\times \times E^{\times} \mid a\bar{a}=b\bar{b}\}\subset E^\times \times E^{\times}.$$ 
For $T\in \CT_{ell}(\GU_{1,1})=\CT_{ell}(\GU_{2,0})$ with $T\neq T_0$, let $$(T\times T)^{0}=\{(t_1,t_2)\in T\times T \mid l(t_1)=l(t_2)\}.$$ 
Up to conjugation, there is a unique embedding from $(T\times T)^{0}$ to $(\GU_{1,1}\times \GU_{1,1})^0$ (resp. $(\GU_{2,0}\times \GU_{0,2})^0$). Combining with the diagonal embedding from $T$ to $(T\times T)^{0}$, we get an embedding from $T$ to $G$ (resp. $G_1(F)$) that factors through $H$ (resp. $H_1$). 
We will denote this embedding by $\nu_T$ (resp. $\nu_{1,T}$).

For $T_0$, in the $p$-adic case up to conjugation there are two embeddings from $(T_0\times T_0)^{0}$ to $(\GU_{1,1}\times \GU_{1,1})^0$ (resp. $(\GU_{2,0}\times \GU_{0,2})^0$). Combining with the diagonal embedding from $T_0$ to $(T_0\times T_0)^{0}$, we get two embeddings  from $T_0$ to $G$ (resp. $G_1$). The centralizer of the image of one of the embeddings is quasi-split,
and we will denote this embedding by $\nu_{T_0}$ (resp. $\nu_{1,T_0}$), while the centralizer of the image of  the other embedding is not quasi-split. In the Archimedean case, we can define the embedding $\nu_{T_0}$ in the same way as in the $p$-adic case. On the other hand, up to conjugation there is only one embedding from $(T_0\times T_0)^{0}$ to $(\GU_{2,0}\times \GU_{0,2})^0$ and this defines the embedding $\nu_{1,T_0}$. Note that in this case the centralizer of the image of $\nu_{1,T_0}$ is still quasi-split.

Meanwhile, consider the following two subgroups of  $(T_0\times T_0)^0$ (we identify $T_0$ with $E^{2,0}:=\{(a,b)\in E^\times \times E^{\times} \mid a\bar{a}=b\bar{b}\}$):
$$T_0'=\{(1,1)\times (1,a)\in (T_0\times T_0)^{0} \mid a\in E^1\},$$
$$T_0''=\{(1,a)\times (1,b)\in (T_0\times T_0)^{0} \mid a,b\in E^1\}.$$
The two embeddings from $(T_0\times T_0)^{0}$ to $(\GU_{1,1}\times \GU_{1,1})^0$ (resp. $(\GU_{1,1}\times \GU_{2,0})^0$) induce two embeddings from $T_0'$ to $G$ (resp. $G_2$) that are conjugate  to each other. Let $\nu_{T_0'}$ (resp. $\nu_{2,T_0'}$) be one of the embeddings. Note that the projection of these embeddings to the first $\GU_{1,1}$ factor is the trivial map. The centralizers of the image of these embeddings are quasi-split.

On the other hand, the two embeddings from $(T_0\times T_0)^{0}$ to $(\GU_{1,1}\times \GU_{1,1})^0$ induce two embeddings from $T_0''$ to $G$. The centralizer of the image of one of the embeddings is quasi-split (we will denote this embedding by $\nu_{T_0''}$) and the centralizer of the image of the other embedding is not quasi-split. Similarly, we can also define the embeddings $\nu_{i,T_0''}$ from $T_0''$ to $G_i$ for $1\leq i\leq 3$.

Now we are ready to define the geometric multiplicity. Let $\theta$ (resp. $\theta_i$) be a quasi-character on $G(F)$ (resp. $G_i(F)$) with trivial central character. For $T\in \CT_{ell}(\GU_{1,1})=\CT_{ell}(\GU_{2,0})$, we use $T^{\ast}(F)$ to denote $T(F)/Z_{\GU_{1,1}}(F)=T(F)/Z_{\GU_{2,0}}(F)$. Define
\begin{eqnarray*}
m_{geom}(\theta)&=&c_\theta(1)+\sum_{T\in \CT_{ell}(H)} |W(H,T)|^{-1} \int_{T(F)/Z_{G,H}(F)} D^H(t)\theta(t)dt\\
&&+\frac{1}{2}\sum_{T\in \CT_{ell}(\GU_{1,1})}\int_{T^\ast(F)}D^H(\nu_T(t))c_\theta(\nu_T(t))dt\\
&&  +\int_{T_{0}'(F)}D^H(\nu_{T_0'}(t))c_\theta(\nu_{T_0'}(t))dt\\
&&+\int_{T_{0}''(F)}D^H(\nu_{T_0''}(t))c_\theta(\nu_{T_0''}(t))dt,
\end{eqnarray*}

\begin{eqnarray*}
m_{geom}(\theta_1)&=&\sum_{T\in \CT_{ell}(H_1)} |W(H_1,T)|^{-1} \int_{T(F)/Z_{G_1,H_1}(F)} D^{H_1}(t)\theta_1(t)dt\\
&&+\frac{1}{2}\sum_{T\in \CT_{ell}(\GU_{2,0})} \int_{T^\ast(F)} D^{H_1}(\nu_{1,T}(t))c_{\theta_1}(\nu_{1,T}(t))dt\\
&&+\int_{T_{0}''(F)}D^{H_1}(\nu_{1,T_0''}(t))c_{\theta_1}(\nu_{1,T_0''}(t))dt,
\end{eqnarray*}

\begin{eqnarray*}
m_{geom}(\theta_2)&=&\sum_{T\in \CT_{ell}(H_2)} |W(H_2,T)|^{-1} \int_{T(F)/Z_{G_2,H_2}(F)} D^{H_2}(t)\theta_{2}(t)dt\\
&&+\int_{T_{0}'(F)}D^{H_2}(\nu_{2,T_0'}(t))c_{\theta_2}(\nu_{2,T_0'}(t))dt\\
&&+\int_{T_{0}''(F)}D^{H_2}(\nu_{2,T_0''}(t))c_{\theta_2}(\nu_{2,T_0''}(t))dt,
\end{eqnarray*}
\begin{eqnarray*}
m_{geom}(\theta_3)&=&\sum_{T\in \CT_{ell}(H_3)} |W(H_3,T)|^{-1} \int_{T(F)/Z_{G_3,H_3}(F)} D^{H_3}(t)\theta_{3}(t)dt\\
&&+\int_{T_{0}''(F)}D^{H_3}(\nu_{3,T_0''}(t))c_{\theta_3}(\nu_{3,T_0''}(t))dt.
\end{eqnarray*}
If we are in the Archimedean case, we also define
$$m_{geom}(\theta_4)=\sum_{T\in \CT_{ell}(H_4)} |W(H_4,T)|^{-1} \int_{T(F)/Z_{G_4,H_4}(F)} D^{H_4}(t)\theta_{4}(t)dt.$$
In our previous paper \cite{WZ2}, we have proved the multiplicity formulas
$$m(\pi)=m_{geom}(\theta_\pi),\;m(\pi_i)=m_{geom}(\theta_{\pi_i})$$
for all tempered representations $\pi$ (resp. $\pi_i$) of $G(F)$ (resp. $G_i(F)$). 

\begin{rmk}
The above integrals need to be regularized, i.e. we replace $D^H(\cdot)$ (resp. $D^{H_i}(\cdot)$) by 
$$D^G(\cdot)^{1/2}(D^{H}(\cdot)^{-2} D^G(\cdot))^{s-1/2}$$ 
$$\text{(resp.}\; D^{G_i}(\cdot)^{1/2}(D^{H_i}(\cdot)^{-2} D^{G_i}(\cdot))^{s-1/2})$$
and take the limit $\lim_{s\rightarrow 0^+}$. Since this regularization does not affect our later computation, to simplify the notation, we will not include this regularization in the expression of the multiplicity formula.
\end{rmk}

Next we study the behavior of the geometric multiplicities under parabolic induction. Let $\theta$ be a quasi-character of $G(F)$, $M=M_1\times M_2$ be a maximal proper Levi subgroup of $G=\GU_{2,2}\times \GU_{1,1}$ and $\theta^M$ be a quasi-character of $M(F)$. If $M_2$ is a maximal quasi-split torus of $\GU_{1,1}$ and $M_1=\GU_{2,2}$, the embedding $\nu_{T_0'}$ from $T_0'$ to $G$ factors through $M$. We then define
$$m_{geom}(\theta^M)=c_{\theta^M}(1)+\int_{T_{0}'(F)}D^H(\nu_{T_0'}(t))c_{\theta^M}(\nu_{T_0'}(t))dt.$$
If $M_1$ is the Siegel Levi subgroup and $M_2=\GU_{1,1}$, up to conjugation we may assume that the embedding $\nu_T$ factors through $M$ for all $T\in \CT_{ell}(\GU_{1,1})$. Then we let
$$m_{geom}(\theta^M)=c_{\theta^M}(1)+\frac{1}{2}\sum_{T\in \CT_{ell}(\GU_{1,1})}\int_{T^\ast(F)}  D^{\GU_{1,1}}(t)c_{\theta^M}(\nu_T(t))dt.$$
The last case is when $M_1\simeq \GU_{1,1}\times Res_{E/F}\GL_1$ and $M_2=\GU_{1,1}$. Let $\iota$ (resp. $\iota'$) be an embedding from $E^{2,0}$ into $\GU_{1,1}(F)$ (resp. $\GU_{2,0}(F)$). We define (note that $M=Res_{E/F}\GL_1\times \GU_{1,1}\times \GU_{1,1}$)
\begin{eqnarray*}
m_{geom}(\theta^M)&=&c_{\theta^M}(1)+\int_{E^1} c_{\theta^M}(1,\iota(1,a),I_2)+c_{\theta^M}(1,aI_2,\iota(a,1)) da\\
&&+\int_{E^{2,0}} c_{\theta^M}(1,\iota(a,b),\iota(1,a))dadb.
\end{eqnarray*}

We also need to discuss the parabolic induction of $G_i$ for $1\leq i\leq 3$. Let $\theta_i$ be a quasi-character of $G_i(F)$, $M^i=M_{1}^{i}\times M_{2}^{i}$ be a maximal proper Levi subgroup of $G_i$ and $\theta^{M^i}$ be a quasi-character of $M^i(F)$. For $(G_1,H_1)$, we must have $M_{2}^{1}=\GU_{2,0}$. If $M_{1}^{1}$ is the Siegel Levi subgroup, up to conjugation we may assume that the embedding $\nu_{1,T}$ factors through $M^1$ for all $T\in \CT_{ell}(\GU_{2,0})$. Then we let
$$m_{geom}(\theta^{M^1})=\frac{1}{2}\sum_{T\in \CT_{ell}(\GU_{2,0})}\int_{T^\ast(F)}  D^{\GU_{2,0}}(t)c_{\theta^{M^1}}(\nu_{1,T}(t))dt.$$
If $M_{1}^{1}\simeq \GU_{1,1}\times Res_{E/F}\GL_1$, we define (note that $M^1=Res_{E/F}\GL_1\times \GU_{1,1}\times \GU_{2,0}$)
\begin{eqnarray*}
m_{geom}(\theta^{M^1})&=&\int_{E^1} c_{\theta^{M^1}}(1,aI_2,\iota'(a,1)) da\\
&&+\int_{E^{2,0}} c_{\theta^{M^1}}(1,\iota(a,b),\iota'(1,a))dadb.
\end{eqnarray*}

For the model $(G_2,H_2)$, if $M_{2}^{2}$ is a maximal quasi-split torus of $\GU_{1,1}$ and $M_{1}^{2}=\GU_{3,1}$, the embedding $\nu_{2,T_0'}$ from $T_0'$ to $G_2$ factors through $M^2$. We then define
$$m_{geom}(\theta^{M^2})=\int_{T_{0}'(F)}D^{H_2}(\nu_{2,T_0'}(t))c_{\theta^{M^2}}(\nu_{2,T_0'}(t))dt.$$
The remaining case is when $M_{2}^{2}=\GU_{1,1}$ and $M_{1}^{2}\simeq \GU_{2,0}\times Res_{E/F}\GL_1$. In this case, we define
\begin{eqnarray*}
m_{geom}(\theta^{M^2})&=&\int_{E^1} c_{\theta^{M^2}}(1,\iota'(1,a),I_2) da\\
&&+\int_{E^{2,0}} c_{\theta^{M^2}}(1,\iota'(a,b),\iota(1,a))dadb.
\end{eqnarray*}

For the model $(G_3,H_3)$, $G_3$ has a unique proper Levi subgroup  $M^3\simeq (Res_{E/F}\GL_1\times \GU_{2,0})\times \GU_{2,0}$. We define
$$m_{geom}(\theta^{M^3})=\int_{E^{2,0}} c_{\theta^{M^3}}(1,\iota'(a,b),\iota'(1,a))dadb.$$

The following proposition is a direct consequence of Proposition \ref{germ parabolic induction}.

\begin{prop}\label{parabolic GU(4)}
Let $\theta$ (resp. $\theta_i$) be a quasi-character on $G(F)$ (resp. $G_i(F)$). Assume that  $\theta$ (resp. $\theta_i$) is the parabolic induction of a quasi-character $\theta^M$ (resp. $\theta^{M^i}$) of a proper maximal Levi subgroup $M(F)$ of $G(F)$ (resp. $M^i(F)$ of $G_i(F)$). We have
$$m_{geom}(\theta)=m_{geom}(\theta^{M}),\;m_{geom}(\theta_i)=m_{geom}(\theta^{M^i}).$$
\end{prop}

Next we study the behavior of the geometric multiplicities under endoscopic transfer. The group $G$ has a unique proper elliptic endoscopic group $G'=G(U_{1,1}\times U_{1,1})\times \GU_{1,1}$. Let $(G',s',{}^L\eta)$ be a proper elliptic extended endoscopic triple with 
$$G'=G(U_{1,1}\times U_{1,1})\times \GU_{1,1},\; s'=(\diag(I_2,-I_2),1,-I_2,1)\in\hat{G},$$
$$\hat{G}=\GL_4(\BC)\times\GL_1(\BC)\times\GL_2(\BC)\times\GL_1(\BC),$$ 
and ${}^L\eta$ is the natural embedding.  This model is different from all the other cases considered in this paper, mainly because it has more than one pure inner form. For all the other cases, we only need to compute $m_{geom}(\theta)$ when $\theta$ is the endoscopic transfer of some stable character of $G'(F)$, which will give us the geometric multiplicity of a model associated to $G'$. 
But for the current model, if we only compute $m_{geom}(\theta)$ (resp. $m_{geom}(\theta_i)$ for $1\leq i\leq 4$), the expression we get does not correspond to the geometric multiplicity of a model associated to $G'$. This is due to the term in $m_{geom}(\theta)$ (resp. $m_{geom}(\theta_i)$) associated to $T\in \CT_{ell}(H)$ (resp. $T\in \CT_{ell}(H_i)$). Instead, we will consider the behavior of some combinations of $m_{geom}(\theta)$ and $m_{geom}(\theta_i)$. By doing this, we can eliminate the terms corresponding to $T\in \CT_{ell}(H)$ and $T\in \CT_{ell}(H_i)$.

Let $\theta'$ be a quasi-character on $G'(F)$. We fix an embedding $\iota$ from $E^{2,0}$ into $\GU_{1,1}$ as before. We define
\begin{eqnarray*}
m_{geom,1}(\theta')&=&c_{\theta'}(1)+\int_{E^1}c_{\theta'}((I_2,\iota(a,1),I_2))+c_{\theta'}((\iota(a,1),I_2,I_2))\\
&&+2c_{\theta'}(aI_2,I_2,\iota(a,1)) da+2\int_{E^1\times E^1}\\
&&c_{\theta'}((\iota(a,b),I_2,\iota(a,1)))+c_{\theta'}((I_2,\iota(a,b),\iota(a,1)))dadb,
\end{eqnarray*}

\begin{eqnarray*}
m_{geom,2}(\theta')&=&\int_{E^1}c_{\theta'}((\iota(a,1),I_2,I_2))-c_{\theta'}((I_2,\iota(a,1),I_2)) da\\
&&+2\int_{E^1\times E^1}c_{\theta'}((\iota(a,b),I_2,\iota(a,1)))\\
&&-c_{\theta'}((I_2,\iota(a,b),\iota(a,1)))dadb.
\end{eqnarray*}

\begin{prop}\label{prop GU(4)}
Let $\theta$ (resp. $\theta_i$) be a quasi-character on $G(F)$ (resp. $G_i(F)$). Assume that $\theta$ (resp. $\theta_i$) is the endoscopic transfer of a stable quasi-character $\theta'$ of $G'(F)$ . We have
$$m_{geom}(\theta)-m_{geom}(\theta_1)-m_{geom}(\theta_4)=m_{geom,1}(\theta'),$$
$$m_{geom}(\theta_2)-m_{geom}(\theta_3)=m_{geom,2}(\theta').$$
\end{prop}

\begin{proof}
We will only prove the first identity, the second identity follows from a similar argument. Recall that 
\begin{eqnarray*}
m_{geom}(\theta)&=&c_\theta(1)+\sum_{T\in \CT_{ell}(H)} |W(H,T)|^{-1} \int_{T(F)/Z_{G,H}(F)} D^H(t)\theta(t)dt\\
&&+\frac{1}{2}\sum_{T\in \CT_{ell}(\GU_{1,1})}\int_{T^\ast(F)}  D^H(\nu_T(t))c_\theta(\nu_T(t))dt\\
&&+\int_{T_{0}'(F)}D^H(\nu_{T_0'}(t))c_\theta(\nu_{T_0'}(t))dt\\
&&+\int_{T_{0}''(F)}D^H(\nu_{T_0''}(t))c_\theta(\nu_{T_0''}(t))dt,
\end{eqnarray*}
\begin{eqnarray*}
m_{geom}(\theta_1)&=&\sum_{T\in \CT_{ell}(H_1)} |W(H_1,T)|^{-1} \int_{T(F)/Z_{G_1,H_1}(F)} D^{H_1}(t)\theta_1(t)dt\\
&&+\frac{1}{2}\sum_{T\in \CT_{ell}(\GU_{2,0})} \int_{T^\ast(F)} D^{H_1}(\nu_{1,T}(t))c_{\theta_1}(\nu_{1,T}(t))dt\\
&&+\int_{T_{0}''(F)}D^{H_1}(\nu_{1,T_0''}(t))c_{\theta_1}(\nu_{1,T_0''}(t))dt,
\end{eqnarray*}
$$m_{geom}(\theta_4)=\sum_{T\in \CT_{ell}(H_4)} |W(H_4,T)|^{-1} \int_{T(F)/Z_{G_4,H_4}(F)} D^{H_4}(t)\theta_{4}(t)dt.$$

First it is easy to see from Proposition \ref{regular germs} and the definition of transfer factors that $c_\theta(1)=c_{\theta'}(1)$. Next we study the term corresponds to $T_0'$. For $t=(1,1)\times (1,a)\in T_0'(F)$, under the notation of Section \ref{section transfer factor}, we know that $D^H(\nu_{T_0'}(t))c_\theta(\nu_{T_0'}(t))$ is equal to the limit of the value of $(D^G)^{1/2}\cdot\theta$ at the conjugacy class (note that $c_i$ is unique in this case and hence we will ignore it)
$$((E\oplus E,E,(\lambda_1,\lambda_{1}^{-1}))\cup (E,F,a)\cup (E,F,1))\times ((E\oplus E,E,(\lambda_2,\lambda_{2}^{-1}))$$
times $\frac{D^{\GU_{1,1}}(\iota(a,1))^{-1/2}}{4}$ as $\lambda_i\in F^{\times}\rightarrow 1$. The value of $(D^G)^{1/2}\cdot\theta$ at the conjugacy class 
$$((E\oplus E,E,(\lambda_1,\lambda_{1}^{-1}))\cup (E,F,a)\cup (E,F,1))\times ((E\oplus E,E,(\lambda_2,\lambda_{2}^{-1}))$$
is equal to the summation of the values of $(D^{G'})^{1/2}\theta'$ at $$(\diag(\lambda_1,\lambda_{1}^{-1}),\iota(a,1),\diag(\lambda_2,\lambda_{2}^{-1}))$$
and
$$(\iota(a,1),\diag(\lambda_1,\lambda_{1}^{-1}),\diag(\lambda_2,\lambda_{2}^{-1})).$$ 
The transfer factor is equal to 1 on $(\iota(a,1),\diag(\lambda_1,\lambda_{1}^{-1}),\diag(\lambda_2,\lambda_{2}^{-1}))$ because the quadratic character $\eta_{F\oplus F/F}$ is trivial. Combining with the argument in Section 1.11 of \cite{Wal}, we know that the transfer factor is also equal to 1 on $(\diag(\lambda_1,\lambda_{1}^{-1}),\iota(a,1),\diag(\lambda_2,\lambda_{2}^{-1}))$.
If we take the limit as $\lambda_i\in F^{\times}\rightarrow 1$, we get 
$$4 D^{G'}((I_2,\iota(a,1),I_2))^{1/2}c_{\theta'}((I_2,\iota(a,1),I_2))$$
$$+4 D^{G'}((\iota(a,1),I_2,I_2))^{1/2}c_{\theta'}((\iota(a,1),I_2,I_2)).$$
Combining with the equations
$$D^{G'}((I_2,\iota(a,1),I_2))^{1/2}=D^{G'}((\iota(a,1),I_2,I_2))^{1/2}=D^{\GU_{1,1}}(\iota(a,1))^{1/2},$$
we know that $D^H(\nu_{T_0'}(t))c_\theta(\nu_{T_0'}(t))$ is equal to 
$$c_{\theta'}((I_2,\iota(a,1),I_2))+c_{\theta'}((\iota(a,1),I_2,I_2)).$$
This implies that
\begin{eqnarray*}
&&\int_{T_{0}'(F)}D^H(\nu_{T_0'}(t))c_\theta(\nu_{T_0'}(t))dt\\
&=&\int_{E^1}c_{\theta'}((I_2,\iota(a,1),I_2))+c_{\theta'}((\iota(a,1),I_2,I_2)) da.
\end{eqnarray*}

Next we study the terms correspond to $T_{0}''$. For $t=(1,a)\times (1,b)\in T_0'(F)$, we know that $D^H(\nu_{T_0''}(t))c_\theta(\nu_{T_0''}(t))$ is equal to the limit of the value of $(D^G)^{1/2}\cdot\theta$ at the conjugacy class (note that $c_i$ is unique in this case and hence we will ignore it)
$$((E\oplus E,E,(\lambda,\lambda^{-1}))\cup (E,F,a)\cup (E,F,b))\times ((E,F,a)\cup (E,F,1))$$
times 
$$\frac{D^{\GU_{1,1}}(\iota(a,1))^{-1/2}\cdot D^{\GU_{1,1}}(\iota(a,b))^{-1/2}}{2}$$ 
as $\lambda\in F^{\times}\rightarrow 1$. The value of $(D^G)^{1/2}\cdot\theta$ at the conjugacy class 
$$((F\oplus F,F,(\lambda,\lambda^{-1}))\cup (E,F,a)\cup (E,F,b))\times ((E,F,a)\cup (E,F,1))$$
is equal to the summation of the values of $(D^{G'})^{1/2}\theta'$ at $$(\diag(\lambda,\lambda^{-1}),\iota(a,b),\iota(a,1))$$
and
$$(\iota(a,b),\diag(\lambda,\lambda^{-1}),\iota(a,1))).$$ 
The transfer factor is equal to 1 on $(\iota(a,b),\diag(\lambda,\lambda^{-1}),\iota(a,1)))$ because the quadratic character $\eta_{F\oplus F/F}$ is trivial. Combining with the argument in Section 1.11 of \cite{Wal}, we know that the transfer factor is also equal to 1 on $(\diag(\lambda,\lambda^{-1}),\iota(a,b),\iota(a,1))$.
If we take the limit as $\lambda\in F^{\times}\rightarrow 1$, we get 
$$2 D^{G'}((I_2,\iota(a,b),\iota(a,1)))^{1/2}c_{\theta'}((I_2,\iota(a,b),\iota(a,1)))$$
$$+2 D^{G'}((\iota(a,b),I_2,\iota(a,1)))^{1/2}c_{\theta'}((\iota(a,b),I_2,\iota(a,1))).$$
Combining with the equations
$$D^{G'}((\iota(a,b),I_2,\iota(a,1)))^{1/2}=D^{G'}((I_2,\iota(a,b),\iota(a,1)))^{1/2}$$
$$=D^{\GU_{1,1}}(\iota(a,1))^{1/2}D^{\GU_{1,1}}(\iota(a,b))^{1/2},$$
we know that $D^H(\nu_{T_0''}(t))c_\theta(\nu_{T_0''}(t))$ is equal to 
$$c_{\theta'}((I_2,\iota(a,b),\iota(a,1)))+c_{\theta'}((\iota(a,b),I_2,\iota(a,1))).$$
Similarly, we can show that $D^{H_1}(\nu_{1,T_0''}(t))c_{\theta_1} (\nu_{1,T_0''}(t))$ is equal to (the negative sign comes from the Kottwitz sign between $\GU_{1,1}$ and $\GU_{2,0}$)
$$-c_{\theta'}((I_2,\iota(a,b),\iota(a,1)))-c_{\theta'}((\iota(a,b),I_2,\iota(a,1))).$$
This implies that
$$\int_{T_{0}''(F)}D^H(\nu_{T_0''}(t))c_\theta(\nu_{T_0''}(t))dt-\int_{T_{0}''(F)}D^{H_1}(\nu_{1,T_0''}(t))c_{\theta_1}(\nu_{1,T_0''}(t))dt$$
is equal to
$$2\int_{E^1\times E^1}c_{\theta'}((\iota(a,b),I_2,\iota(a,1)))+c_{\theta'}((I_2,\iota(a,b),\iota(a,1)))dadb.$$

Next we study the terms correspond to $T\in \CT_{ell}(\GU_{1,1})$ and $T\in \CT_{ell}(\GU_{2,0})$. We know that there is a natural bijection $T\leftrightarrow F_T$ between $\CT_{ell}(\GU_{1,1})$ and the set of quadratic extensions of $F$. If $F_T\neq E$, then $E_T=F_T\otimes_F E$ is a quadratic extension of $E$. For $t\in T(F)$, we can identify it with an element in $E_{T}^{\times}$ (by abusing of notation we still denote it by $t$), and $D^H(\nu_{T}(t))c_\theta(\nu_{T}(t))$ is equal to the limit of the value of $(D^G)^{1/2}\cdot\theta$ at the conjugacy class ($\bar{t}$ is the conjugation of $t$ by the nontrivial element in $Gal(E_T/F_T)$)
$$(E_T\oplus E_T,E_T,\lambda t,\lambda^{-1}\bar{t})\times (E_T,F_T,t)$$
times $\frac{D^{\GU_{1,1}}(t)}{2}$ as $\lambda\rightarrow 1$. Note that $c_i$ is unique in this case and hence we will ignore it. It is easy to see that the above conjugacy class does not correspond to a conjugacy class of $G'(F)$. Hence we know that $D^H(\nu_{T}(t))c_\theta(\nu_{T}(t))=0$ and 
$$\int_{T^{\ast}(F)} D^H(\nu_T(t))c_\theta(\nu_T(t))dt=0.$$

If $F_T=E$, $T^{\ast}(F)$ is isomorphic to $E^1$. For $t\in T^{\ast}(F)$ corresponding to $a\in E^1$, $D^H(\nu_{T}(t))c_\theta(\nu_{T}(t))$ is equal to the limit of the value of $\frac{1}{4}(D^G)^{1/2}\cdot\theta$ at the conjugacy class (note that $c_i$ is unique in this case and hence we will ignore it)
$$((E\oplus E,E,\lambda a,\lambda^{-1}a^{-1})\cup (E\oplus E,E,\lambda ,\lambda^{-1})) \times ((E,F,a)\cup (E,F,1))$$
times $\frac{D^{\GU_{1,1}}(\iota(a,1))^{-1/2}}{4}$ as $\lambda\rightarrow 1$. The value of $(D^G)^{1/2}\cdot\theta$ at the conjugacy class 
$$((E\oplus E,E,\lambda a,\lambda^{-1}a^{-1})\cup (E\oplus E,E,\lambda ,\lambda^{-1})) \times ((E,F,a)\cup (E,F,1))$$
is equal to the summation of the values of $(D^{G'})^{1/2}\theta'$ at $$(\diag(\lambda,\lambda^{-1}),\diag(a\lambda,a\lambda^{-1}),\iota(a,1))$$
and
$$(\diag(a\lambda,a\lambda^{-1}),\diag(\lambda,\lambda^{-1}),\iota(a,1)).$$
The transfer factor is equal to 1 because the quadratic character $\eta_{F\oplus F/F}$ is trivial. If we take the limit as $\lambda\in F^{\times}\rightarrow 1$, we get 
$$4 D^{G'}((I_2,aI_2,\iota(a,1)))^{1/2}c_{\theta'}((I_2,aI_2,\iota(a,1)))$$
$$+4 D^{G'}((aI_2,I_2,\iota(a,1)))^{1/2}c_{\theta'}((aI_2,I_2,\iota(a,1))).$$
Hence $D^H(\nu_{T}(t))c_\theta(\nu_{T}(t))$ is equal to 
$$c_{\theta'}((I_2,aI_2,\iota(a,1)))+c_{\theta'}((aI_2,I_2,\iota(a,1))).$$
This implies that 
\begin{eqnarray*}
&&\frac{1}{2}\sum_{T\in \CT_{ell}(\GU_{1,1})}\int_{T^\ast(F)} D^H(\nu_T(t))c_\theta(\nu_T(t))dt\\
&=&\frac{1}{2}\int_{E^1} c_{\theta'}((I_2,aI_2,\iota(a,1)))+c_{\theta'}((aI_2,I_2,\iota(a,1))) da\\
&=&\int_{E^1} c_{\theta'}((aI_2,I_2,\iota(a,1)))da.
\end{eqnarray*}
Similarly, we can show that (the negative sign comes from the Kottwitz sign between $\GU_{1,1}$ and $\GU_{2,0}$)
\begin{eqnarray*}
&&\frac{1}{2}\sum_{T\in \CT_{ell}(\GU_{2,0})} \int_{T^\ast(F)} D^{H_1}(\nu_{1,T}(t))c_{\theta_1}(\nu_{1,T}(t))dt\\
&=&-\int_{E^1} c_{\theta'}((aI_2,I_2,\iota(a,1)))da.
\end{eqnarray*}

This recovers all the terms in $m_{geom,1}(\theta')$. It remains to show that the summations of the terms correspond to $\CT_{ell}(H),\CT_{ell}(H_1)$ and $\CT_{ell}(H_4)$ are equal to zero. Fix a stable regular elliptic conjugacy class $t_{st}$
of $H$. We use $H(t_{st})$ (resp. $H_i(t_{st})$) to denote the rational conjugacy classes of $H(F)$ (resp. $H_i(F)$) in $t_{st}$.  We only need to show that
\begin{equation}\label{unitary 1}
\sum_{t\in H(t_{st})}\theta(t)-\sum_{i\in \{1,4\},\;t_i\in H_i(t_{st})} \theta_i(t_i)=0.
\end{equation}
We have three cases (the second and third cases only happen in the $p$-adic case).
\begin{itemize}
\item $t_{st}$ is of the form
$$t_{st}=((E,F,a_1)\cup (E,F,a_2))\times ((E,F,b_1)\cup (E,F,b_2))$$ with $a_i,b_i\in E^{\times}$ and 
$$\eta_{E/F}(a_1)=\eta_{E/F}(a_2)=\eta_{E/F}(b_1)=\eta_{E/F}(b_2).$$
\item $t_{st}$ is of the form
$$t_{st}=((E,F,a_1)\cup (E,F,a_2))\times (E',F',b)$$ 
or 
$$t_{st}=(E',F',b) \times ((E,F,a_1)\cup (E,F,a_2))$$ 
where $F'\neq E$ is a quadratic extension of $F$, $E'=F'\otimes_F E$, $a_i\in E^{\times}$ and $b\in (E')^{\times}$ such that 
$$\eta_{E/F}(a_1)=\eta_{E/F}(a_2)=\eta_{E'/F'}(b).$$
\item $t_{st}$ is of the form
$$t_{st}=(E',F',a)\times (E'',F'',b)$$ 
where $F'\neq E$ (resp. $F''\neq E$) is a quadratic extension of $F$, $E'=F'\otimes_F E$, $E''=F''\otimes_F E$, $a\in (E')^{\times}$ and $b\in (E'')^{\times}$ such that 
$$\eta_{E'/F'}(a)=\eta_{E''/F''}(b)\in F^{\times}.$$
\end{itemize}
We will only consider the first case, the remaining two cases follow from a similar and easier argument. 

From now on, assume that 
$$t_{st}=((E,F,a_1)\cup (E,F,a_2))\times ((E,F,b_1)\cup (E,F,b_2))$$ 
with $a_i,b_i\in E^{\times}$ and 
$$\eta_{E/F}(a_1)=\eta_{E/F}(a_2)=\eta_{E/F}(b_1)=\eta_{E/F}(b_2).$$ 
Then the sets $H(t_{st})$ and $H_1(t_{st})\cup H_4(t_{st})$ each contain two elements. As a stable conjugacy class of $G=\GU_{2,2}\times \GU_{1,1}$, $t_{st}$ corresponds to 
$$((E,F,a_1)\cup (E,F,a_2)\cup (E,F,b_1)\cup (E,F,b_2))\times ((E,F,a_1)\cup (E,F,a_2)).$$
Let $\varepsilon_1,\varepsilon_2$ be two elements in $ker(\tr_{E/F})\cap E^{\times}$ that belong to different $Im(N_{E/F})$-orbits. If $\eta_{E/F}(-1)=1$, then the two elements in $H(t_{st})$ are of the form (viewed as conjugacy classes of $\GU_4\times \GU_2$)
\begin{equation}\label{conjugacy class 1}
((E,F,a_1,\varepsilon_1)\cup (E,F,a_2,\varepsilon_1)\cup (E,F,b_1,\varepsilon_1)\cup (E,F,b_2,\varepsilon_1))
\end{equation}
$$\times ((E,F,a_1,\varepsilon_1)\cup (E,F,a_2,\varepsilon_1)),$$
\begin{equation}\label{conjugacy class 2}
((E,F,a_1,\varepsilon_1)\cup (E,F,a_2,\varepsilon_1)\cup (E,F,b_1,\varepsilon_2)\cup (E,F,b_2,\varepsilon_2))
\end{equation}
$$\times ((E,F,a_1,\varepsilon_1)\cup (E,F,a_2,\varepsilon_1))$$
and 
the two elements in $H_1(t_{st})\cup H_4(t_{st})$ are of the form (viewed as conjugacy classes of $\GU_4\times \GU_2$)
\begin{equation}\label{conjugacy class 3}
((E,F,a_1,\varepsilon_1)\cup (E,F,a_2,\varepsilon_2)\cup (E,F,b_1,\varepsilon_1)\cup (E,F,b_2,\varepsilon_2))
\end{equation}
$$\times ((E,F,a_1,\varepsilon_1)\cup (E,F,a_2,\varepsilon_2)),$$
\begin{equation}\label{conjugacy class 4}
((E,F,a_1,\varepsilon_1)\cup (E,F,a_2,\varepsilon_2)\cup (E,F,b_1,\varepsilon_2)\cup (E,F,b_2,\varepsilon_1))
\end{equation}
$$\times ((E,F,a_1,\varepsilon_1)\cup (E,F,a_2,\varepsilon_2)).$$
The stable conjugacy class $t_{st}$ corresponds to 6 stable conjugacy classes of $G'(F)$:
\begin{equation}\label{conjugacy class 5}
((E,F,a_1)\cup (E,F,a_2))\times ((E,F,b_1)\cup (E,F,b_2))\times ((E,F,a_1)\cup (E,F,a_2)),
\end{equation}
\begin{equation}\label{conjugacy class 6}
((E,F,b_1)\cup (E,F,b_2))\times ((E,F,a_1)\cup (E,F,a_2))\times ((E,F,a_1)\cup (E,F,a_2)),
\end{equation}
\begin{equation}\label{conjugacy class 7}
((E,F,a_1)\cup (E,F,b_1))\times ((E,F,a_2)\cup (E,F,b_2))\times ((E,F,a_1)\cup (E,F,a_2)),
\end{equation}
\begin{equation}\label{conjugacy class 8}
((E,F,a_1)\cup (E,F,b_2))\times ((E,F,a_2)\cup (E,F,b_1))\times ((E,F,a_1)\cup (E,F,a_2)),
\end{equation}
\begin{equation}\label{conjugacy class 9}
((E,F,a_2)\cup (E,F,b_1))\times ((E,F,a_1)\cup (E,F,b_2))\times ((E,F,a_1)\cup (E,F,a_2)),
\end{equation}
\begin{equation}\label{conjugacy class 10}
((E,F,a_2)\cup (E,F,b_2))\times ((E,F,a_1)\cup (E,F,b_1))\times ((E,F,a_1)\cup (E,F,a_2)).
\end{equation}

By our definition of the transfer factor, we know that the transfer factors between the conjugacy classes \eqref{conjugacy class 1}, \eqref{conjugacy class 2}, \eqref{conjugacy class 3}, \eqref{conjugacy class 4} and the conjugacy class \eqref{conjugacy class 5} are equal to each other. This implies that  the value of $\theta'$ at the conjugacy class \eqref{conjugacy class 5}
does not contribute to the left hand side of \eqref{unitary 1}. Similarly, we can also show that the value of $\theta'$ at the conjugacy class \eqref{conjugacy class 6} does not contribute to the left hand side of \eqref{unitary 1}.

On the other hand, by our definition of the transfer factor, we know that the transfer factors between the conjugacy classes \eqref{conjugacy class 1}, \eqref{conjugacy class 2} (resp. \eqref{conjugacy class 3}, \eqref{conjugacy class 4}) and the conjugacy class \eqref{conjugacy class 7} are opposite to each other. This implies that  the value of $\theta'$ at the conjugacy class \eqref{conjugacy class 7}
does not contribute to the left hand side of \eqref{unitary 1}. Similarly, we can also show that the values of $\theta'$ at the conjugacy class \eqref{conjugacy class 8}, \eqref{conjugacy class 8} and \eqref{conjugacy class 9} do not contribute to the left hand side of \eqref{unitary 1}. This proves \eqref{unitary 1}.

If $\eta_{E/F}(-1)=-1$, then the two elements in $H_1(t_{st})\cup H_4(t_{st})$ are of the form (viewed as conjugacy classes of $\GU_4\times \GU_2$)
$$((E,F,a_1,\varepsilon_1)\cup (E,F,a_2,\varepsilon_1)\cup (E,F,b_1,\varepsilon_1)\cup (E,F,b_2,\varepsilon_1))$$
$$\times ((E,F,a_1,\varepsilon_1)\cup (E,F,a_2,\varepsilon_1)),$$
$$((E,F,a_1,\varepsilon_1)\cup (E,F,a_2,\varepsilon_1)\cup (E,F,b_1,\varepsilon_2)\cup (E,F,b_2,\varepsilon_2))$$
$$\times ((E,F,a_1,\varepsilon_1)\cup (E,F,a_2,\varepsilon_1))$$
and 
the two elements in $H(t_{st})$ are of the form (viewed as conjugacy classes of $\GU_4\times \GU_2$)
$$((E,F,a_1,\varepsilon_1)\cup (E,F,a_2,\varepsilon_2)\cup (E,F,b_1,\varepsilon_1)\cup (E,F,b_2,\varepsilon_2))$$
$$\times ((E,F,a_1,\varepsilon_1)\cup (E,F,a_2,\varepsilon_2)),$$
$$((E,F,a_1,\varepsilon_1)\cup (E,F,a_2,\varepsilon_2)\cup (E,F,b_1,\varepsilon_2)\cup (E,F,b_2,\varepsilon_1))$$
$$\times ((E,F,a_1,\varepsilon_1)\cup (E,F,a_2,\varepsilon_2)).$$
We can prove \eqref{unitary 1} by a similar argument. This finishes the proof of the proposition.
\end{proof}

\subsection{Some preparation}\label{sec GU(4) preparation}
Let $\Pi_\phi(G)$ be a tempered local $L$-packet of $G(F)$ whose central character is trivial on $Z_{G,H}(F)$ and $\theta_{\Pi_\phi(G)}=\sum_{\pi\in \Pi_\phi(G)}\theta_\pi$ be the distribution character of $\Pi_\phi(G)$. Recall the $\rho_X=\rho_1\oplus \rho_2$ with 
$$\rho_1=\wedge^2\otimes std_2,\;\rho_2=std_4\oplus std_{4}^{\vee}.$$ 
We consider another model associated to $G$. Let $Q=LU$ be a parabolic subgroup of $G$ with $L\simeq (\GL_2(E)\times \GL_1(F))\times \GU_{1,1}$. Up to conjugation there are two generic characters of $U(F)$ corresponding to two Hermitian forms of dimension 2. 
Let $\psi_+$ be the character whose centralizer in $L(F)$ is isomorphic to $\GU_{1,1}(F)\times \GU_{1,1}(F)$. We can diagonally embed  $\GU_{1,1}$ into this centralizer and we will denote its image by $H_0'$. The model $(G,H'=H_0'\ltimes U,\psi_+)$ is essentially the Gan--Gross--Prasad model $(\SO_6\times \SO_3,\SO_3\ltimes U)$ and it is an analogy of the model $(\GL_4\times \GL_2,\GL_2\ltimes U)$ of the previous section for unitary similitude groups. This model has a unique pure inner form $(G_D,H_D'=H_{0,D}'\ltimes U,\psi_-)$ with $G_D=G_1=\GU_{2,2}\times \GU_{2,0}$ and $H_{0,D}'\simeq \GU_{2,0}$. One can easily prove the multiplicity formulas
$$m(\pi)'=c_\pi(1)+\frac{1}{2}\sum_{T\in \CT_{ell}(H_0')}\int_{T(F)/Z_{H_0'}(F)}D^{H'}(t)c_{\pi,\CO_{T,+}}(t)dt,$$
$$m(\pi_D)'=\frac{1}{2}\sum_{T\in \CT_{ell}(H_{0,D}')}\int_{T(F)/Z_{H_{0,D}'}(F)}D^{H_D'}(t)c_{\pi_D,\CO_{T,-}}(t)dt$$
for these models by the same arguments as in the orthogonal Gan--Gross--Prasad models case in \cite{Wal1} and \cite{Wal2}. Here for $t\in T_{reg}(F)$, $\CO_{T,+}$ (resp. $\CO_{T,-}$) is the regular nilpotent orbit in $\Fg_{t}(F)$ (resp. $(\Fg_D)_{t}(F)$) corresponding to the character $\psi_+$ (resp. $\psi_-$). We will denote this model by $(\GU_4\times \GU_2,\GU_2\ltimes U)$.

The above multiplicity formulas imply that 
$$\sum_{\pi\in \Pi_\phi(G)}m(\pi)'+\sum_{\pi_D\in \Pi_\phi(G_D)}m(\pi_D)'$$
is equal to 
$$c_{\theta_{\Pi_\phi(G)}}(1)+\frac{1}{2}\sum_{T\in \CT_{ell}(H_0')}\int_{T(F)/Z_{H_0'}(F)}D^{H'}(t)c_{\theta_{\Pi_\phi(G)},\CO_{T,+}}(t)dt$$
$$+\frac{1}{2}\sum_{T\in \CT_{ell}(H_{0,D}')}\int_{T(F)/Z_{H_{0,D}'}(F)}D^{H_D'}(t)c_{\theta_{\Pi_\phi(G_D)},\CO_{T,-}}(t)dt=c_{\theta_{\Pi_\phi(G)}}(1)=1,$$
i.e. the summation of the multiplicities is equal to 1 over every tempered local $L$-packet. In particular, we know that for an irreducible tempered representation of $G(F)$, the expression $$\frac{1}{2}\sum_{T\in \CT_{ell}(H_0')}\int_{T(F)/Z_{H_0'}(F)}D^{H'}(t)c_{\pi,\CO_{T,+}}(t)dt$$
is equal to 0 or $-1$.

The goal of this section is to prove the following lemma (recall the $\chi_\phi$ is the central character of the $L$-packet of $\GU_4$ obtained from $\Pi_\phi(G)$).

\begin{lemma}\label{lemma-for-GU(4)-GU(2)}
\begin{enumerate}
\item  We have  
$$\int_{T_{0}'(F)}D^H(\nu_{T_0'}(t))c_{\theta_{\Pi_\phi(G)}}(\nu_{T_0'}(t))dt=\frac{\chi_\phi(-1)\epsilon(\frac{1}{2},\Pi_\phi,\rho_2)-1}{2}.$$ 
\item If $\Pi_\phi(G)$ is not a discrete $L$-packet with only one element, or if the central character of $\Pi_\phi(G)$ is trivial, then 
\begin{eqnarray*}
&&\frac{1}{2}\sum_{T\in \CT_{ell}(H_0')}\int_{T(F)/Z_{H_0'}(F)}D^{H'}(t)c_{\theta_{\Pi_\phi(G)}}(t)dt\\
&=&\frac{\chi_\phi(-1)\eta_{E/F}(-1)\epsilon(\frac{1}{2},\Pi_\phi,\rho_1)-1}{2}.
\end{eqnarray*}
\item If $\Pi_\phi(G)$ is not a discrete $L$-packet with only one element, or if the central character of $\Pi_\phi(G)$ is trivial, then 
\begin{eqnarray*}
&&\frac{1}{2}\sum_{T\in \CT_{ell}(H_0')}\int_{T(F)/Z_{H_0'}(F)}D^{H'}(t)c_{\theta_{\Pi_\phi(G)},\CO_{T,+}}(t)dt\\
&=&\frac{\chi_\phi(-1)\eta_{E/F}(-1)\epsilon(\frac{1}{2},\Pi_\phi,\rho_1)-1}{2}.
\end{eqnarray*}
\end{enumerate}
\end{lemma}

\begin{proof}
The first part follows from the multiplicity formulas and the epsilon dichotomy for the Gan--Gross--Prasad model $(U_4\times U_1,U_1\ltimes U)$ proved in  \cite{Beu1}, \cite{Beu3}, \cite{Xue} (we just need to first restrict the $L$-packet to $\GU_{2,2}$, then further restrict to $U_{2,2}$). Note that 
$$\chi_\phi(-1)\epsilon(\frac{1}{2},\Pi_\phi,\rho_2)$$ 
is equal to the epsilon factor of the base change of the packet to $\GL_4(E)$.

The second part is equivalent to the third part since $c_{\theta_{\Pi_\phi(G)},\CO_{T,+}}(t)=c_{\theta_{\Pi_\phi(G)}}(t)$ (this is because the character $\theta_{\Pi_\phi(G)}$ is stable). We first consider the case when the packet is induced from a maximal Levi subgroup $M$. In this case we will prove the second part of the lemma. We use $\Pi_\phi(M)$ to denote the corresponding $L$-packet of $M$ such that the packet $\Pi_\phi(G)$ is induced from the packet $\Pi_\phi(M)$. If $M(F)$ is equal to $\GU_{2,2}(F)\times (\GL_1(E)\times \GL_1(F))$, by Proposition \ref{germ parabolic induction}, we have 
$$\frac{1}{2}\sum_{T\in \CT_{ell}(H_0')}\int_{T(F)/Z_{H_0'}(F)}D^{H'}(t)c_{\theta_{\Pi_\phi(G)}}(t)dt$$
$$=0=\frac{\chi_\phi(-1)\eta_{E/F}(-1)\epsilon(\frac{1}{2},\Pi_\phi,\rho_1)-1}{2}.$$
Note that in this case we can decompose $\rho_1\circ\phi$ as the direct sum of a 6-dimensional representation with its dual and the determinant of the 6-dimensional representation is equal to $\chi_\phi(-1)\eta_{E/F}(-1)$ at $-1$. This implies that $\chi_\phi(-1)\eta_{E/F}(-1)\epsilon(\frac{1}{2},\Pi_\phi,\rho_1)=1$.

If $M(F)$ is isomorphic to $\GL_2(E)\times \GL_1(F)\times \GU_{1,1}(F)$, we may assume that $M=L$. By Proposition \ref{germ parabolic induction}, we have
\begin{eqnarray*}
&&\frac{1}{2}\sum_{T\in \CT_{ell}(H_0')}\int_{T(F)/Z_{H_0'}(F)}D^{H'}(t)c_{\theta_{\Pi_\phi(G)}}(t)dt\\
&=&\frac{1}{2}\sum_{T\in \CT_{ell}(H_0')}\int_{T(F)/Z_{H_0'}(F)}D^{H_0'}(t)\theta_{\Pi_\phi(M)}(t)dt.
\end{eqnarray*}
Meanwhile, the epsilon factor $\chi_\phi(-1)\epsilon(\frac{1}{2},\Pi_\phi,\rho_1)$ is equal to the tensor product epsilon factor of $\GL_2(E)\times \GL_2(F)\times \GL_1(F)$ as in Theorem D of \cite{P92}. To be specific, we can decompose $\rho_1\circ\phi$ as the direct sum of an 8-dimensional representation with two 2-dimensional representations that are dual to each other and whose determinant is equal to $\chi_\phi(-1)$ at $-1$. The epsilon factor corresponds to the 8-dimensional representation is equal to the tensor product epsilon factor of $\GL_2(E)\times \GL_2(F)\times \GL_1(F)$ while the epsilon factor of the remaining 4-dimensional representation is equal to $\chi_\phi(-1)$. Then the second part follows from the epsilon dichotomy (Theorem D of \cite{P92}) and the multiplicity formula (Section 4.5 of \cite{WZ1}) for the generalized trilinear model $(\GL_2(E)\times \GL_2(F),\GL_2(F))$. Note that Section 4.5 of \cite{WZ1} and Theorem D of \cite{P92} only considered the $p$-adic case,  but the multiplicity formula in the Archimedean case follows from a very similar argument as in the $p$-adic case, while the epsilon dichotomy in the Archimedean case follows from the epsilon dichotomy of the Waldspurger model (because all tempered representations of $\GL_2(\BC)$ are principal series).

If $M(F)$ is isomorphic to $\GU_{1,1}(F)\times \GL_1(E)\times \GU_{1,1}(F)$, by Proposition \ref{germ parabolic induction}, we have ($\iota$ is an embedding from $E^{2,0}$ to $\GU_{1,1}(F)$)
\begin{eqnarray*}
&&\frac{1}{2}\sum_{T\in \CT_{ell}(H_0')}\int_{T(F)/Z_{H_0'}(F)}D^{H'}(t)c_{\theta_{\Pi_\phi(G)}}(t)dt\\
&=&\int_{E^1}\theta_{\Pi_\phi(M)}(I_2,a,\iota(a,1))da.
\end{eqnarray*}
By first restricting the packet to the second and third components of $M(F)$ and then further restricting to $U_1(F)\times U_{1,1}(F)$ we get an $L$-packet of $U_1(F)\times U_{1,1}(F)$ and we let $\Pi$ be its base change to $\GL_1(E)\times \GL_2(E)$.  
We can decompose $\rho_1\circ\phi$ as the direct sum of a 4-dimensional representation with two 4-dimensional representations that are dual to each other and whose determinant is equal to $1$ at $-1$. The epsilon factor corresponding to the first 4-dimensional representation is equal to the tensor product epsilon factor of $\Pi$ times $\eta_{E/F}(-1)\chi_\phi(-1)$. Hence $\chi_\phi(-1)\eta_{E/F}(-1)\epsilon(\frac{1}{2},\Pi_\phi,\rho_1)$ is equal to the tensor product epsilon factor of $\Pi$. Then the second part follows from the epsilon dichotomy and the multiplicity formula for the Gan--Gross--Prasad model $(U_2,U_1)$.

It remains to prove the third part for the remaining cases (i.e. when the packet is discrete). We just need to show that if 
$$\chi_\phi(-1)\eta_{E/F}(-1)\epsilon(\frac{1}{2},\Pi_\phi,\rho)=1,$$
$$ \text{(resp.}\; \chi_\phi(-1)\eta_{E/F}(-1)\epsilon(\frac{1}{2},\Pi_\phi,\rho)=-1)$$ 
then the distinguished element for the model $(\GU_4\times \GU_2,\GU_2\ltimes U)$ belongs to $\Pi_\phi(G)$ (resp. $\Pi_\phi(G_D)$).

If the packet is discrete with one element, then $F$ must be $p$-adic and by assumption we know the central character is trivial. In this case under the lower rank isomorphisms $\mathrm{PGU}_{2,2}\simeq \mathrm{PGSO}_6$ and $\mathrm{PGU}_{1,1}\simeq \mathrm{SO}_3$, $\Pi_\phi(G)$ induces an $L$-packet of $\SO_6\times \SO_3$ (here $\SO_3$ is the split odd special orthogonal group of rank $1$ and $\SO_6$ is a quasi-split but not split even special orthogonal group of rank 3). Then the second part follows from the multiplicity formulas and the epsilon dichotomy for the Gan--Gross--Prasad model $(\SO_6\times \SO_3,\SO_3\ltimes U)$ (\cite{Wal1}, \cite{Wal2}, \cite{Wal3}).

The last case is when $\Pi_\phi(G)$ is discrete and contains more than one element. In this case, the centralizer $Z_\phi$ contains an element that does not belong to the center. Hence we can find a proper elliptic extended endoscopic triple $(G',s',{}^L\eta)$ such that $\phi$ factors through ${}^L\eta$, $G'=G(U_{1,1}\times U_{1,1})\times \GU_{1,1}$, $s'=(s_1,1,I_2,1)\in Z_\phi$ with $s_1$ being conjugated to $\diag(I_2,-I_2)$. We use $\Pi_\phi(G')$ to denote the associated $L$-packet of $G'$. We know that the character 
$$\sum_{\pi\in \Pi_\phi(G)} \chi_\pi(s') \theta_\pi$$ 
of $G(F)$ is the transfer of the stable character $$\theta_{\Pi_\phi(G')}=\sum_{\tau\in \Pi_\phi(G')}\theta_\tau$$ 
of $G'(F)$.

By the same argument as in the orthogonal Gan--Gross--Prasad model case in Section 3.3 of \cite{Wal3}, we can prove that
$$\sum_{\pi\in \Pi_\phi(G)} \chi_\pi(s')\cdot(c_{\theta_\pi}(1)+ \frac{1}{2}\sum_{T\in \CT_{ell}(H_0')}\int_{T(F)/Z_{H_0'}(F)}D^{H'}(t)c_{\pi,\CO_{T,+}}(t)dt)$$
is equal to (recall that $T^\ast(F)=T(F)/Z_{\GU_{1,1}}(F)=T(F)/Z_{\GU_{2,0}}(F)$)
\begin{align*}
& c_{\theta_{\Pi_\phi(G')}}(1)
+\int_{E^1} c_{\theta_{\Pi_\phi(G')}}(aI_2,I_2,\iota(a,1))da\\
&+\frac{1}{2}\sum_{T\in \CT_{ell}(\GU_{1,1})}\int_{T^\ast(F)} D^{\GU_{1,1}}(t) \theta_{\Pi_\phi(G')}(t,t,t)dt.
\end{align*}

Meanwhile, as in Section \ref{sec epsilon factor}, we can decompose $\rho_1\circ \phi$ as 
$$\rho_1\circ \phi=\rho_{1,s',\phi}\oplus \rho_{1,s',\phi}'.$$ The packet $\Pi_\phi(G')$ induces a packet of $U_{1,1}$ by restricting to the third copy of $G'$. By the epsilon dichotomy and the multiplicity formula for the Gan--Gross--Prasad  model $(U_2,U_1)$, we have
$$\int_{E^1} c_{\theta_{\Pi_\phi(G')}}(aI_2,I_2,\iota(a,1))da=\frac{\chi_\phi(-1)\eta_{E/F}(-1)\epsilon(\frac{1}{2},\rho_{1,s',\phi}')-1}{2}.$$

On the other hand, by Theorem 8.1 of \cite{La}, the packet $\Pi_\phi(G')$ is the restriction of an irreducible tempered representation $\Pi'$ of $\GU_{1,1}(F)\times \GU_{1,1}(F)\times \GU_{1,1}(F)$ to $G'(F)$ (the choice of $\Pi'$ is not unique) and $\Pi'$ induces an irreducible tempered representation $\Pi$ of $\GL_2(F)\times \GL_2(F)\times \GL_2(F)$. 
By the epsilon dichotomy (\cite{P90}, \cite{L01}) and the multiplicity formula (\cite{Wan16}) for the trilinear $\GL_2$ model $((\GL_2)^3,\GL_2)$ , we have 
$$\frac{1}{2}\sum_{T\in \CT_{ell}(\GU_{1,1})}\int_{T^\ast(F)} D^{\GU_{1,1}}(t) \theta_{\Pi_\phi(G')}(t,t,t)dt=\frac{\epsilon(\frac{1}{2},\rho_{1,s',\phi})-1}{2}.$$
This implies that 
$$\sum_{\pi\in \Pi_\phi(G)} \chi_\pi(s')m'(\pi)=\frac{\epsilon(\frac{1}{2},\rho_{1,s',\phi})+\chi_\phi(-1)\eta_{E/F}(-1)\epsilon(\frac{1}{2},\rho_{1,s',\phi}')}{2}.$$
If $\chi_\phi(-1)\eta_{E/F}(-1)\epsilon(\frac{1}{2},\Pi_\phi,\rho_1)=1$, then 
$$\frac{\epsilon(\frac{1}{2},\rho_{1,s',\phi})+\chi_\phi(-1)\eta_{E/F}(-1)\epsilon(\frac{1}{2},\rho_{1,s',\phi}')}{2}= \epsilon(\frac{1}{2},\rho_{1,s',\phi})\neq 0$$
and hence the unique distinguished element belongs to $\Pi_\phi(G)$. If $\chi_\phi(-1)\eta_{E/F}(-1)\epsilon(\frac{1}{2},\Pi_\phi,\rho_1)=-1$, then 
$$\frac{\epsilon(\frac{1}{2},\rho_{1,s',\phi})+\chi_\phi(-1)\eta_{E/F}(-1)\epsilon(\frac{1}{2},\rho_{1,s',\phi}')}{2}= 0$$
and hence the unique distinguished element belongs to $\Pi_\phi(G_D)$. This proves the lemma. 
\end{proof}

It is clear that for a tempered $L$-packet $\Pi_\phi(G)$ of $G(F)$ whose central character is trivial on $Z_{G,H}(F)$, the identity
\begin{eqnarray*}
&&\frac{1}{2}\sum_{T\in \CT_{ell}(H_0')}\int_{T(F)/Z_{H_0'}(F)}D^{H'}(t)c_{\theta_{\Pi_\phi(G)}}(t)dt\\
&=&\frac{\chi_\phi(-1)\eta_{E/F}(-1)\epsilon(\frac{1}{2},\Pi_\phi,\rho_1)-1}{2}
\end{eqnarray*}
is equivalent to the following conjecture which is an analogy of Conjecture \ref{weak conjecture} for the model $(G,H_0\ltimes U)=(\GU_4\times \GU_2,\GU_2\ltimes U)$.

\begin{conj}\label{weak conjecture GU(4)}
The unique distinguished element for the model $(\GU_4\times \GU_2,\GU_2\ltimes U)$ belongs to $\Pi_\phi(G)$ (resp. $\Pi_\phi(G_D)$) if and only if  $$\chi_\phi(-1)\eta_{E/F}(-1)\epsilon(\frac{1}{2},\Pi_\phi,\rho_1)=1$$ 
$$\text{(resp.}\; \chi_\phi(-1)\eta_{E/F}(-1)\epsilon(\frac{1}{2},\Pi_\phi,\rho_1)=-1).$$
\end{conj}

\subsection{The proof of Theorem \ref{main theorem} for $(\GU_4\times \GU_2,(\GU_2\times \GU_2)^0)$}\label{sec GU(4) 2}
In this subsection we will prove Theorem \ref{main theorem} for $(\GU_4\times \GU_2,(\GU_2\times \GU_2)^0)$. Let $\Pi_\phi=\Pi_\phi(G)\cup \Pi_\phi(G_i)$ be a tempered $L$-packet of $G/Z_{G,H}$ and let $\theta_{\Pi_\phi(G)}$ (resp. $\theta_{\Pi_\phi(G_i)}$) be the distribution character of the packet $\Pi_\phi(G)$ (resp. $\Pi_\phi(G_i)$). Then we know that $\theta_{\Pi_\phi(G)}$ is the transfer of $a_i\theta_{\Pi_\phi(G_i)}$ where $a_3=1$ and $a_1=a_2=a_4=-1$. Combining with the multiplicity formula, we have (note that $c_{\theta_{\Pi_\phi(G)}}(1)=1$ since there is a unique generic element in the packet)
\begin{eqnarray*}
&&m_{geom}(\theta_{\Pi_\phi(G)})+m_{geom}(\theta_{\Pi_\phi(G_1)})+m_{geom}(\theta_{\Pi_\phi(G_4)})\\
&=&1+\int_{T_{0}'(F)}D^H(\nu_{T_0'}(t))c_{\theta_{\Pi_\phi(G)}}(\nu_{T_0'}(t))dt,\\
&&m_{geom}(\theta_{\Pi_\phi(G)})+m_{geom}(\theta_{\Pi_\phi(G_2)})\\
&=&1+\frac{1}{2}\sum_{T\in \CT_{ell}(\GU_{1,1})}\int_{T^\ast(F)}D^H(\nu_T(t))c_{\theta_{\Pi_\phi(G)}}(\nu_T(t))dt,\\
&&m_{geom}(\theta_{\Pi_\phi(G_3)})+m_{geom}(\theta_{\Pi_\phi(G_1)})+m_{geom}(\theta_{\Pi_\phi(G_4)})\\
&=&-\frac{1}{2}\sum_{T\in \CT_{ell}(\GU_{1,1})}\int_{T^\ast(F)}D^H(\nu_T(t))c_{\theta_{\Pi_\phi(G)}}(\nu_T(t))dt,
\end{eqnarray*}
$$m_{geom}(\theta_{\Pi_\phi(G_3)})+m_{geom}(\theta_{\Pi_\phi(G_2)})=-\int_{T_{0}'(F)}D^H(\nu_{T_0'}(t))c_{\theta_{\Pi_\phi(G)}}(\nu_{T_0'}(t))dt.$$

From now on, assume that either $\Pi_\phi$ is not a discrete L-packet with $|\Pi_\phi(G)|=1$ or the central character of $\Pi_\phi(G)$ is trivial. Combining the above equations with Lemma \ref{lemma-for-GU(4)-GU(2)}, we have
\begin{equation}\label{GU(4) equation 1}
m_{geom}(\theta_{\Pi_\phi(G)})=1\iff 
\end{equation}
$$\chi_\phi(-1)\eta_{E/F}(-1)\epsilon(\frac{1}{2},\Pi_\phi,\rho_1)=\chi_\phi(-1)\epsilon(\frac{1}{2},\Pi_\phi,\rho_2)=1,$$
$$m_{geom}(\theta_{\Pi_\phi(G_1)})+m_{geom}(\theta_{\Pi_\phi(G_4)})=1\iff$$
$$-\chi_\phi(-1)\eta_{E/F}(-1)\epsilon(\frac{1}{2},\Pi_\phi,\rho_1)=\chi_\phi(-1)\epsilon(\frac{1}{2},\Pi_\phi,\rho_2)=1,$$
$$m_{geom}(\theta_{\Pi_\phi(G_2)})=1\iff$$ 
$$\chi_\phi(-1)\eta_{E/F}(-1)\epsilon(\frac{1}{2},\Pi_\phi,\rho_1)=-\chi_\phi(-1)\epsilon(\frac{1}{2},\Pi_\phi,\rho_2)=1,$$
$$m_{geom}(\theta_{\Pi_\phi(G_3)})=1\iff $$
$$\chi_\phi(-1)\eta_{E/F}(-1)\epsilon(\frac{1}{2},\Pi_\phi,\rho_1)=\chi_\phi(-1)\epsilon(\frac{1}{2},\Pi_\phi,\rho_2)=-1.$$
This proves Conjecture \ref{weak conjecture GU(4)xGU(2)} for the packet $\Pi_\phi$.

Let $\omega_\phi$ be the character of $S_\phi$ corresponding to the unique distinguished element of the $L$-packet and we also view it as a character of $Z_\phi$. Fix $s\in S_\phi$. By Lemma \ref{lem extended endoscopic triple}, there exists an elliptic extended endoscopic triple $(G',s',{}^L\eta)$ of $G/Z_{G,H}$ such that $s'\in sZ_{\phi}^{\circ}$ and $\phi$ factors through ${}^L\eta$. We need to show that $\omega_{\phi,H}(s)=\omega_\phi(s')$. 

By the above relation and the definition of $\omega_{\phi,H}$ we know that $\omega_{\phi,H}(s)=\chi_\phi(s')$ if $s'$ belongs to the center of the dual group. Then we consider the case when $s'$ does not belong to the center of the dual group. We will only study the case when 
$$s'=(s_1,I_2,1)\in \widehat{G/Z_{G,H}}=\GL_4(\BC)\times \GL_2(\BC)\times \GL_1(\BC).$$ 
The case when $s'=(s_1,-I_2,1)$ follows from a similar argument (recall that $s_1$ is an element in $\GL_4(\BC)$ which is conjugate to $\diag(I_2,-I_2)$).

We first consider the case when the distinguished element belongs to $\Pi_\phi(G)\cup \Pi_\phi(G_1)\cup \Pi_\phi(G_4)$. This implies that 
$$\chi_\phi(-1)\eta_{E/F}(-1)\epsilon(\frac{1}{2},\Pi_\phi,\rho_1)\in \{\pm 1\}, \;\chi_\phi(-1)\epsilon(\frac{1}{2},\Pi_\phi,\rho_2)=1.$$
As in Section \ref{sec epsilon factor} we have a decomposition (we refer the reader to Section \ref{sec epsilon factor} for various notation)
$$\rho_1\circ\phi=\rho_{1,s',\phi}\oplus \rho_{1,s',\phi}',\;\rho_2\circ \phi=\rho_{2,s',\phi,+}\oplus \rho_{2,s',\phi,-}$$
and the equation
$$\omega_{\phi,H}(s)=\eta_{E/F}(-1)\chi_{\phi,s',2}(-1)\epsilon(\frac{1}{2},\rho_{1,s',\phi}\oplus \rho_{2,s',\phi,-}).$$
The Langlands parameter $\phi$ induces a parameter $\phi$ of $G''=G(U_{1,1}\times U_{1,1})\times \GU_{1,1}$ and we let $\Pi_\phi(G'')$ be the corresponding $L$-packet. 

By the formula of endoscopy in Proposition \ref{prop GU(4)} together with the epsilon dichotomy and the multiplicity formula for the Gan--Gross--Prasad model $(U_2,U_1)$, we know that 
$$\sum_{\pi\in \Pi_\phi(G)}\chi_\pi(s')m(\pi)-\sum_{\pi_1\in \Pi_\phi(G_1)} \chi_{\pi_1}(s')m(\pi_1)-\sum_{\pi_4\in \Pi_\phi(G_4)} \chi_{\pi_4}(s')m(\pi_4)$$
is equal to
$$1+\frac{\eta_{E/F}(-1)\chi_{\phi,s',1}(-1)\epsilon(\frac{1}{2},\rho_{2,s',\phi,+})-1}{2}$$
$$+\frac{\eta_{E/F}(-1)\chi_{\phi,s',2}(-1)\epsilon(\frac{1}{2},\rho_{2,s',\phi,-})-1}{2}$$
$$+(\eta_{E/F}(-1)\chi_{\phi}(-1)\epsilon(\frac{1}{2},\rho_{1,s',\phi}')-1)$$
$$+\frac{(\eta_{E/F}(-1)\chi_{\phi}(-1)\epsilon(\frac{1}{2},\rho_{1,s',\phi}')-1)(\eta_{E/F}(-1)\chi_{\phi,s',1}(-1)\epsilon(\frac{1}{2},\rho_{2,s',\phi,+})-1)}{2}$$
$$+\frac{(\eta_{E/F}(-1)\chi_{\phi}(-1)\epsilon(\frac{1}{2},\rho_{1,s',\phi}')-1)(\eta_{E/F}(-1)\chi_{\phi,s',2}(-1)\epsilon(\frac{1}{2},\rho_{2,s',\phi,-})-1)}{2}$$
$$=\frac{\chi_{\phi}(-1)\epsilon(\frac{1}{2},\rho_{1,s',\phi}')\cdot (\chi_{\phi,s',1}(-1)\epsilon(\frac{1}{2},\rho_{2,s',\phi,+})+\chi_{\phi,s',2}(-1)\epsilon(\frac{1}{2},\rho_{2,s',\phi,-}))}{2}.$$
Since $\chi_\phi(-1)\epsilon(\frac{1}{2},\Pi_\phi,\rho_2)=1$, we have $$\chi_{\phi,s',1}(-1)\epsilon(\frac{1}{2},\rho_{2,s',\phi,+})=\chi_{\phi,s',2}(-1)\epsilon(\frac{1}{2},\rho_{2,s',\phi,-})$$
which implies that 
\begin{equation}\label{model GU4 equation 1}
\sum_{\pi\in \Pi_\phi(G)}\chi_\pi(s')m(\pi)-\sum_{\pi_1\in \Pi_\phi(G_1)} \chi_{\pi_1}(s')m(\pi_1)-\sum_{\pi_4\in \Pi_\phi(G_4)} \chi_{\pi_4}(s')m(\pi_4)
\end{equation}
is equal to
$$\chi_{\phi}(-1)\epsilon(\frac{1}{2},\rho_{1,s',\phi}')\cdot \chi_{\phi,s',2}(-1)\epsilon(\frac{1}{2},\rho_{2,s',\phi,-}).$$

If $\chi_\phi(-1)\eta_{E/F}(-1)\epsilon(\frac{1}{2},\Pi_\phi,\rho_1)=1$, the distinguished element belongs to $\Pi_\phi(G)$ and we have $$\chi_\phi(-1)\epsilon(\frac{1}{2},\rho_{1,s',\phi}')=\eta_{E/F}(-1)\epsilon(\frac{1}{2},\rho_{1,s',\phi}).$$ 
This implies that $\omega_\phi(s')$ is equal to \eqref{model GU4 equation 1} which is equal to
$$\chi_{\phi,s',2}(-1)\eta_{E/F}(-1)\epsilon(\frac{1}{2},\rho_{1,s',\phi})\cdot \epsilon(\frac{1}{2},\rho_{2,s',\phi,-})=\omega_{\phi,H}(s).$$

If $\chi_\phi(-1)\eta_{E/F}(-1)\epsilon(\frac{1}{2},\Pi_\phi,\rho_1)=-1$, the distinguished element belongs to $\Pi_\phi(G_1)\cup \Pi_\phi(G_4)$  and we have $$\chi_\phi(-1)\epsilon(\frac{1}{2},\rho_{1,s',\phi}')=-\eta_{E/F}(-1)\epsilon(\frac{1}{2},\rho_{1,s',\phi}).$$ 
This implies that $\omega_\phi(s')$ is equal to $-1$ times \eqref{model GU4 equation 1} which is equal to
$$\chi_{\phi,s',2}(-1)\eta_{E/F}(-1)\epsilon(\frac{1}{2},\rho_{1,s',\phi})\cdot \epsilon(\frac{1}{2},\rho_{2,s',\phi,-})=\omega_{\phi,H}(s).$$
This proves the identity $\omega_\phi(s')=\omega_{\phi,H}(s)$ when the distinguished element belongs to $\Pi_\phi(G)\cup \Pi_\phi(G_1)\cup \Pi_\phi(G_4)$. The argument for the case when the distinguished element belongs to $\Pi_\phi(G_2)\cup \Pi_\phi(G_3)$ is similar (we just need use the second equation in Proposition \ref{prop GU(4)}) and we will skip it here. This finishes the proof of the theorem.

\subsection{The proof of Theorem \ref{main theorem} and \ref{thm weak conjecture smaller models} for $(\GU_6,\GU_2\ltimes U)$}\label{sec GU(6)}
In this subsection we will prove Theorem \ref{main theorem} for the model $(\GU_6,\GU_2\ltimes U)$. We first recall the definition of the model. Let $G=\GU_{3,3}$, and $P=LU$ be the standard parabolic subgroup of $G$ with
\begin{align*}
L(F)=&\{m(g,h)=\begin{pmatrix}
g&&\\&h&\\&&l(h)g^*	\end{pmatrix} \mid \\
&g\in \GL_2(E),\; g^*=w_2{}^t\bar{g}^{-1}w_2,~h\in\GU_{1,1}(F)\},
\end{align*}
Let $\xi$ be a generic character of $U(F)$ given by
$$
\xi(u(X,Y))=\psi(\lambda(u(X,Y))),\;\lambda(u(X,Y))=\tr_{E/F}(\tr(X)).
$$
Then the stabilizer of $\xi$ under the adjoint action of $L(F)$ is
$$
H_{0}(F):=\{m(h,h)\mid h\in \GU_{1,1}(F)\}=\{diag(h,h,h) \mid h\in \GU_{1,1}(F)\}.
$$
Let $H=H_0\ltimes U$ and we extend the character $\xi$ to $H(F)$ by making it trivial on $H_0(F)$. The model $(G,H,\xi)$ is the analogue of the Ginzburg--Rallis model. We can also define the quaternion (non quasi-split) version of this model by letting $G_D=\GU_{4,2}$ be the non quasi-split unitary similitude group. In this case, we have $H_{0,D}=\GU_{2,0}$.

For a quasi-character $\theta$ (resp. $\theta_D$) of $G(F)$ (resp. $G_D(F)$), define the geometric multiplicities
$$m_{geom}(\theta)=c_\theta(1)+\sum_{T\in \CT_{ell}(H_0)}|W(H_0,T)|^{-1} \int_{T(F)/Z_{G,H}(F)}D^H(t) c_{\theta}(t)dt,$$
\begin{eqnarray*}
m_{geom}(\theta_D)&=&\sum_{T_D\in \CT_{ell}(H_{0,D})} |W(H_{0,D},T_D)|^{-1}\\
&&\cdot\int_{T_D(F)/Z_{G_D,H_D}(F)} D^{H_D}(t) c_{\theta_D}(t)dt.
\end{eqnarray*}
In our previous paper \cite{WZ1}, we have proved the multiplicity formulas
$$m(\pi)=m_{geom}(\theta_\pi),\;m(\pi_D)=m_{geom}(\theta_{\pi_D})$$
for all tempered representations $\pi$ (resp. $\pi_D$) of $G(F)$ (resp. $G_D(F)$) in the $p$-adic case. For the rest of this section we will assume that the multiplicity formulas hold for both the $p$-adic case and the real case. Now we study the behaviors of the geometric multiplicities under endoscopy and under parabolic induction. We first define the geometric multiplicities associated to parabolic subgroups and elliptic endoscopic groups.

Let $M$ be a proper maximal Levi subgroup of $G$ and $\theta^M$ be a quasi-character on $M(F)$. If $M$ is the Siegel Levi subgroup of $G$, define $m_{geom}(\theta^M)=c_{\theta^M}(1)$. Otherwise, $M$ corresponds to a proper Levi subgroup $M_D$ of $G_D$. Let $\theta_{D}^{M_D}$ be a quasi-character on $M_D(F)$. If $M=L$ (and hence $M_D=L_D$), define
\begin{eqnarray*}
m_{geom}(\theta^M)&=&c_{\theta^M}(1)+\sum_{T\in \CT_{ell}(H_0)}|W(H_0,T)|^{-1}\\
&&\cdot \int_{T(F)/Z_{G,H}(F)} D^{H_0}(t) \theta^M(t)dt,\\
m_{geom}(\theta_{D}^{M_D})&=&\sum_{T_D\in \CT_{ell}(H_{0,D})} |W(H_{0,D},T_D)|^{-1} \\
&&\cdot \int_{T_D(F)/Z_{G_D,H_D}(F)} D^{H_{0,D}}(t) \theta_{D}^{M_D}(t)dt.
\end{eqnarray*}

If $M(F)\simeq \GU_{2,2}(F)\times \GL_1(E)$, let $T_0(F)\simeq E^{2,0}$ be a maximal elliptic torus of $\GU_{1,1}(F)$ and we fix an isomorphism $T_0(F)\simeq E^{2,0}$. We embed it into $M(F)$ via the map:
$$\nu(t)=\diag(a,b,t,b,a),\;t\in T_0(F) \text{ corresponds to } (a,b)\in  E^{2,0}.$$
Note that the image of the embedding is contained in the Levi subgroup $\GU_{1,1}(F)\times \GL_1(E)\times \GL_1(E)$. Similarly, we can also define an embedding from a maximal elliptic torus $T_{0,D}(F)\simeq E^{2,0}$ of $\GU_{2,0}(F)$ into $M_D(F)$ which is denoted by $\nu_D$. 
Define
$$m_{geom}(\theta^M)=c_{\theta^M}(1)+\int_{T_0(F)/Z_{\GU_{1,1}}(F)} D^{\GU_{1,1}}(t)c_{\theta^M}(\nu(t))dt,$$
$$m_{geom}(\theta_{D}^{M_D})=\int_{T_{0,D}(F)/Z_{\GU_{2,0}}(F)} D^{\GU_{2,0}}(t)c_{\theta_{D}^{M_D}}(\nu_D(t))dt.$$
The following proposition is a direct consequence of Proposition \ref{germ parabolic induction}.

\begin{prop}\label{parabolic GU(6)}
Let $\theta$ (resp. $\theta_D$) be a quasi-character on $G(F)$ (resp. $G_D(F)$). Assume that  $\theta$ (resp. $\theta_D$) is the parabolic induction of a quasi-character $\theta^M$ (resp. $\theta_{D}^{M_D}$) of a proper maximal Levi subgroup $M$ of $G$ (resp. $M_D$ of $G_D$). We have
$$m_{geom}(\theta)=m_{geom}(\theta^M),\;m_{geom}(\theta_D)=m_{geom}(\theta_{D}^{M_D}).$$
\end{prop}

Next we study the behavior of the geometric multiplicities under endoscopic transfer. Let $(G',s',{}^L\eta)$ be a proper extended endoscopic triple with $G'=G(U_{1,1}\times U_{2,2})$, $s'=(\diag(I_2,-I_4),1)\in \hat{G}=\GL_6(\BC)\times \GL_1(\BC)$ and ${}^L\eta$ be the natural embedding. 
Let $\theta'$ be a quasi-character on $G'(F)$. Using the diagonal embedding from $\GU_{1,1}$ to $\GU_{2,2}$ in the previous case we get a diagonal embedding from $\GU_{1,1}$ to $G'$ (denoted by $\nu'$). Like in the previous case, we have an embedding, denoted by $\nu$, from $E^1$ into $G'(F)$ given by $a\mapsto \diag(1,\iota(a,1),1)\times aI_2$. 
We define
\begin{eqnarray*}
m_{geom}(\theta')&=&c_{\theta'}(1)+\sum_{T\in \CT_{ell}(\GU_{1,1})} \frac{1}{2}\int_{T^\ast(F)} D^{\GU_{1,1}}(t)^2 c_{\theta'}(\nu'(t))dt\\
&&+\int_{E^1} D^{\GU_{1,1}}(\iota(a,1))c_{\theta'}(\nu(a))da,\\
m_{geom,D}(\theta')&=&\sum_{T\in \CT_{ell}(\GU_{1,1})} \frac{1}{2}\int_{T^\ast(F)} D^{\GU_{1,1}}(t)^2 c_{\theta'}(\nu'(t))dt\\
&&-\int_{E^1} D^{\GU_{1,1}}(\iota(a,1))c_{\theta'}(\nu(a))da.
\end{eqnarray*}

\begin{prop}\label{prop GU(6)}
Let $\theta$ (resp. $\theta_D$) be a quasi-character on $G(F)$ (resp. $G_D(F)$). Assume that $\theta$ (resp. $\theta_D$) is the endoscopic transfer of a stable quasi-character $\theta'$ of $G'(F)$ . We have
$$m_{geom}(\theta)=m_{geom}(\theta'),\;m_{geom}(\theta_D)=m_{geom,D}(\theta').$$
\end{prop}

\begin{proof}
We will only prove the quasi-split case, the quaternion case follows from a similar argument. Recall that 
$$m_{geom}(\theta)=c_\theta(1)+\sum_{T\in \CT_{ell}(H_0)}|W(H_0,T)|^{-1} \int_{T(F)/Z_{G,H}(F)}D^H(t) c_{\theta}(t)dt.$$
The proof of the identity $c_\theta(1)=c_{\theta'}(1)$ is easy and we will skip it here. Now we fix $T\in \CT_{ell}(H_0)$ and we will study the term corresponding to $T$ in $m_{geom}(\theta)$. The element $t\in T(F)\subset G(F)$ is of the form $\diag(t_0,t_0,t_0)$ with $t_0\in \GU_{1,1}(F)$ belongs to the torus that is isomorphic to $T$. There is a natural bijection $T\leftrightarrow F_T$ between $\CT_{ell}(H_0)$ and the set of quadratic extensions of $F$. 

If $F_T\neq E$, then $E_T=F_T\otimes_F E$ is a quadratic extension of $E$ and we can identify $t_0$ with an element in $E_{T}^{\times}$ whose norm (with respect to the quadratic extension $E_T/F_T$) belongs to $F^{\times}$. By Proposition \ref{germ parabolic induction}, we know that (assume that $t\in T_{reg}(F)$, i.e. $t$ does not belong to the center) $D^{H}(t)c_{\theta}(t)$ is equal to
$$D^{\GU_{1,1}}(t_0)^{-1/2}D^{G}(t)^{1/2}c_{\theta}(t)=\frac{1}{2}\lim_{\lambda\in F^{\times},\lambda\rightarrow 1} D^{\GU_{1,1}}(t_0)^{-1/2}$$
$$D^{G}(\diag(\lambda t_0,t_0,\lambda^{-1}t_0))^{1/2} \theta(\diag(\lambda t_0, t_0,\lambda^{-1}t_0)).$$
Under the notation of Section \ref{section transfer factor}, the conjugacy class $\diag(\lambda t_0, t_0,\lambda^{-1}t_0)$ is of the form (note that $c_i$ is unique in this case and hence we will ignore it)
$$(E_T,F_T,t_0)\cup (E_T\oplus E_T, E_T,(\lambda t_0,\lambda^{-1}\bar{t}_0)).$$
This conjugacy class corresponds to a unique conjugacy class in $G'=G(U_{1,1}\times U_{2,2})$ given by 
$$(E_T,F_T,t_0)\times  (E_T\oplus E_T, E_T,(\lambda t_0,\lambda^{-1}\bar{t}_0))$$
and the transfer factor is trivial since the quadratic character $\eta_{E_T\oplus E_T/E_T}$ is trivial. As a result, we know that $D^{H}(t)c_{\theta}(t)$ is equal to 
\begin{eqnarray*}
&&\frac{1}{2}\lim_{\lambda\in F^{\times},\lambda\rightarrow 1} D^{\GU_{1,1}}(t_0)^{-1/2} D^{G}(\diag(\lambda t_0,t_0,\lambda^{-1}t_0))^{1/2}\\
&&\cdot \theta(\diag(\lambda t_0, t_0,\lambda^{-1}t_0))\\
&=&\frac{1}{2}\lim_{\lambda\in F^{\times},\lambda\rightarrow 1} D^{\GU_{1,1}}(t_0)^{-1/2}D^{G'}(t_0\times \diag(\lambda t_0,\lambda^{-1}t_0))^{1/2}\\
&&\cdot\theta'(t_0\times \diag(\lambda t_0,\lambda^{-1}t_0))\\
&=&D^{\GU_{1,1}}(t_0)^{-1/2}D^{G'}(t_0\times \diag(t_0,t_0))^{1/2} c_{\theta'}(t_0\times \diag(t_0,t_0))\\
&=&D^{\GU_{1,1}}(t_0)^2 c_{\theta'}(\nu'(t_0))
\end{eqnarray*}
where the second equality follows from Proposition \ref{germ parabolic induction}. This implies that the terms in $m_{geom}(\theta)$ and $m_{geom}(\theta')$ associated to $T$ are equal to each other.

If $E_T=E$, then we can identify $t_0$ with an element in $(a,b)\in E^{2,0}$. By Proposition \ref{germ parabolic induction}, we know that (assume that $t\in T_{reg}(F)$) $D^{H}(t)c_{\theta}(t)$ is equal to $\frac{1}{4}D^{\GU_{1,1}}(t_0)^{-1/2}$ times the limit of $D^G(\cdot)^{1/2}c_\theta(\cdot)$ at the conjugacy class 
$$(E,F,a)\cup (E,F,b)\cup (E\oplus E, E,(\lambda_1 a,\lambda_{1}^{-1}\bar{a}))\cup (E\oplus E, E,(\lambda_2 b,\lambda_{2}^{-1}\bar{b}))$$
as $\lambda_i\rightarrow 1$. In this case $c_i$ is again unique and we will ignore it. This conjugacy class corresponds to three conjugacy classes in $G'=G(U_{1,1}\times U_{2,2})$ given by 
$$((E,F,a)\cup (E,F,b))\times ((E\oplus E, E,(\lambda_1 a,\lambda_{1}^{-1}\bar{a}))\cup (E\oplus E, E,(\lambda_2 b,\lambda_{2}^{-1}\bar{b}))), $$
$$(E\oplus E, E,(\lambda_1 a,\lambda_{1}^{-1}\bar{a}))\times ((E,F,a)\cup (E,F,b)\cup (E\oplus E, E,(\lambda_2 b,\lambda_{2}^{-1}\bar{b}))),$$
$$(E\oplus E, E,(\lambda_2 b,\lambda_{2}^{-1}\bar{b}))\times ((E,F,a)\cup (E,F,b)\cup (E\oplus E, E,(\lambda_1 a,\lambda_{1}^{-1}\bar{a}))).$$
The transfer factor for the first conjugacy class is trivial since the quadratic character $\eta_{E\oplus E/E}$ is trivial. Moreover, by the same argument as in the previous case, this recovers the term in $m_{geom}(\theta')$ associated to $T$. 

For the second and third conjugacy classes,  the transfer factor is still trivial by Section 1.11 of \cite{Wal}. The second and third conjugacy classes give  us the expression
$$\frac{1}{4}\lim_{\lambda_i\in F^{\times},\lambda_i\rightarrow 1} D^{\GU_{1,1}}(t_0)^{-1/2}D^{G'}(\diag(a\lambda_1,a\lambda_{1}^{-1})\times \diag(b\lambda_2, t_0,\lambda_{2}^{-1}b))^{1/2}$$
$$\cdot \theta'(\diag(a\lambda_1,a\lambda_{1}^{-1})\times \diag(b\lambda_2, t_0,\lambda_{2}^{-1}b))$$
$$+\frac{1}{4}\lim_{\lambda_i\in F^{\times},\lambda_i\rightarrow 1} D^{\GU_{1,1}}(t_0)^{-1/2}D^{G'}(\diag(b\lambda_2,b\lambda_{2}^{-1})\times \diag(a\lambda_1, t_0,\lambda_{1}^{-1}a))^{1/2}$$ 
$$\theta'(\diag(b\lambda_2,b\lambda_{2}^{-1})\times \diag(a\lambda_1, t_0,\lambda_{1}^{-1}a))$$
$$=D^{\GU_{1,1}}(t_0)^{-1/2}D^{G'}(aI_2\times \diag(b, t_0,b))^{1/2}c_{\theta'}(aI_2\times \diag(b, t_0,b))$$ 
$$+D^{\GU_{1,1}}(t_0)^{-1/2}D^{G'}(bI_2\times \diag(a, t_0,a))^{1/2} c_{\theta'}(bI_2\times \diag(a, t_0,a))$$
$$=D^{\GU_{1,1}}(t_0)c_{\theta'}(aI_2\times \diag(b, t_0,b))+D^{\GU_{1,1}}(t_0) c_{\theta'}(bI_2\times \diag(a, t_0,a)).$$
Up to modulo the center, this recovers the term associated to $E^1$ in $m_{geom}(\theta')$. This finishes the proof of the proposition.
\end{proof}

Let $\Pi_\phi=\Pi_\phi(G)\cup \Pi_\phi(G_D)$ be a tempered $L$-packet of $G/Z_{G,H}$. We assume that $\Pi_\phi(G)$ is not discrete with one element. We first prove that the distinguished element belongs to $\Pi_\phi(G)$ if and only if $\eta_{E/F}(-1)\epsilon(\frac{1}{2},\Pi_\phi,\rho_X)=1$, which is equivalent to the equation
\begin{equation}\label{GU(6) epsilon formula}
\sum_{T\in \CT_{ell}(H_0)}|W(H_0,T)|^{-1} \int_{T(F)/Z_{G,H}(F)}D^H(t) c_{\theta_{\Pi_\phi(G)}}(t)dt
\end{equation}
$$=\frac{\eta_{E/F}(-1)\epsilon(\frac{1}{2},\Pi_\phi,\rho_X)-1}{2}.$$

There are two cases. The first case is when the packet is induced from a maximal Levi subgroup $M$ of $G$. If $M$ is the Siegel parabolic subgroup, we have
$$\sum_{T\in \CT_{ell}(H_0)}|W(H_0,T)|^{-1} \int_{T(F)/Z_{G,H}(F)}D^H(t) c_{\theta_{\Pi_\phi(G)}}(t)dt$$
$$=0=\frac{\eta_{E/F}(-1)\epsilon(\frac{1}{2},\Pi_\phi,\rho_X)-1}{2}.$$
Note that in this case we can decompose the representation $\rho_X\circ\phi$ into the direct sum of a 10-dimensional representation with its dual such that the determinant of the 10-dimensional representation is equal to $\eta_{E/F}(-1)$ at $-1$. This implies that $\eta_{E/F}(-1)\epsilon(\frac{1}{2},\Pi_\phi,\rho_X)=1$.

If $M\simeq \GU_{2,2}\times \GL_1(E)$ (resp. $M\simeq \GU_{1,1}\times \GL_2(E)$), \eqref{GU(6) epsilon formula} follows from Proposition \ref{parabolic GU(6)} and Lemma \ref{lemma-for-GU(4)-GU(2)}~(1) (resp. Theorem D of \cite{P92}).

Next we consider the case when the packet is discrete. By our assumption the packet $\Pi_\phi(G)$ contains more than one element. Hence there exists a proper elliptic extended endoscopic triple $(G',s',{}^L\eta)$ of $G$ such that $\phi$ factors through ${}^L\eta$, $G'=G(U_{1,1}\times U_{2,2})$ and $s'\in Z_\phi$. The $L$-parameter $\phi$ of $G$ induces an $L$-parameter (still denoted by $\phi$) of the endoscopic group $G'$. As in Section \ref{sec:pre}, we can decompose $\rho_X\circ \phi$ as $\rho_{1,\phi,s'}\oplus \rho_{2,\phi,s'}$ where $\dim(\rho_{1,\phi,s'})=12$ and $\dim(\rho_{2,\phi,s'})=8$. By the endoscopic relation in Proposition \ref{prop GU(6)} and Lemma \ref{lemma-for-GU(4)-GU(2)} (note that since we have assumed that Conjecture \ref{weak conjecture GU(4)} holds, the three identities in Lemma \ref{lemma-for-GU(4)-GU(2)} hold for all tempered $L$-packets), we have
$$\sum_{\pi\in \Pi_\phi(G)} \chi_\pi(s')\cdot(c_\pi(1)+\sum_{T\in \CT_{ell}(H_0)}|W(H_0,T)|^{-1} \int_{T(F)/Z_{G,H}(F)}D^H(t) c_{\pi}(t)dt)$$
$$=c_{\theta_{\Pi_\phi(G')}}(1)+\sum_{T\in \CT_{ell}(\GU_{1,1})} \frac{1}{2}\int_{T^\ast(F)} D^{\GU_{1,1}}(t)^2 c_{\theta_{\Pi_\phi(G')}}(\nu'(t))dt$$
$$+\int_{E^1} D^{\GU_{1,1}}(t_a)c_{\theta_{\Pi_\phi(G')}}(\nu(a))da$$
$$=1+\frac{\chi_{\phi,s'}(-1)\eta_{E/F}(-1)\epsilon(\frac{1}{2},\rho_{1,\phi,s'})-1}{2}+\frac{\chi_{\phi,s'}(-1)\epsilon(\frac{1}{2},\rho_{2,\phi,s'})-1}{2}$$
$$=\frac{\chi_{\phi,s'}(-1)\eta_{E/F}(-1)\epsilon(\frac{1}{2},\rho_{1,\phi,s'})+\chi_{\phi,s'}(-1)\epsilon(\frac{1}{2},\rho_{2,\phi,s'})}{2}.$$
Here $\chi_{\phi,s'}$ was defined in Section \ref{sec epsilon factor}. Note that by Theorem 8.1 of \cite{La}, the packet $\Pi_\phi(G')$ is the restriction of an $L$-packet $\Pi'$ of $\GU_{2,2}\times \GU_{1,1}$ which allows us to apply Lemma \ref{lemma-for-GU(4)-GU(2)}. In particular, we know that the above summation is nonzero ($\iff$ the distinguished element belongs to $\Pi_\phi(G)$) if and only if $\eta_{E/F}(-1)\epsilon(\frac{1}{2},\Pi_\phi,\rho_X)=1$. This proves \eqref{GU(6) epsilon formula}.

Now we prove Theorem \ref{main theorem}. Let $\Pi_\phi=\Pi_\phi(G)\cup \Pi_\phi(G_D)$ be as above and let $\omega_\phi\in \hat{S_\phi}$ correspond to the unique distinguished element in the packet and we also view $\omega_\phi$ as a character of $Z_\phi$. For $s\in S_\phi$, by Lemma \ref{lem extended endoscopic triple}, there exists an elliptic extended endoscopic triple $(G',s',{}^L\eta)$ of $G/Z_{G,H}$ such that $s'\in sZ_{\phi}^{\circ}$ and $\phi$ factors through ${}^L\eta$. We need to show that $\omega_\phi(s')=\omega_{\phi,H}(s)$. 

The above discussion implies that $\omega_\phi(s')=\omega_{\phi,H}(s)$ if $s'$ belongs to the center. It remains to consider the case when $s'$ does not belong to the center. We only consider the case when the $-1$ eigenspace of $s'$ is 4 dimensional. The argument for the case when the $-1$ eigenspace of $s'$ is 2 dimensional is similar (note that $s'\in \widehat{G/Z_{G,H}}=\SL_6(\BC)$). In this case the $L$-parameter $\phi$ induces an $L$-parameter of $G''=G(U_{1,1}\times U_{2,2})$ (still denoted by $\phi$). We still decompose $\rho_X\circ \phi$ as $\rho_{1,\phi,s'}\oplus \rho_{2,\phi,s'}$ where $\dim(\rho_{1,\phi,s'})=12$ and $\dim(\rho_{2,\phi,s'})=8$. By our discussion above, we know that 
$$\sum_{\pi\in \Pi_\phi(G)}\chi_\pi(s')m(\pi)$$
is equal to
$$\frac{\chi_{\phi,s'}(-1)\eta_{E/F}(-1)\epsilon(\frac{1}{2},\rho_{1,\phi,s'})+\chi_{\phi,s'}(-1)\epsilon(\frac{1}{2},\rho_{2,\phi,s'})}{2}.$$
Similarly, we can also show that 
$$\sum_{\pi_D\in \Pi_\phi(G_D)} \chi_{\pi_D}(s')m(\pi_D)$$
is equal to 
$$\frac{-\chi_{\phi,s'}(-1)\eta_{E/F}(-1)\epsilon(\frac{1}{2},\rho_{1,\phi,s'})+\chi_{\phi,s'}(-1)\epsilon(\frac{1}{2},\rho_{2,\phi,s'})}{2}.$$

If the distinguished element belongs to $\Pi_\phi(G)$, then 
$$\eta_{E/F}(-1)\epsilon(\frac{1}{2},\Pi_\phi,\rho_X)=1$$ 
and we have
\begin{eqnarray*}
\omega_\phi(s')&=&\sum_{\pi\in \Pi_\phi(G)}\chi_\pi(s')m(\pi)\\
&=&\frac{\chi_{\phi,s'}(-1)\eta_{E/F}(-1)\epsilon(\frac{1}{2},\rho_{1,\phi,s'})+\chi_{\phi,s'}(-1)\epsilon(\frac{1}{2},\rho_{2,\phi,s'})}{2}\\
&=&\chi_{\phi,s'}(-1)\epsilon(\frac{1}{2},\rho_{2,\phi,s'})=\omega_{\phi,H}(s).
\end{eqnarray*}

If the distinguished element belongs to $\Pi_\phi(G_D)$, then 
$$\eta_{E/F}(-1)\epsilon(\frac{1}{2},\Pi_\phi,\rho_X)=-1$$ 
and we have
\begin{eqnarray*}
\omega_\phi(s')&=&\sum_{\pi_D\in \Pi_\phi(G_D)} \chi_{\pi_D}(s')m(\pi_D)\\
&=&\frac{-\chi_{\phi,s'}(-1)\eta_{E/F}(-1)\epsilon(\frac{1}{2},\rho_{1,\phi,s'})+\chi_{\phi,s'}(-1)\epsilon(\frac{1}{2},\rho_{2,\phi,s'})}{2}\\
&=&\chi_{\phi,s'}(-1)\epsilon(\frac{1}{2},\rho_{2,\phi,s'})=\omega_{\phi,H}(s).
\end{eqnarray*}
This completes the proof of Theorem \ref{main theorem} for the model $(\GU_6,\GU_2\ltimes U)$.

Lastly, we prove Theorem \ref{thm weak conjecture smaller models} for the model $(\GU_6,\GU_2\ltimes U)$. Assume that Conjecture \ref{weak conjecture} holds for the model $(\GU_6,\GU_2\ltimes U)$, the goal is to prove Conjecture \ref{weak conjecture GU(4)xGU(2)} for the smaller model $(\GU_4\times \GU_2,\GU_2\ltimes U)$. Let $\Pi_\phi(\GU_4\times \GU_2)$ be a tempered L-packet whose central character is trivial on $\{(aI_4,aI_2)|\;a\in E^{\times}\}$. By restriction it induces a tempered L-packet $\Pi_\phi(G'')$ of $G''=G(U_{1,1}\times U_{2,2})$ with trivial central character and hence a tempered L-packet of $\Pi_\phi(G')$ with $G'=G''/Z_{G''}$. By the endoscopic transfer this induces a L-packet $\Pi_\phi(G/Z_{G,H})$ of $G/Z_{G,H}$. Let $(G',s',{}^L\eta)$ be the elliptic extended endoscopic triple as above. The endoscopic relation in Proposition \ref{prop GU(6)} implies that 
$$\sum_{\pi\in \Pi_\phi(G)} \chi_\pi(s')m(\pi)$$ 
is equal to
\begin{eqnarray*}
&&c_{\theta_{\Pi_\phi(G')}}(1)+\sum_{T\in \CT_{ell}(\GU_{1,1})} \frac{1}{2}\int_{T^\ast(F)} D^{\GU_{1,1}}(t)^2 c_{\theta_{\Pi_\phi(G')}}(\nu'(t))dt\\
&&+\int_{E^1} D^{\GU_{1,1}}(t_a)c_{\theta_{\Pi_\phi(G')}}(\nu(a))da\\
&=&1+\sum_{T\in \CT_{ell}(\GU_{1,1})} \frac{1}{2}\int_{T^\ast(F)} D^{\GU_{1,1}}(t)^2 c_{\theta_{\Pi_\phi(G')}}(\nu'(t))dt\\
&&+\frac{\chi_{\phi}(-1)\epsilon(\frac{1}{2},\Pi_\phi,\rho_{2})-1}{2}.
\end{eqnarray*}
By Conjecture \ref{weak conjecture} together with the fact that $\chi_\pi\in \{\pm 1\}$ for all $\pi\in \Pi_\phi(G)$, we know that the above expression is equal to 
$$\pm \frac{\eta_{E/F}(-1)\epsilon(\frac{1}{2},\Pi_\phi,\rho_X)+1}{2}.$$ By our discussion in the previous subsection, we know that 
$$\sum_{T\in \CT_{ell}(\GU_{1,1})} \frac{1}{2}\int_{T^\ast(F)} D^{\GU_{1,1}}(t)^2 c_{\theta_{\Pi_\phi(G')}}(\nu'(t))dt\in \{0,-1\}$$
and our goal is to show that it is equal to 
$$\frac{\chi_\phi(-1)\eta_{E/F}(-1)\epsilon(\frac{1}{2},\Pi_\phi,\rho_1)-1}{2},$$ 
which is equivalent to show that 
$$\sum_{T\in \CT_{ell}(\GU_{1,1})} \frac{1}{2}\int_{T^\ast(F)} D^{\GU_{1,1}}(t)^2 c_{\theta_{\Pi_\phi(G')}}(\nu'(t))dt$$
is nonzero if and only if 
$$\chi_\phi(-1)\eta_{E/F}(-1)\epsilon(\frac{1}{2},\Pi_\phi,\rho_1)=-1.$$

If 
$$1+\sum_{T\in \CT_{ell}(\GU_{1,1})} \frac{1}{2}\int_{T^\ast(F)} D^{\GU_{1,1}}(t)^2 c_{\theta_{\Pi_\phi(G')}}(\nu'(t))dt$$
$$+\frac{\chi_{\phi}(-1)\epsilon(\frac{1}{2},\Pi_\phi,\rho_{2})-1}{2}=\frac{\eta_{E/F}(-1)\epsilon(\frac{1}{2},\Pi_\phi,\rho_X)+1}{2},$$ 
we have 
$$\sum_{T\in \CT_{ell}(\GU_{1,1})} \frac{1}{2}\int_{T^\ast(F)} D^{\GU_{1,1}}(t)^2 c_{\theta_{\Pi_\phi(G')}}(\nu'(t))dt$$
$$=\frac{\eta_{E/F}(-1)\epsilon(\frac{1}{2},\Pi_\phi,\rho_X)-\chi_{\phi}(-1)\epsilon(\frac{1}{2},\Pi_\phi,\rho_{2})}{2}.$$
Since $\eta_{E/F}(-1)\epsilon(\frac{1}{2},\Pi_\phi,\rho_X),\chi_{\phi}(-1)\epsilon(\frac{1}{2},\Pi_\phi,\rho_{2})\in \{\pm 1\}$ and 
$$\eta_{E/F}(-1)\epsilon(\frac{1}{2},\Pi_\phi,\rho_X)=\chi_{\phi}(-1)\epsilon(\frac{1}{2},\Pi_\phi,\rho_{2})\cdot \chi_\phi(-1)\eta_{E/F}(-1)\epsilon(\frac{1}{2},\Pi_\phi,\rho_{1}),$$
we know that 
$$\sum_{T\in \CT_{ell}(\GU_{1,1})} \frac{1}{2}\int_{T^\ast(F)} D^{\GU_{1,1}}(t)^2 c_{\theta_{\Pi_\phi(G')}}(\nu'(t))dt$$
is nonzero if and only if $\chi_\phi(-1)\eta_{E/F}(-1)\epsilon(\frac{1}{2},\Pi_\phi,\rho_1)=-1$.

If 
$$1+\sum_{T\in \CT_{ell}(\GU_{1,1})} \frac{1}{2}\int_{T^\ast(F)} D^{\GU_{1,1}}(t)^2 c_{\theta_{\Pi_\phi(G')}}(\nu'(t))dt$$
$$+\frac{\chi_{\phi}(-1)\epsilon(\frac{1}{2},\Pi_\phi,\rho_{2})-1}{2}=-\frac{\eta_{E/F}(-1)\epsilon(\frac{1}{2},\Pi_\phi,\rho_X)+1}{2},$$ 
we have 
$$\sum_{T\in \CT_{ell}(\GU_{1,1})} \frac{1}{2}\int_{T^\ast(F)} D^{\GU_{1,1}}(t)^2 c_{\theta_{\Pi_\phi(G')}}(\nu'(t))dt$$
$$=\frac{-\eta_{E/F}(-1)\epsilon(\frac{1}{2},\Pi_\phi,\rho_X)-\chi_{\phi}(-1)\epsilon(\frac{1}{2},\Pi_\phi,\rho_{2})}{2}-1.$$
Since the left hand side is either 0 or -1, we have
$$\eta_{E/F}(-1)\epsilon(\frac{1}{2},\Pi_\phi,\rho_X),\chi_{\phi}(-1)\epsilon(\frac{1}{2},\Pi_\phi,\rho_{2})\in \{\pm 1\}$$ 
and 
$$\eta_{E/F}(-1)\epsilon(\frac{1}{2},\Pi_\phi,\rho_X)=\chi_{\phi}(-1)\epsilon(\frac{1}{2},\Pi_\phi,\rho_{2})\cdot \chi_\phi(-1)\eta_{E/F}(-1)\epsilon(\frac{1}{2},\Pi_\phi,\rho_{1}).$$
This implies that
$$\sum_{T\in \CT_{ell}(\GU_{1,1})} \frac{1}{2}\int_{T^\ast(F)} D^{\GU_{1,1}}(t)^2 c_{\theta_{\Pi_\phi(G')}}(\nu'(t))dt$$
is nonzero if and only if $\chi_\phi(-1)\eta_{E/F}(-1)\epsilon(\frac{1}{2},\Pi_\phi,\rho_1)=-1$. This proves Conjecture \ref{weak conjecture GU(4)xGU(2)}.

\section{The models $(\GSO_8\times \GL_2,\GL_2\ltimes U)$, $(\GSO_{12},\GL_2\ltimes U)$, $(\GSp_6\times \GL_2,\GL_2\ltimes U)$ and $(\GSp_{10},\GL_2\ltimes U)$}\label{sec GSO and GSp}
In this section, we consider the models $(\GSO_8\times \GL_2,\GL_2\ltimes U)$, $(\GSO_{12},\GL_2\ltimes U)$, $(\GSp_6\times \GL_2,\GL_2\ltimes U)$ and $(\GSp_{10},\GL_2\ltimes U)$. In Section \ref{sec GSO GSp 1}, we will define the models and the multiplicity formulas.  In Sections \ref{sec GSO} and \ref{sec GSp 1}, we will prove the main theorem for these models. The proofs for these four models are very similar to each other.

Each of the two models associated to even special orthogonal groups has two versions (corresponding to the two Siegel parabolic subgroups) and they are differed by the outer automorphism of even special orthogonal groups. We will only consider one of them, the other one can be studied by the same argument.

\subsection{The models and the multiplicity formulas}\label{sec GSO GSp 1}

We start with the model $(\GSO_8\times \GL_2,\GL_2\ltimes U)$. Let $$J_2'=\begin{pmatrix}0&-1\\1&0\end{pmatrix},\; J_{2n}'=\begin{pmatrix}0&J_{2n-2}'\\J_2&0 \end{pmatrix},\; L_{4}=\begin{pmatrix}0&J_2'\\-J_2'&0 \end{pmatrix}$$ and
$L_{4n}=\begin{pmatrix}0&0&J_2'\\0&L_{4n-4}&0\\-J_2'&0&0 \end{pmatrix}$. Define
$$\GSO_{4n}=\{g\in \GL_{4n} \mid {}^tgL_{4n}g =l(g)L_{4n},\;\det(g)=l(g)^{2n}\},$$
$$\GSO_{2n}(D)=\{g\in \GL_{2n}(D) \mid {}^t\bar{g} J_{2n}'g=l(g)J_{2n}'\}.$$

Let $G=\GSO_{8}\times \GL_2$, $H=H_0\ltimes U$ with 
$$H_0=\{\diag(h,h,h,h)\times h \mid h\in \GL_2\}$$
and $U$ be the unipotent radical of the standard parabolic subgroup $P=LU$ of $G=\GSO_8\times \GL_2$ where 
$$L=\{\diag(h_1,h_2, th_{2}^{\ast}, th_{1}^{\ast})\times h_3 \mid h_i\in \GL_2,t\in \GL_1\},\;h^\ast=J_2'{}^th^{-1} (J_2')^{-1}.$$ 
We define a generic character $\xi$ on $U(F)$ to be $\xi(u)=\psi(\lambda(u))$ where 
$$\lambda(u)=\tr(X)+\tr(Y),\;u=\begin{pmatrix}I_2&X&\ast &\ast  \\ 0&I_2&Y&\ast\\ 0&0&I_2& \ast \\ 0&0&0&I_2\end{pmatrix}.$$ 
Similarly we can also define the quaternion algebra version of this model $(G_D,H_D)$ with $G_D=\GSO_{4}(D)\times \GL_1(D)$ and $H_{0,D}=\GL_1(D)$.

For the model $(\GSO_{12},\GL_2\ltimes U)$, let $G=\GSO_{12}$, $H=H_0\ltimes U$ with 
$$H_0=\{\diag(h,h,h,h,h,h) \mid h\in \GL_2\}$$
and $U$ is the unipotent radical of the standard parabolic subgroup $P=LU$ of $G$ where 
$$L=\{\diag(h_1,h_2, h_3,t h_{3}^{\ast},th_{2}^{\ast}, th_{1}^{\ast}) \mid h_i\in \GL_2,t\in \GL_1\}.$$ 
We define a generic character $\xi$ on $U(F)$ to be $\xi(u)=\psi(\lambda(u))$ where 
$$\lambda(u)=\tr(X)+\tr(Y)+\tr(Z),\;u=\begin{pmatrix}I_2&X&\ast &\ast &\ast&\ast \\ 0&I_2&Y&\ast&\ast &\ast\\ 0&0&I_2& Z & \ast&\ast \\ 0&0&0&I_2 &\ast&\ast \\ 0&0&0&0&I_2&\ast  \\ 0&0&0&0&0&I_2\end{pmatrix}.$$  
Similarly we can also define the quaternion algebra version of this model $(G_D,H_D)$ with $G_D=\GSO_{6}(D)$ and $H_{0,D}=\GL_1(D)$.

For the model $(\GSp_6\times \GL_2,\GL_2\ltimes U)$, define 
$$\GSp_{2n}=\{g\in \GL_{2n} \mid {}^t gJ_{2n}'g =l(g)J_{2n}'\},$$ 
$$\GSp_{n}(D)=\{g\in \GL_n(D)\mid {}^t\bar{g}w_ng=l(g)w_n\}
$$
where $w_n$ is the longest Weyl element of $\GL_n$. Let $G=\GSp_{6}\times \GL_2$, $H=H_0\ltimes U$ with 
$$H_0=\{\diag(h,h,h)\times h \mid h\in \GL_2\}$$
and the unipotent radical $U$ of the standard parabolic subgroup $P=LU$ of $G=\GSp_6\times \GL_2$ where 
$$L=\{(h_1,h_2, \det(h_2)h_{1}^{\ast})\times h_3 \mid h_i\in \GL_2\}.$$ 
We define a generic character $\xi$ on $U(F)$ to be $\xi(u)=\psi(\lambda(u))$ where 
$$\lambda(u)=\tr(X),\;u=\begin{pmatrix}I_2&X&\ast  \\ 0&I_2&\ast \\ 0&0&I_2 \end{pmatrix}.$$ 
Similarly, we can also define the quaternion algebra version of this model $(G_D,H_D)$ with $G_D=\GSp_{3}(D)\times \GL_1(D)$ and $H_{0,D}=\GL_1(D)$.

For the model $(\GSp_{10},\GL_2\ltimes U)$, let $G=\GSp_{10}$, $H=H_0\ltimes U$ with 
$$H_0=\{\diag(h,h,h,h,h) \mid h\in \GL_2,\;h^\ast=J_2'{}^th^{-1} (J_2')^{-1}\}$$
and   the unipotent radical $U$ of the standard parabolic subgroup $P=LU$ of $G$ where 
$$L=\{\diag(h_1,h_2, h_3,\det(h_3) h_{2}^{\ast}, \det(h_3)h_{1}^{\ast}) \mid h_i\in \GL_2\}.$$ 
We define a generic character $\xi$ on $U(F)$ to be $\xi(u)=\psi(\lambda(u))$ where 
$$\lambda(u)=\tr(X)+\tr(Y),\;u=\begin{pmatrix}I_2&X&\ast &\ast &\ast \\ 0&I_2&Y&\ast&\ast \\ 0&0&I_2& \ast & \ast \\ 0&0&0&I_2 &\ast \\ 0&0&0&0&I_2 \end{pmatrix}.$$ 
Similarly, we can also define the quaternion algebra version of this model $(G_D,H_D)$ with $G_D=\GSp_{5}(D)$ and $H_{0,D}=\GL_1(D)$.

Let $(G,H)$ be one of the four models above. For a quasi-character $\theta$ (resp. $\theta_D$) of $G(F)$ (resp. $G_D(F)$), define the geometric multiplicities
$$m_{geom}(\theta)=c_\theta(1)+\sum_{T\in \CT_{ell}(H_0)} \frac{1}{2}\int_{T(F)/Z_{G,H}(F)}D^H(t) c_{\theta}(t)dt,$$
$$m_{geom}(\theta_D)=\sum_{T_D\in \CT_{ell}(H_{0,D})}\frac{1}{2} \int_{T_D(F)/Z_{G_D,H_D}(F)} D^{H_D}(t) c_{\theta_D}(t)dt.$$
In our previous paper \cite{WZ2}, we have proved the multiplicity formulas
$$m(\pi)=m_{geom}(\theta_\pi),\;m(\pi_D)=m_{geom}(\theta_{\pi_D})$$
for all tempered representations $\pi$ (resp. $\pi_D$) of $G(F)$ (resp. $G_D(F)$) in the $p$-adic case. For the rest of this section we will assume that the multiplicity formulas hold for both the $p$-adic case and the real case.  

To end this subsection, we will discuss the behavior of the geometric multiplicities under parabolic induction. 
Let $M$ be a proper Levi subgroup of $G$ and $\theta^M$ be a quasi-character on $M(F)$. If $M$ does not contain the Levi subgroup $L$ up to conjugation, define $m_{geom}(\theta^M)=c_{\theta^M}(1)$. Otherwise, $M$ corresponds to a proper Levi subgroup $M_D$ of $G_D$. Moreover, up to conjugation we may assume that $L\subset M$ and $L_D\subset M_D$. Let $\theta_{D}^{M_D}$ be a quasi-character on $M_D(F)$. Define
\begin{eqnarray*}
m_{geom}(\theta^M)&=&c_{\theta^M}(1)+\sum_{T\in \CT_{ell}(H_0)}\frac{1}{2}\int_{T(F)/Z_{G,H}(F)}\\
&&D^M(t)^{1/2}(t)D^{H_0}(t)^{-1/2} c_{\theta^M}(t)dt,\\
m_{geom}(\theta_{D}^{M_D})&=&\sum_{T_D\in \CT_{ell}(H_{0,D})} \frac{1}{2}\int_{T_D(F)/Z_{G_D,H_D}(F)}\\
&& D^{M_D}(t)^{1/2}D^{H_{0,D}}(t)^{-1/2} c_{\theta_{D}^{M_D}}(t)dt.
\end{eqnarray*}
The following proposition is a direct consequence of Proposition \ref{germ parabolic induction} (one just needs to use the fact that $D^H(t)=D^G(t)^{1/2}D^{H_0}(t)^{-1/2}$ for $t\in H_0(F)$). 

\begin{prop}\label{GSO GSp parabolic induction}
Let $\theta$ (resp. $\theta_D$) be a quasi-character on $G(F)$ (resp. $G_D(F)$). Assume that $\theta$ (resp. $\theta_D$) is the parabolic induction of a quasi-character $\theta^M$ (resp. $\theta_{D}^{M_D}$) of a proper Levi subgroup $M$ of $G$ (resp. $M_D$ of $G_D$). We have
$$m_{geom}(\theta)=m_{geom}(\theta^M),\;m_{geom}(\theta_D)=m_{geom}(\theta_{D}^{M_D}).$$
\end{prop}

\subsection{The proof of Theorem \ref{main theorem} and \ref{thm weak conjecture smaller models} for $(\GSO_{12},\GL_2\ltimes U)$ and $(\GSO_{8}\times\GL_2,\GL_2\ltimes U)$}\label{sec GSO}

In this subsection, we will prove Theorem \ref{main theorem} for $(\GSO_{12},\GL_2\ltimes U)$ and $(\GSO_{8}\times\GL_2,\GL_2\ltimes U)$. The proof for the model $(\GSO_{8}\times\GL_2,\GL_2\ltimes U)$ is very similar (and easier) to the proof for the model $(\GSO_{12},\GL_2\ltimes U)$. So we will only consider the model $(\GSO_{12},\GL_2\ltimes U)$.

Let $(G,H)$ be the model $(\GSO_{12},\GL_2\ltimes U)$ defined in the previous subsection. We will first study the behaviors of the geometric multiplicities under endoscopic transfer. Then we will prove Theorem \ref{main theorem}.

Let $(G',s',{}^L\eta)$ be a proper extended elliptic endoscopic triple of $G/Z_{G,H}$. The projection of $s'\in \widehat{G/Z_{G,H}}=\Spin_{12}(\BC)$ to $\SO_{12}(\BC)$ is conjugated to 
$$\diag(I_6,-I_6),\;\diag(I_8,-I_4)\;\;\;\text{or}\;\;\;\diag(I_4,-I_8).$$ 
We will only consider the first two cases  since the third case is very similar to the second case.

If $G'=G(\SO_6\times \SO_6)/\GL_{1}^{\diag}$, we just let 
$$m_{geom}(\theta')=c_{\theta'}(1),\;\;m_{geom,D}(\theta')=0.$$ 

If $G'=G(\SO_8\times \SO_4)/\GL_{1}^{\diag}$, as explained in Section \ref{sec epsilon factor}, when we restrict the representation $\rho_X$ to $\hat{G}'=(\widehat{G/Z_{G,H}})_{s'}$,  we can decompose it as $$\rho_X=\rho_{s',+}\oplus \rho_{s',-}$$ 
where $\rho_{s',+}$ (resp. $\rho_{s',-}$) is the tensor product of a Half-Spin representation of $\Spin_8(\BC)$ with a Half-Spin representation of $\Spin_4(\BC)$ and it is the $+1$ (resp. $-1$) eigenspace of $\rho_X(s')$. 

Up to conjugation the group $G'$ has 4 Levi subgroups that are isomorphic to $\GL_2\times \GL_2\times \GL_2\times \GL_1/\GL_{1}^{\diag}$. There are exactly 2 of these 4 Levi subgroups whose elliptic conjugacy classes correspond to the elliptic conjugacy classes of $L/Z_{G,H}$ (the elliptic conjugacy classes of the other 2 Siegel Levi subgroups correspond to the elliptic conjugacy classes of $\sigma(L)$ where $\sigma$ is the outer automorphism of $G$), we denote them by $L_1,L_2$. By switching $L_1$ and $L_2$ we may assume that the $\GSO_8$-component of $L_1$ (resp. $L_2$) corresponds to the Half-Spin representation of $\Spin_8(\BC)$ appeared in $\rho_{s',+}$ (resp. $\rho_{s',-}$). 

\begin{rmk}
For each Levi subgroup of $\GSO_8(\BC)$ that is isomorphic to $\GL_2\times \GL_2\times \GL_1$, we can construct the model $(\GSO_8\times \GL_2,\GL_2\ltimes U)$ as in the previous subsection. Hence it corresponds to a Half-Spin representation of $\Spin_8(\BC)$. Up to conjugation there are two such Levi subgroups differed by the outer automorphisms.
\end{rmk}

We can embed $H_0/Z_{G,H}\simeq \PGL_2$ diagonally into $L_i$ (denoted by $\nu_i$). We define
\begin{eqnarray*}
m_{geom}(\theta')&=&c_{\theta'}(1)+\sum_{T\in \CT_{ell}(H_0)}\frac{1}{2}\int_{T(F)/Z_{G,H}(F)} \\
&&\sum_{i=1}^{2} D^{H_0}(t)^3c_{\theta'}(\nu_i(t))dt,\\
m_{geom,D}(\theta')&=&\sum_{T\in \CT_{ell}(H_0)}\frac{1}{2} \int_{T(F)/Z_{G,H}(F)} \\
&&\sum_{i=1}^{2} (-1)^{i-1}D^{H_0}(t)^3c_{\theta'}(\nu_i(t))dt.
\end{eqnarray*}

\begin{prop}\label{prop GSO(12)}
Let $\theta$ (resp. $\theta_D$) be a quasi-character on $G/Z_{G,H}(F)$. Assume that $\theta$ (resp. $\theta_D$) is the endoscopic transfer of a quasi-character $\theta'$ of $G'(F)$ . We have
$$m_{geom}(\theta)=m_{geom}(\theta'),\;m_{geom}(\theta_D)=m_{geom,D}(\theta').$$
\end{prop}

\begin{proof}
We first consider the split case, the proof of $c_\theta(1)=c_{\theta'}(1)$ is straightforward and we will skip it here. We only need to study the term $D^H(t)c_\theta(t)$ for $t=\diag(t_0,t_0,t_0,t_0,t_0,t_0)$ and $t_0\in \GL_2(F)_{ell,reg}$. Let $E/F$ be the quadratic extension associated to $t_0$ and we can identify $t_0$ with an element of $E^\times$. By Proposition \ref{regular germs}, we know that 
$$D^H(t)c_\theta(t)$$ 
is equal to 
$$\frac{1}{48}D^{\GL_2}(t_0)^{-1/2}\cdot \lim_{\lambda_j\in F^\times \rightarrow 1} D^G(g(t_0,\lambda_1,\lambda_2,\lambda_3))^{1/2}\theta(g(t_0,\lambda_1,\lambda_2,\lambda_3))$$
where $g(t_0,\lambda_1,\lambda_2,\lambda_3)$ is the conjugacy class of $G$ corresponding to ($\bar{t}_0$ is the conjugation of $t_0$ under the nontrivial element in $Gal(E/F)$)
$$(E\oplus E,E,(\lambda_1t_0,\lambda_{1}^{-1}\bar{t}_0))\cup (E\oplus E,E,(\lambda_2t_0,\lambda_{2}^{-1}\bar{t}_0))$$
$$\cup (E\oplus E,E,(\lambda_3t_0,\lambda_{3}^{-1}\bar{t}_0))$$
and has a nonempty intersection with $L$ (recall that in the even special orthogonal group case each data above gives two conjugacy classes differed by the outer automorphism). Here $48$ is the cardinality of the Weyl group of the centralizer of $t_0$, which is of type $C_3$. Note that 
$$\Pi_{i\in I} F_{\pm i}^{\times}/Im(N_{F_{i}/F_{\pm i}})/\sim$$ 
is trivial for this conjugacy class.

When $G'=G(\SO_6\times \SO_6)/\GL_{1}^{diag}$, there is no conjugacy class in $G'$ corresponding to the conjugacy class $g(t_0,\lambda_1,\lambda_2,\lambda_3)$ of $G$. This implies that $c_{\theta}(t)=0$. In particular, we have
$$m_{geom}(\theta)=c_\theta(1)=c_{\theta'}(1)=m_{geom}(\theta').$$

When $G'=G(\SO_8\times \SO_4)/\GL_{1}^{diag}$, there are six conjugacy classes in $G'$ corresponding to the conjugacy class $g(t_0,\lambda_1,\lambda_2,\lambda_3)$ of $G$:
$$((E\oplus E,E,(\lambda_1t_0,\lambda_{1}^{-1}\bar{t}_0))\cup (E\oplus E,E,(\lambda_2t_0,\lambda_{2}^{-1}\bar{t}_0)))$$
$$\times (E\oplus E,E,(\lambda_3t_0,\lambda_{3}^{-1}\bar{t}_0)),$$
$$((E\oplus E,E,(\lambda_1t_0,\lambda_{1}^{-1}\bar{t}_0))\cup (E\oplus E,E,(\lambda_3t_0,\lambda_{3}^{-1}\bar{t}_0)))$$
$$\times (E\oplus E,E,(\lambda_2t_0,\lambda_{2}^{-1}\bar{t}_0)),$$
$$((E\oplus E,E,(\lambda_2t_0,\lambda_{2}^{-1}\bar{t}_0))\cup (E\oplus E,E,(\lambda_3t_0,\lambda_{3}^{-1}\bar{t}_0)))$$
$$\times (E\oplus E,E,(\lambda_1t_0,\lambda_{1}^{-1}\bar{t}_0)).$$

Note that each data above gives 4 conjugacy classes of $G'$ and two of them correspond to the conjugacy class $g(t_0,\lambda_1,\lambda_2,\lambda_3)$ of $G$. The other two correspond to the image of $g(t_0,\lambda_1,\lambda_2,\lambda_3)$ under the outer automorphism. For the two conjugacy classes corresponding to $g(t_0,\lambda_1,\lambda_2,\lambda_3)$, one of them has a nonempty intersection with the Levi subgroup $L_1$ and the other one has a nonempty intersection with $L_2$. The transfer factor is trivial in this case since the character $\eta_{E\oplus E/E}$ is trivial.

Let $g_i'(t_0,\lambda_1,\lambda_2,\lambda_3)$ ($1\leq i\leq 2$) be the two conjugacy classes of $G'$ associated to one of the three data above that corresponds to the conjugacy class $g(t_0,\lambda_1,\lambda_2,\lambda_3)$ of $G$ and has a nonempty intersection with the Levi subgroup $L_i$. By Proposition \ref{regular germs}, we have
\begin{eqnarray*}
&&\lim_{\lambda_j\in F^\times \rightarrow 1} D^{G'}(g_j'(t_0,\lambda_1,\lambda_2,\lambda_3))^{1/2}\theta'(g_j'(t_0,\lambda_1,\lambda_2,\lambda_3))\\
&=&16\cdot D^{G'}(\nu_i(t_0))^{1/2}c_{\theta'}(\nu_i(t_0))
\end{eqnarray*}
where $16$ comes from the cardinality of the Weyl group, which is of type $C_2\times A_1$.
Then the proposition follows from the fact that 
$D^{G'}(\nu_i(t_0))^{1/2}D^{\GL_2}(t_0)^{-1/2}=D^{H_0}(t_0)^3.$
This proves the split case.

For the quaternion case, the proof is exactly the same as the split case with one exception. The only difference is that the transfer factor between $g(t_0,\lambda_1,\lambda_2,\lambda_3)$ and $g_i'(t_0,\lambda_1,\lambda_2,\lambda_3)$ is equal to $1$ when $i=1$ and is equal to $-1$ when $i=2$. This difference comes from the extra pairing $\langle inv[z](\delta^\ast,\delta),s^\Fe \rangle$ in the definition of the transfer factor in Section 2.3 of \cite{K}. With the same notation as in loc. cit., for the conjugacy classes we considered here, both $H^1(\Gamma, S)$ and $\pi_0(\hat{S}^\Gamma)$ are isomorphic to $\BZ/2\BZ$.  
Also it is easy to see that $inv[z](\delta^\ast,\delta)$ is the nontrivial element in $H^1(\Gamma, S)$ and the element $s^\Fe$ belongs to the identity component (resp. the non-identity component) of $\hat{S}^\Gamma$ when $i=1$ (resp. $i=2$). In particular the pairing $\langle inv[z](\delta^\ast,\delta),s^\Fe \rangle$ is equal to 1 (resp. $-1$) if $i=1$ (resp. $i=2$). This finishes the proof of the proposition.
\end{proof}

Let $\Pi_\phi=\Pi_\phi(G)\cup \Pi_\phi(G_D)$ be a tempered $L$-packet with trivial central character. Assume that $\Pi_\phi$ is not discrete with $|\Pi_\phi(G)|=1$. We first show that the unique distinguished element belongs to $\Pi_\phi(G)$ if and only if $\epsilon(\frac{1}{2},\Pi_\phi,\rho_X)=1$, which is equivalent to say that
\begin{equation}\label{epsilon formula GSO(12)}
\sum_{T\in \CT_{ell}(H_0)}\frac{1}{2}\int_{T(F)/Z_{G,H}(F)} D^H(t)c_{\theta_{\Pi_\phi(G)}}(t)dt=\frac{\epsilon(\frac{1}{2},\Pi_\phi,\rho_X)-1}{2}.    
\end{equation}

There are two situations. The first situation is when the packet $\Pi_\phi(G)$ is induced from a packet $\Pi_\phi(M)$ of a maximal Levi subgroup $M(F)$ of $G(F)$. If $M=\GL_6\times \GL_1$ (resp. $M=\GSO_{8}\times \GL_2$), \eqref{epsilon formula GSO(12)} follows from Proposition \ref{GSO GSp parabolic induction} and Conjecture \ref{weak conjecture} for the model $(\GL_6,\GL_2\ltimes U)$ (resp. $(\GSO_8\times \GL_2,\GL_2\ltimes U)$). To be specific, in this case, Proposition \ref{GSO GSp parabolic induction} implies that the multiplicity of the packet $\Pi_\phi(G)$ is equal to the multiplicity of the packet $\Pi_\phi(M)$ with respect to the model $(\GL_6,\GL_2\ltimes U)$ (resp. $(\GSO_8\times \GL_2,\GL_2\ltimes U)$). Meanwhile, the epsilon factor $\epsilon(\frac{1}{2},\Pi_\phi,\rho_X)$ is equal to the epsilon factor of the packet $\Pi_\phi(M)$ associated to the model $(\GL_6,\GL_2\ltimes U)$ (resp. $(\GSO_8\times \GL_2,\GL_2\ltimes U)$). This proves \eqref{epsilon formula GSO(12)} (note that both the model $(\GL_6,\GL_2\ltimes U)$ and the model $(\GSO_8\times \GL_2,\GL_2\ltimes U)$ are smaller than the model $(G,H)$).

If $M=\GL_4\times \GSO_4=\GL_4\times (\GL_2\times \GL_2)/\GL_{1}^{anti-diag}$, as in the previous two cases, \eqref{epsilon formula GSO(12)} follows from Proposition \ref{GSO GSp parabolic induction} and Conjecture \ref{weak conjecture} for the model $(\GL_4,\GL_2\times \GL_2)$. By Remark \ref{GL(6) implies GL(4)}, we know that Conjecture \ref{weak conjecture} for the model $(\GL_4,\GL_2\times \GL_2)$ follows from Conjecture \ref{weak conjecture} for the model $(\GL_6,\GL_2\ltimes U)$.

If $M=\GSO_6\times \GL_3$ or $\GSO_{10}\times \GL_1$, then it is easy to see that both sides of \eqref{epsilon formula GSO(12)} are equal to 0 and this proves the equation.

\begin{rmk}\label{SO(8) implies GL(4)}
By the above discussion, we know that if the weak conjecture (Conjecture \ref{weak conjecture}) holds for the model $(\GSO_{12},\GL_2\ltimes U)$, then it also holds for the models $(\GL_6,\GL_2\ltimes U)$, $(\GSO_8\times \GL_2,\GL_2\ltimes U)$,\;$(\GL_4\times \GL_2,\GL_2\times \GL_2)$. Similarly, if the weak conjecture holds for the model $(\GSO_8\times \GL_2,\GL_2\ltimes U)$, then it also holds for the model $(\GL_4\times \GL_2,\GL_2\times \GL_2)$. This proves Theorem \ref{thm weak conjecture smaller models} for the models $(\GSO_{12},\GL_2\ltimes U)$ and $(\GSO_8\times \GL_2,\GL_2\ltimes U)$.
\end{rmk}

The second situation is when the packet $\Pi_\phi(G)$ is discrete. By our assumption, we must have $|\Pi_\phi(G)|>1$. Hence there exists a proper elliptic extended endoscopic triple $(G',s',{}^L\eta)$ of $G$ such that $\phi$ factors through ${}^L\eta$ and $s'\in Z_\phi$. If the order of $s'$ is equal to 2, let $W=W_{s,+}\oplus W_{s,-}$ be the decomposition of the 12-dimensional quadratic space as in Section \ref{sec epsilon factor}. 
Up to multiplying $s'$ by an element in the center, we may assume that $\dim(W_{s',+})=8$. Then the $L$-parameter $\phi$ of $G$ can be viewed as an $L$-parameter (still denoted by $\phi$) of the endoscopic group $G'=G(\SO_8\times \SO_4)/\GL_{1}^{\diag}$. As in Section \ref{sec epsilon factor}, we can decompose $\rho_X\circ\phi$ as $\rho_{s',\phi,+}\oplus \rho_{s',\phi,-}$ where the underlying vector space of $\rho_{s',\phi,+}$ (resp. $\rho_{s',\phi,-}$) is the tensor product of a Half-Spin representation of the even Spin group associated to $W_{s',+}$ with a Half-Spin representation of the even Spin group associated to $W_{s',-}$ and it is the $+1$ (resp. $-1$) eigenspace of $\rho_X(s')$.

By our formula of endoscopy in Proposition \ref{prop GSO(12)}, we know that (note that the Kottwitz sign between $G$ and $G_D$ is $-1$)
$$\sum_{\pi\in \Pi_\phi(G)} \tr(\chi_\pi(s'))m(\pi)=m_{geom}(\theta'),$$
$$\sum_{\pi_D\in \Pi_\phi(G_D)} \tr(\chi_{\pi_D}(s'))m(\pi_D)=-m_{geom,D}(\theta')$$
where $\theta'=\theta_{\Pi_{\phi}(G')}$ is the distribution character of the $L$-packet $\Pi_{\phi}(G')$. Note that in this case $S_\phi$ is not necessarily abelian.

The $L$-packet $\Pi_{\phi}(G')$ induces an $L$-packet of $G(\SO_8\times \SO_4)$. By Theorem 8.1 of \cite{La}, such an $L$-packet is the restriction of an $L$-packet of the group $\GSO_8\times \GSO_4=\GSO_8\times (\GL_2\times \GL_2/\GL_{1}^{anti-diag})$ to $G'$. This gives an $L$-packet of the group $G''=\GSO_8\times \GL_2\times \GL_2$,  denoted by $\Pi_\phi(G'')$. 
We use $\theta''$ to denote the distribution character of this packet. Up to switching the two $\GL_2$ copies we assume that the Half-Spin representation of $\Spin_4(\BC)=\Spin_4(W_{s',-})$ appeared in $\rho_{s',\phi,+}$ (resp. $\rho_{s',\phi,-}$) is the 2-dimensional standard representation of the dual group of the first (resp. second) copy of $\GL_2$.

In this case, the  embedding $\nu_1$ (resp. $\nu_2$) of $\PGL_2$ into $G'$ induces a diagonal embedding of $\GL_2$ into $\GSO_8$ and the first copy (resp. second copy) of $\GL_2$. Note that the restriction of these two embeddings to $\GSO_8$ are differed by the outer automorphism.  Combining with the multiplicity formula and Conjecture \ref{weak conjecture} for the model $(\GSO_8\times \GL_2,\GL_2\ltimes U)$, we know that 
$$m_{geom}(\theta')=\frac{\epsilon(\frac{1}{2},\rho_{s',\phi,+})+\epsilon(\frac{1}{2},\rho_{s',\phi,-})}{2},$$
$$m_{geom,D}(\theta')=\frac{\epsilon(\frac{1}{2},\rho_{s',\phi,+})-\epsilon(\frac{1}{2},\rho_{s',\phi,-})}{2}.$$
In particular, this implies that 
$$\sum_{\pi\in \Pi_\phi(G)} \tr(\chi_\pi(s'))m(\pi)=\frac{\epsilon(\frac{1}{2},\rho_{s',\phi,+})+\epsilon(\frac{1}{2},\rho_{s',\phi,-})}{2},$$

$$\sum_{\pi_D\in \Pi_\phi(G_D)} \tr(\chi_{\pi_D}(s'))m(\pi_D)=\frac{-\epsilon(\frac{1}{2},\rho_{s',\phi,+})+\epsilon(\frac{1}{2},\rho_{s',\phi,-})}{2}.$$
As a result, we know that the unique distinguished element belongs to $\Pi_\phi(G)$ if and only if $\epsilon(\frac{1}{2},\rho_{s',\phi,+})=\epsilon(\frac{1}{2},\rho_{s',\phi,-})$ which is equivalent to $\epsilon(\frac{1}{2},\Pi_\phi,\rho_X)=1$.

If the order of $s'$ is equal to 4, let $W=W_{s',+}\oplus W_{s',-}$ be the decomposition of the 12-dimensional quadratic space as in Section \ref{sec:pre} with $\dim(W_{s',+})=\dim(W_{s',-})=6$. Then the $L$-parameter $\phi$ of $G$ can be viewed as an $L$-parameter (still denoted by $\phi$) of the endoscopic group $G'=G(\SO_6\times \SO_6)/\GL_{1}^{diag}$. Let $\theta'$ be the distribution character of $\Pi_\phi(G')$. Then Proposition \ref{prop GSO(12)} implies that 
$$\sum_{\pi\in \Pi_\phi(G)} \tr(\chi_\pi(s'))m(\pi)=c_{\theta'}(1)=1,\;\sum_{\pi_D\in \Pi_\phi(G_D)} \tr(\chi_{\pi_D}(s'))m(\pi_D)=0,$$
i.e. the unique distinguished element belongs to $\Pi_\phi(G)$. In this case, by our discussion in Section \ref{sec epsilon factor}, we also know that $\epsilon(\frac{1}{2},\Pi_\phi,\rho_X)=1$. This proves \eqref{epsilon formula GSO(12)}.

Now we are ready to prove the theorem. Let $\omega_\phi\in \hat{S}_\phi$ correspond  to the unique distinguished element in the packet. By Remark \ref{distinguished is character} we know that $\omega_{\phi}$ is a character and we view it as a character of $Z_\phi$. For $s\in S_\phi$, by Lemma \ref{lem extended endoscopic triple}, there exists an elliptic extended endoscopic triple $(G',s',{}^L\eta)$ of $G/Z_{G,H}$ such that $s'\in sZ_{\phi}^{\circ}$ and $\phi$ factors through ${}^L\eta$. We need to show that $\omega_\phi(s')=\omega_{\phi,H}(s)$. The above discussion implies that $\omega_\phi(s')=\omega_{\phi,H}(s)$ if $s'$ belongs to the center of the dual group.

If $s'$ does not belong to the center of the dual group, there are two cases. If the order of $s'$ is equal to 4, by the discussion above we know that the unique distinguished element belongs to $\Pi_\phi(G)$. By the definition of $\omega_{\phi,H}$ we know that $\omega_{\phi,H}(s)=1$. This implies that
$$\omega_\phi(s')=\tr(\omega_\phi(s'))=\sum_{\pi\in \Pi_\phi(G)} \tr(\chi_\pi(s'))m(\pi)=1=\omega_{\phi,H}(s).$$

If the order of $s'$ is equal to 2, let $W=W_{s',+}\oplus W_{s',-}$ be the decomposition of the 12 dimensional quadratic space as in Section \ref{sec epsilon factor} with $\dim(W_{s',+})\in \{4,8\}$. We will only consider the case when $\dim(W_{s',+})=8$, the other case follows from a similar argument. By our discussion above, we have 
$$\sum_{\pi\in \Pi_\phi(G)} \tr(\chi_\pi(s'))m(\pi)=\frac{\epsilon(\frac{1}{2},\rho_{s',\phi,+})+\epsilon(\frac{1}{2},\rho_{s',\phi,-})}{2},$$
$$\sum_{\pi_D\in \Pi_\phi(G_D)} \tr(\chi_{\pi_D}(s'))m(\pi_D)=\frac{-\epsilon(\frac{1}{2},\rho_{s',\phi,+})+\epsilon(\frac{1}{2},\rho_{s',\phi,-})}{2}.$$
By the definition in Section \ref{sec epsilon factor}, we have
$$\omega_{\phi,H}(s)=\epsilon(\frac{1}{2},\rho_{s',\phi,-}).$$

We have two cases. If the unique distinguished element belongs to $\Pi_\phi(G)$, we have
$$\epsilon(\frac{1}{2},\Pi_\phi,\rho_X)=1,\;\epsilon(\frac{1}{2},\rho_{s',\phi,+})=\epsilon(\frac{1}{2},\rho_{s',\phi,-}).$$
This implies that
$$\omega_\phi(s')=\tr(\omega_\phi(s'))=\sum_{\pi\in \Pi_\phi(G)} \tr(\chi_\pi(s'))m(\pi)$$
$$=\frac{\epsilon(\frac{1}{2},\rho_{s',\phi,+})+\epsilon(\frac{1}{2},\rho_{s',\phi,-})}{2}=\epsilon(\frac{1}{2},\rho_{s',\phi,-}).$$

If the unique distinguished element belongs to $\Pi_\phi(G_D)$, we have
$$\epsilon(\frac{1}{2},\Pi_\phi,\rho_X)=-1,\;\epsilon(\frac{1}{2},\rho_{s',\phi,+})=-\epsilon(\frac{1}{2},\rho_{s',\phi,-}).$$
This implies that
$$\omega_\phi(s')=\tr(\omega_\phi(s'))=\sum_{\pi_D\in \Pi_\phi(G_D)} \tr(\chi_{\pi_D}(s'))m(\pi_D)$$
$$=\frac{-\epsilon(\frac{1}{2},\rho_{s',\phi,+})+\epsilon(\frac{1}{2},\rho_{s',\phi,-})}{2}=\epsilon(\frac{1}{2},\rho_{s',\phi,-}).$$
This finishes the proof of Theorem \ref{main theorem} for the model $(\GSO_{12},\GL_2\ltimes U)$.

\subsection{The proof of Theorem \ref{main theorem} for $(\GSp_{10},\GL_2\ltimes U)$ and $(\GSp_{6}\times \GL_2,\GL_2\ltimes U)$}\label{sec GSp 1}

In this subsection, we will prove Theorem \ref{main theorem} for $(\GSp_{10},\GL_2\ltimes U)$ and $(\GSp_{6}\times \GL_2,\GL_2\ltimes U)$. The proof for the model $(\GSp_{6}\times \GL_2,\GL_2\ltimes U)$ is very similar (and easier) to the proof for the model $(\GSp_{10},\GL_2\ltimes U)$. So we will only consider the model $(\GSp_{10},\GL_2\ltimes U)$.

Let $(G,H)$ be the model $(\GSp_{10},\GL_2\ltimes U)$ defined in Section \ref{sec GSO GSp 1}. We will first study the behavior of the geometric multiplicities under  endoscopic transfer. Then we will prove Theorem \ref{main theorem}.

Let $(G',s',{}^L\eta)$ be a proper elliptic extended endoscopic triple of $G/Z_{G,H}$. If $G'=\PGSO_{10}$ or $G(\Sp_4\times \SO_6)/\GL_{1}^{\diag}$, we just let 
$$m_{geom}(\theta')=c_{\theta'}(1),\;m_{geom,D}(\theta')=0.$$ 
If $G'=G(\Sp_6\times \SO_4)/\GL_{1}^{\diag}$ (resp. $G'=G(\Sp_2\times \SO_8)/\GL_{1}^{\diag}$), the projection of $s'\in \widehat{G/Z_{G,H}}=\Spin_{11}(\BC)$ to $\SO_{11}(\BC)$  is conjugated to $(I_7,-I_4)$ (resp. $(I_3,-I_8)$). For such $s'$, as explained in Section \ref{sec epsilon factor}, when we restrict the representation $\rho_X$ to $\hat{G}'=(\widehat{G/Z_{G,H}})_{s'}$,  we can decompose it as $$\rho_X=\rho_{s',+}\oplus \rho_{s',-}$$ 
where  $\rho_{s',+}$ (resp. $\rho_{s',-}$) is the tensor product of  the Spin representation of $\Spin_7(\BC)$ (resp. $\Spin_3(\BC)$) with a Half-Spin representation of $\Spin_4(\BC)$ (resp. $\Spin_8(\BC)$). Moreover the underlying vector space of $\rho_{s',+}$ (resp. $\rho_{s',-}$) is the $+1$ (resp. $-1$) eigenspace of $\rho_X(s')$.

Up to conjugation the group $G'$ has 2 Levi subgroups that are are isomorphic to $\GL_2\times\GL_2\times  \GL_2\times \GL_1/\GL_{1}^{\diag}$, we denote them by $L_1,L_2$. When $G'=G(\Sp_6\times \SO_4)/\GL_{1}^{\diag}$, we assume that the Half-Spin representation of $\Spin_4(\BC)$ appeared in $\rho_{s,+}$ ($\rho_{s,-}$) corresponds to the $\GSO_4$-component of $L_1$ (resp. $L_2$). When $G'=G(\Sp_3\times \SO_8)/\GL_{1}^{\diag}$, we assume that the Half-Spin representation of $\Spin_8(\BC)$ appeared in $\rho_{s,+}$ (resp. $\rho_{s,-}$) corresponds to the $\GSO_8$-component of $L_1$ (resp. $L_2$).

\begin{rmk}
We already explained how to relate the Levi subgroup of $\GSO_8$ to the Half-Spin representation in the previous subsection. For the $\GSO_4$ case, a Levi subgroup that is isomorphic to $\GL_2\times \GL_1$ corresponds to a Levi subgroup of its dual group $\GSpin_4(\BC)=(\GL_2(\BC)\times \GL_2(\BC))^0$ that is isomorphic to $\GL_2(\BC)\times \GL_1(\BC)$. Then it corresponds to the Half-Spin representation of $\GSpin_4(\BC)$ whose restriction to $\GL_2(\BC)$ is the standard representation.
\end{rmk}

Like the orthogonal group cases in the previous subsection, we can embed $H_0/Z_{G,H}\simeq \PGL_2$ diagonally into $L_i$ (denoted by $\nu_i$). We define
\begin{eqnarray*}
m_{geom}(\theta')&=&c_{\theta'}(1)+\sum_{T\in \CT_{ell}(H_0)}\frac{1}{2} \int_{T(F)/Z_{G,H}(F)} \\
&& \sum_{i=1}^{2} D^{H_0}(t)^3c_{\theta'}(\nu_i(t))dt,\\
m_{geom,D}(\theta')&=&\sum_{T\in \CT_{ell}(H_0)}\frac{1}{2} \int_{T(F)/Z_{G,H}(F)}\\
&&\sum_{i=1}^{2}(-1)^{i-1} D^{H_0}(t)^3c_{\theta'}(\nu_i(t))dt.
\end{eqnarray*}

\begin{prop}\label{prop GSp(10)}
Let $\theta$ (resp. $\theta_D$) be a quasi-character on $G/Z_{G,H}(F)$. Assume that $\theta$ (resp. $\theta_D$) is the endoscopic transfer of a quasi-character $\theta'$ of $G'(F)$ . We have
$$m_{geom}(\theta)=m_{geom}(\theta').$$
\end{prop}

\begin{proof}
The proof is very similar to the orthogonal group case in the previous subsection, we will skip it here.
\end{proof}

\begin{rmk}\label{GSp10 weak conjecture remark}
When $G'=G(Sp_2\times \SO_8)/\GL_{1}^{diag}$, $m_{geom}(\theta')$ contains the geometric multiplicity of the two $(\GSO_8\times \GL_2,\GL_2\ltimes U)$-models (note that these two models are differed by the outer automorphism). This is why Conjecture \ref{weak conjecture} for the model $(\GSp_{10},\GL_2\ltimes U)$ cannot imply Conjecture \ref{weak conjecture} for the model $(\GSO_8\times \GL_2,\GL_2\ltimes U)$. This is different from the $(\GU_6,\GU_2\ltimes U)$ case in the previous section. In that case, the endoscopic relation in Proposition \ref{prop GU(6)} contains the geometric multiplicity of the model $(\GU_4\times \GU_2,\GU_2\ltimes U)$ and the Gan-Gross-Prasad model for $U_4\times U_1$ (whose epsilon dichotomy is already known). Hence we can use the endoscopic identity to show that Conjecture \ref{weak conjecture} for the model  $(\GU_6,\GU_2\ltimes U)$ will imply Conjecture \ref{weak conjecture GU(4)} for the model $(\GU_4\times \GU_2,\GU_2\ltimes U)$.
\end{rmk}

Now we are ready to prove Theorem \ref{main theorem} for this case. The proof is very similar to the previous case and we will only give a sketch of it. Let $\Pi_\phi=\Pi_\phi(G)\cup \Pi_\phi(G_D)$ be a tempered $L$-packet with trivial central character. Assume that $\Pi_\phi(G)$ is not discrete with $|\Pi_\phi(G)|=1$. We first show that the unique distinguished element belongs to $\Pi_\phi(G)$ if and only if $\epsilon(\frac{1}{2},\Pi_\phi,\rho_X)=1$.

There are two situations. The first situation is when the packet $\Pi_\phi(G)$ is induced from a packet $\Pi_\phi(M)$ of a maximal Levi subgroup $M(F)$ of $G(F)$. If $M=\GL_4\times \GL_2$ (resp. $M=\GSp_6\times \GL_2$), this follows from Proposition \ref{GSO GSp parabolic induction} and Conjecture \ref{weak conjecture} for the model 
$$(\GL_4\times \GL_2,\GL_2\times \GL_2), \text{(resp.} \;(\GSp_6\times \GL_2,\GL_2\ltimes U)).$$ 
The model $(\GSp_6\times \GL_2,\GL_2\ltimes U)$ is smaller than $(G,H)$, by our assumption we know that Conjecture \ref{weak conjecture} holds for $(\GSp_6\times \GL_2,\GL_2\ltimes U)$. As for the model $(\GL_4\times \GL_2,\GL_2\times \GL_2)$, by our assumption and the fact that the model $(\GSO_8\times \GL_2,\GL_2\ltimes U)$ is smaller than $(G,H)$, we know that Conjecture \ref{weak conjecture} holds for $(\GSO_8\times \GL_2,\GL_2\ltimes U)$. Combining with Remark \ref{SO(8) implies GL(4)}, we know that Conjecture \ref{weak conjecture} holds for $(\GL_4\times \GL_2,\GL_2\times \GL_2)$. This also proves that if Conjecture \ref{weak conjecture} holds for the model $(\GSp_{10},\GL_2\ltimes U)$, then it also holds for the model $(\GSp_6\times \GL_2,\GL_2\ltimes U)$.

If $M=\GSp_8\times \GL_1$ or $\GSp_4\times \GL_3$, Proposition \ref{GSO GSp parabolic induction} implies that the unique distinguished element in the packet belongs to $\Pi_\phi(G)$. Also it is easy to see that the epsilon factor $\epsilon(\frac{1}{2},\Pi_\phi,\rho_X)$ is equal to 1 in this case.

Next we consider the case when the packet is discrete. By our assumption we have $|\Pi_\phi(G)|>1$. Hence there exists a proper elliptic extended endoscopic triple $(G',s',{}^L\eta)$ of $G$ such that $\phi$ factors through ${}^L\eta$ and $s'\in Z_\phi$. By the endoscopic identity in Proposition \ref{prop GSp(10)} and the same argument as in the orthogonal group case, we know that 
$$\sum_{\pi\in \Pi_\phi(G)} \tr(\chi_\pi(s'))m(\pi)=c_{\theta'}(1)=1,\;\sum_{\pi_D\in \Pi_\phi(G_D)} \tr(\chi_{\pi_D}(s'))m(\pi_D)=0$$
if the order of $s'$ is equal to 4, and
$$\sum_{\pi\in \Pi_\phi(G)} \tr(\chi_\pi(s'))m(\pi)=\frac{\epsilon(\frac{1}{2},\rho_{s',\phi,+})+\epsilon(\frac{1}{2},\rho_{s',\phi,-})}{2},$$
$$\sum_{\pi_D\in \Pi_\phi(G_D)} \tr(\chi_{\pi_D}(s'))m(\pi_D)=\frac{\epsilon(\frac{1}{2},\rho_{s',\phi,-})-\epsilon(\frac{1}{2},\rho_{s',\phi,+})}{2}$$
if the order of $s'$ is equal to 2. Here we need to use Conjecture \ref{weak conjecture} for the models $(\GSp_6\times \GL_2,\GL_2\ltimes U)$ and $(\GSO_8\times \GL_2,\GL_2\ltimes U)$ (both of them are smaller than the model $(\GSp_{10},\GL_2\ltimes U)$).

If the order of $s'$ is equal to 4, the endoscopic relation implies that the unique distinguished element belongs to $\Pi_\phi(G)$. By our discussion in Section \ref{sec epsilon factor} we know that $\epsilon(\frac{1}{2},\Pi_\phi,\rho_X)=1$.

If the order of $s'$ is equal to 2,  we know that the unique distinguished element belongs to $\Pi_\phi(G)$ if and only if $\epsilon(\frac{1}{2},\rho_{s',\phi,+})=\epsilon(\frac{1}{2},\rho_{s',\phi,-})$ which is equivalent to $\epsilon(\frac{1}{2},\Pi_\phi,\rho_X)=1$.

Now we prove Theorem \ref{main theorem}. Let $\omega_\phi\in \hat{S}_\phi$ correspond  to the unique distinguished element in the packet. By Remark \ref{distinguished is character} we know that $\omega_\phi$ is a character we view it as a character of $Z_\phi$. For $s\in S_\phi$, by Lemma \ref{lem extended endoscopic triple}, there exists an elliptic extended endoscopic triple $(G',s',{}^L\eta)$ of $G/Z_{G,H}$ such that $s'\in sZ_{\phi}^{\circ}$ and $\phi$ factors through ${}^L\eta$. We need to show that $\omega_\phi(s')=\omega_{\phi,H}(s)$. The above discussion implies that $\omega_\phi(s')=\omega_{\phi,H}(s)$ if $s'$ belongs to the center of the dual group.

If $s'$ does not belong to the center of the dual group, there are two cases. If the order of $s'$ is equal to 4, by the discussion above we know that the unique distinguished element belongs to $\Pi_\phi(G)$. By the definition of $\omega_{\phi,H}$ we know that $\omega_{\phi,H}(s)=1$. This implies that
$$\omega_\phi(s')=\tr(\omega_\phi(s'))=\sum_{\pi\in \Pi_\phi(G)} \tr(\chi_\pi(s'))m(\pi)=1=\omega_{\phi,H}(s).$$

If the order of $s'$ is equal to 2, by our discussion above, we have 
$$\sum_{\pi\in \Pi_\phi(G)} \tr(\chi_\pi(s'))m(\pi)=\frac{\epsilon(\frac{1}{2},\rho_{s',\phi,+})+\epsilon(\frac{1}{2},\rho_{s',\phi,-})}{2},$$
$$\sum_{\pi_D\in \Pi_\phi(G_D)} \tr(\chi_{\pi_D}(s'))m(\pi_D)=\frac{-\epsilon(\frac{1}{2},\rho_{s',\phi,+})+\epsilon(\frac{1}{2},\rho_{s',\phi,-})}{2}.$$
By the definition in Section \ref{sec epsilon factor}, we have
$$\omega_{\phi,H}(s)=\epsilon(\frac{1}{2},\rho_{s',\phi,-}).$$

We have two cases. If the unique distinguished element belongs to $\Pi_\phi(G)$, we have
$$\epsilon(\frac{1}{2},\Pi_\phi,\rho_X)=1,\;\epsilon(\frac{1}{2},\rho_{s',\phi,+})=\epsilon(\frac{1}{2},\rho_{s',\phi,-}).$$
This implies that
$$\omega_\phi(s')=\tr(\omega_\phi(s'))=\sum_{\pi\in \Pi_\phi(G)} \tr(\chi_\pi(s'))m(\pi)$$
$$=\frac{\epsilon(\frac{1}{2},\rho_{s',\phi,+})+\epsilon(\frac{1}{2},\rho_{s',\phi,-})}{2}=\epsilon(\frac{1}{2},\rho_{s',\phi,-}).$$

If the unique distinguished element belongs to $\Pi_\phi(G_D)$, we have
$$\epsilon(\frac{1}{2},\Pi_\phi,\rho_X)=-1,\;\epsilon(\frac{1}{2},\rho_{s',\phi,+})=-\epsilon(\frac{1}{2},\rho_{s',\phi,-}).$$
This implies that
$$\omega_\phi(s')=\tr(\omega_\phi(s'))=\sum_{\pi_D\in \Pi_\phi(G_D)} \tr(\chi_{\pi_D}(s'))m(\pi_D)$$
$$=\frac{-\epsilon(\frac{1}{2},\rho_{s',\phi,+})+\epsilon(\frac{1}{2},\rho_{s',\phi,-})}{2}=\epsilon(\frac{1}{2},\rho_{s',\phi,-}).$$
This finishes the proof of Theorem \ref{main theorem} for the model $(\GSp_{10},\GL_2\ltimes U)$.

\section{The model $(\GSp_6\times \GSp_4, (\GSp_4\times \GSp_2)^0)$}\label{sec GSp}
In this section, we consider the model $(\GSp_6\times \GSp_4, (\GSp_4\times \GSp_2)^0)$. In Section \ref{sec GSp 3}, we will define the models and the multiplicity formulas. We will also study the behaviors of the geometric multiplicities under endoscopic transfer and under parabolic induction. In Section \ref{sec:model-GSp4-GL2} we will discuss the smaller model $(\GSp_4\times \GL_2\times \GL_2,(\GL_2\times \GL_2)^0)$. In Section \ref{sec GSp 4}, we will prove the main theorem for this model.

\subsection{The models and the multiplicity formulas}\label{sec GSp 3}
Let $G=\GSp_6\times \GSp_4$ and 
$$H=(\GSp_4\times \GSp_2)^0=\{(g_1,g_2)\in \GSp_4\times \GSp_2 \mid  l(g_1)=l(g_2)\}.$$ 
There is a natural embedding from the group $(\GSp_4\times \GSp_2)^0$ into $\GSp_6$. Together with the projection map from $(\GSp_4\times \GSp_2)^0$ to $\GSp_4$, we get an embedding from $H$ to $G$. Similarly, we can define the model $(G_D,H_D)$ with 
$$G_D=\GSp_3(D)\times \GSp_2(D),\;H_D=(\GSp_2(D)\times \GL_1(D))^0.$$ 

Next we recall the definition of the geometric multiplicities from Section 9 of \cite{WZ2}. For $T\in \CT_{ell}(\GSp_2)$, let
$$T^{n,0}=\{(t_1,\cdots,t_n)\in T^n \mid \det(t_i)=\det(t_j)\;\text{for all}\;1\leq i,j\leq n\}.$$
We use $\iota_n$ to denote the diagonal embedding from $T$ to $T^{n,0}$. We can view $T^{n,0}$ as a maximal elliptic torus of $\GSp_{2n}$. Moreover, up to $\GSp_{2n}$-conjugation, there are $2^{n-1}$ distinct embeddings from $T^{n,0}$ to $\GSp_{2n}$.

When $n=2$, there are two embeddings $\nu_2,\nu_2'$ from $T^{2,0}$ to $\GSp_{4}$ and the centralizer of the image of $ \nu_2\circ\iota_2$ (resp. $\nu_2'\circ\iota_2$) in $\GSp_4$ is the quasi-split (resp. non quasi-split) unitary similitude group of two variables. Meanwhile, there are four embeddings from $T^{3,\circ}$ to $(\GSp_{4}\times \GSp_2)^0$ and there are two of them whose projection to $\GSp_4$ coincides with $\nu_2$. 
Composing with the embedding from $(\GSp_{4}\times \GSp_2)^0$ to $\GSp_6$, we get  two embeddings $\nu_{31},\nu_{32}$ from $T^{3,0}$ to $\GSp_6$. We use 
$$\nu_{T,i}=(\nu_{3i} \circ \iota_3)\times ( \nu_2\circ \iota_2)$$  
to denote the two embeddings from $T$ to $G$ (both factor through $H$). It is easy to see that these two embeddings are conjugated to each other in $H$ and we will use $\nu_T$ to denote one of it.

Meanwhile, let $\iota_{1,2}$ be the embedding from $T^{2,0}$ to $T^{3,0}$ given by 
$$(t_1,t_2)\mapsto (t_1,t_2,t_2).$$ 
Among the four embeddings from $T^{3,0}$ to $\GSp_{6}$, there are two of them (denoted by $\nu_3,\nu_3'$) such that the centralizers in $\GSp_{6}$ of the image of $\nu_3\circ \iota_{1,2}$ and $\nu_3'\circ \iota_{1,2}$ are quasi-split (the centralizer is the quasi-split unitary similitude group of two variables times an abelian group). Up to conjugation we may assume that $\nu_3,\nu_3'$ factor  through $(\GSp_{4}\times \GSp_2)^0$ and the projection to $\GSp_4$ of $\nu_3\circ \iota_{1,2}$ (resp. $\nu_3'\circ \iota_{1,2}$) is equal to $\nu_2$ (resp. $\nu_2'$). We use $$\nu_{T^{2,0},1}=(\nu_{3}\circ \iota_{1,2})\times \nu_2,\;\nu_{T^{2,0},2}=(\nu_{3}'\circ \iota_{1,2})\times \nu_2'$$ 
to denote the two embeddings from $T^{2,0}$ to $G$. Both of them factor through $H$.

Finally, for $T_1,T_2\in \CT_{ell}(\GSp_2)$ with $T_1\neq T_2$ (this will not happen in the Archimedean case), let
$$(T_1\times T_2)^0=\{(t_1,t_2)\in T_1\times T_2 \mid \det(t_1)=\det(t_2)\}. $$
Similarly, we can define $(T_1\times T_2\times T_2)^0$. Up to conjugation, there is only one embedding from $(T_1\times T_2)^0$ to $\GSp_4$ and there are two embeddings from $(T_1\times T_2\times T_2)^0$ to $\GSp_6$. The two embeddings induce  two embeddings from $(T_1\times T_2)^0$ to $\GSp_6$ (we first map $T_2$ diagonally into $(T_2\times T_2)^{0}$). We let $\nu$ be the embedding such that the centralizer of its image is quasi-split (the centralizer of the other embedding is not quasi-split). 
Up to conjugation we may assume that $\nu$ factors through $(\GSp_4\times \GSp_2)^0$ and its projection to $\GSp_4$ is equal to the embedding from $(T_1\times T_2)^0$ to $\GSp_4$. This gives us an embedding $\nu_{T_1,T_2}$ from $(T_1\times T_2)^{0}$ to $G$ that factors through $H$. 

Let $\theta$ be a quasi-character on $G(F)$. Define the geometric multiplicity to be 
$$m_{geom}(\theta)=c_{\theta}(1)+\sum_{T\in \CT_{ell}(H)}|W(H,T)|^{-1}\int_{T(F)/Z_{G,H}(F)} D^H(t)\theta(t)d t$$
\begin{eqnarray*} 
&&+\frac{1}{2} \sum_{T\in \CT_{ell}(\GSp_2)}\big(\int_{T(F)/Z_{\GL_2}(F)} D^H(\nu_{T}(t))c_{\theta}(\nu_{T}(t))d t\\
&&+\sum_{i\in \{1,2\}}\int_{T^{2,0}(F)/Z_{\GL_2}(F)} D^H(\nu_{T^{2,0},i}(t))c_{\theta}(\nu_{T^{2,0},i}(t))d t \big)\\
&&+\frac{1}{4} \sum_{T_1,T_2\in \CT_{ell}(\GSp_2),T_1\neq T_2}\int_{(T_1\times T_2)^0(F)/Z_{\GL_2}(F)^{\diag}}\\
&&D^H(\nu_{T_1,T_2}(t)) c_{\theta}(\nu_{T_1,T_2}(t))d t.
\end{eqnarray*}

Similarly, for the quaternion version $(G_D,H_D)$, we can also define the embeddings $\nu_{T_D}$, $\nu_{T_{D}^{2,0},i}$, $\nu_{T_{1,D},T_{2,D}}$ for $T_D$, $T_{1,D}$, $T_{2,D}\in \CT_{ell}(\GL_1(D))=\CT_{ell}(\GSp_2)$ with $T_{1,D}\neq T_{2,D}$. Let $\theta_D$ be a quasi-character on $G_D(F)$. Define the geometric multiplicity $m_{geom}(\theta_D)$ to be
\begin{eqnarray*} 
&&\sum_{T_D\in \CT_{ell}(H_D)}|W(H_D,T_D)|^{-1}\int_{T_D(F)/Z_{G_D,H_D}(F)} D^{H_D}(t)\theta_{D}(t)d t\\
&&+\frac{1}{2}\sum_{T_D\in \CT_{ell}(\GL_1(D))}\big(\int_{T_D(F)/Z_{\GL_1(D)}(F)} D^{H_D}(\nu_{T_D}(t))c_{\theta_D}(\nu_{T_D}(t))d t\\
&&+\sum_{i\in \{1,2\}}\int_{T_{D}^{2,0}(F)/Z_{\GL_1(D)}(F)} D^{H_D}(\nu_{T_{D}^{2,0},i}(t))c_{\theta_D}(\nu_{T_{D}^{2,0},i}(t))d t \big)\\
&&+\frac{1}{4}\sum_{T_{1,D},T_{2,D}\in \CT_{ell}(\GL_1(D)),T_{1,D}\neq T_{2,D}}\int_{(T_{1,D}\times T_{2,D})^0(F)/Z_{\GL_1(D)}(F)^{\diag}}\\
&&D^{H_D}(\nu_{T_{1,D},T_{2,D}}(t)) c_{\theta_D}(\nu_{T_{1,D},T_{2,D}}(t))d t.
\end{eqnarray*} 
In our previous paper \cite{WZ2}, we have proved the multiplicity formulas
$$m(\pi)=m_{geom}(\theta_\pi),\;m(\pi_D)=m_{geom}(\theta_{\pi_D})$$
for all tempered representations.

\begin{rmk}
If $F=\BR$, the above integrals need to be regularized, i.e. we replace $D^H(\cdot)$ (resp. $D^{H_D}(\cdot)$) by $D^G(\cdot)^{1/2}(D^{H}(\cdot)^{-2} D^G(\cdot))^{s-1/2}$ (resp. $D^{G_D}(\cdot)^{1/2}(D^{H_D}(\cdot)^{-2} D^{G_D}(\cdot))^{s-1/2}$) and take the limit $\lim_{s\rightarrow 0^+}$. Since this regularization does not affect our later computation, to simplify the notation, we will not include this regularization in the expression of the multiplicity formula.
\end{rmk}

Next we study the behavior of the geometric multiplicities under parabolic induction. Let $M$ be a proper Levi subgroup of $G$ and $\theta^M$ be a quasi-character on $M(F)$. Let $L(F)$ (resp. $L_D(F)$) be the standard Levi subgroup of $G(F)$ (resp. $G_D(F)$) that is isomorphic to $$(\GL_2(F)\times \GL_2(F))\times (\GL_2(F)\times \GL_1(F))$$ 
$$(\text{resp.}\; (\GL_1(D)\times \GL_1(D))\times (\GL_1(D)\times \GL_1(F))).$$ 
If $M$ does not contain the Levi subgroup $L$ up to conjugation, define 
$$m_{geom}(\theta^M)=c_{\theta^M}(1).$$ 
Otherwise, $M$ corresponds to a proper Levi subgroup $M_D$ of $G_D$.
Moreover, up to conjugation we may assume that $L\subset M$ and $L_D\subset M_D$. Let $\theta_{D}^{M_D}$ be a quasi-character on $M_D(F)$. We have a natural embedding $\iota$ (resp. $\iota_D$) of $\GSp_2(F)$ (resp. $\GL_1(D)$) into $L(F)$ (resp. $L_D(F)$) given by $h\mapsto \diag(h,h,h)\times \diag(h,h)$. 
When the Levi subgroup $M$ is not isomorphic to $(\GL_2(F)\times \GL_2(F))\times \GSp_4(F)$, we define
\begin{eqnarray*}
m_{geom}(\theta^M)&=&c_{\theta^M}(1)+\frac{1}{2}\sum_{T\in \CT_{ell}(\GSp_2)}\int_{T(F)/Z_{\GSp_2}(F)}\\
&& D^M(\iota(t))^{1/2}D^{\GSp_2}(t)^{-1/2} c_{\theta^M}(\iota(t))dt,\\
m_{geom}(\theta_{D}^{M_D})&=&\frac{1}{2}\sum_{T_D\in \CT_{ell}(\GL_1(D))} \int_{T_D(F)/Z_{\GL_1(D)}(F)}\\
&&  D^{M_D}(\iota_D(t))^{1/2}D^{\GL_1(D)}(\iota_D(t))^{-1/2} c_{\theta_{D}^{M_D}}(\iota_D(t))dt.
\end{eqnarray*}

When $M$ is isomorphic to $(\GL_2(F)\times \GL_2(F))\times \GSp_4(F)$, for  $T_1$, $T_2\in \CT_{ell}(\GSp_2)$ (resp. $T_{1,D}$, $T_{2,D}\in \CT_{ell}(\GL_1(D))$), there is a natural embedding from $(T_1\times T_2)^{0}$ (resp. $(T_{1,D}\times T_{2,D})^{0}$) into the Levi subgroup $\GL_2\times \GL_2$ of $\GSp_6$ (resp. $\GL_1(D)\times \GL_1(D)$ of $\GSp_3(D)$) given by $(t_1,t_2)\mapsto \diag(t_1,t_2,t_1)$. When $T_1\neq T_2$ (resp. $T_{1,D}\neq T_{2,D}$), up to conjugation there is a unique embedding from $(T_1\times T_2)^{0}$ (resp. $(T_{1,D}\times T_{2,D})^{0}$) into $\GSp_4$ (resp. $\GSp_2(D)$). This gives us an embedding, denoted by $\iota_{T_1,T_2}$ (resp. $\iota_{T_{1,D},T_{2,D}}$), from $(T_1\times T_2)^{0}$ (resp. $(T_{1,D}\times T_{2,D})^{0}$) into $M$ (resp. $M_D$). 

When $T=T_1= T_2$ (resp. $T_D=T_{1,D}= T_{2,D}$), up to conjugation there are two embeddings from $T^{2,0}=(T_1\times T_2)^{0}$ (resp. $T_{D}^{2,0}=(T_{1,D}\times T_{2,D})^{0}$) into $\GSp_4$ (resp. $\GSp_2(D)$). This gives us two embeddings, denoted by $\iota_{T^2,i}$ (resp. $\iota_{T_{D}^2,i}$) with $1\leq i\leq 2$, from $T^{2,0}$ (resp. $T_{D}^{2,0}$) into $M$ (resp. $M_D$). We define $m_{geom}(\theta^M)$ to be

\begin{eqnarray}\label{GSpin multiplicity 1}
&&c_{\theta^M}(1)+\frac{1}{2}\sum_{T\in \CT_{ell}(GSp_2)}\big(\int_{T(F)/Z_{\GSp_2}(F)} D^{\GSp_2}(t)^{2} c_{\theta^M}(\iota(t))dt  \\
&&+\sum_{i=1}^{2}\int_{T^{2,0}(F)/Z_{\GSp_2}(F)} D^{\GSp_2}(t_1)D^{\GSp_2}(t_2) c_{\theta^M}(\iota_{T^{2,0},i}(t_1,t_2))dt_1dt_2\big) \nonumber \\
&&+\frac{1}{4}\sum_{T_1,T_2\in \CT_{ell}(\GSp_2),T_1\neq T_2} \int_{(T_1\times T_2)^{0}(F)/Z_{\GSp_2}(F)^{diag}} \nonumber \\
&&D^{\GSp_2}(t_1)D^{\GSp_2}(t_2) c_{\theta^M}(\iota_{T_1,T_2}(t_1,t_2))dt_1dt_2,\nonumber
\end{eqnarray}

For the quaternion side, we define $m_{geom}(\theta_{D}^{M_D})$ to be
\begin{equation}\label{GSpin multiplicity 2}
\frac{1}{2}\sum_{T_D\in \CT_{ell}(\GL_1(D))} \big(\int_{T_D(F)/Z_{\GL_1(D)}(F)}D^{\GL_1(D)}(\iota_D(t))^{2}c_{\theta_{D}^{M_D}}(\iota_D(t))dt 
\end{equation}
$$+\sum_{i=1}^{2}\int_{T_{D}^{2,0}(F)/Z_{\GL_1(D)}(F)} D^{\GL_1(D)}(t_1) D^{\GL_1(D)}(t_2)c_{\theta_{D}^{M_D}}(\iota_{T_{D}^{2,0},i}(t_1,t_2))dt_1dt_2\big)$$
\begin{eqnarray}
&&+\frac{1}{4}\sum_{T_{1,D},T_{2,D}\in \CT_{ell}(\GL_1(D)),T_{1,D}\neq T_{2,D}}\int_{(T_{1,D}\times T_{2,D})^{0}(F)/Z_{\GL_1(D)}(F)^{diag}} \nonumber \\
&&D^{\GL_1(D)}(t_1)D^{\GL_1(D)}(t_2) c_{\theta_{D}^{M_D}}(\iota_{T_{1,D},T_{2,D}}(t_1,t_2))dt_1dt_2.\nonumber
\end{eqnarray}
 
The following proposition is a direct consequence of Proposition \ref{germ parabolic induction}. 

\begin{prop}\label{GSp(6)xGSp(4) parabolic induction}
Let $\theta$ (resp. $\theta_D$) be a quasi-character on $G(F)$ (resp. $G_D(F)$). Assume that $\theta$ (resp. $\theta_D$) is the parabolic induction of a quasi-character $\theta^M$ (resp. $\theta_{D}^{M_D}$) of a proper Levi subgroup $M$ of $G$ (resp. $M_D$ of $G_D$). We have
$$m_{geom}(\theta)=m_{geom}(\theta^M),\;m_{geom}(\theta_D)=m_{geom}(\theta_{D}^{M_D}).$$
\end{prop}

Next we study the behavior of the geometric multiplicities under endoscopic transfer. Let $(G',s',{}^L\eta)$ be a proper elliptic extended endoscopic triple of $G/Z_{G,H}$. Up to multiplying $s'$ by  an element in the neutral component of the center of the dual group, we may assume that 
$$s'=(s_1,s_2)\in \Spin_7(\BC)\times \Spin_5(\BC)\subset \widehat{G/Z_{G,H}},$$
$$\widehat{G/Z_{G,H}}=\{(g_1,g_2)\in \GSpin_7(\BC)\times \GSpin_5(\BC) \mid  l(g_1)l(g_2)=1\}.$$
We will only consider the case when one of $s_i$ is the identity element. Under this assumption, we have $G'=\GSO_6\times \GSp_4/\GL_{1}^{\diag},\;\GSp_6\times \GSO_4/\GL_{1}^{\diag}$ or $G(\Sp_2\times \SO_4)\times \GSp_4/\GL_{1}^{\diag}$. If $G'=\GSO_6\times \GSp_4/\GL_{1}^{\diag}$, we define
$$m_{geom}(\theta')=c_{\theta'}(1),\;m_{geom,D}(\theta')=0.$$

If $G'=\GSp_6\times \GSO_4/\GL_{1}^{\diag}$, we have $s_1=1$ and the projection of $s_2\in \Spin_5(\BC)$ to $\SO_5(\BC)$ is conjugated to $\diag(1,-I_4)$. As we explained in Section \ref{sec epsilon factor}, when we restrict the representation $\rho_X$ to $\hat{G}'=(\widehat{G/Z_{G,H}})_{s'}$, we can decompose it as $\rho_{s,+}\oplus \rho_{s,-}$ where $\rho_{s,+}$ (resp. $\rho_{s,-}$) is the tensor product of the Spin representation of $\Spin_7(\BC)$ with a Half-Spin representation of $\Spin_4(\BC)$ and it is the $+1$ (resp. $-1$) eigenspace of $\rho_X(s')$.

We have an embedding from $\GSp_2=\GL_2$ into $\GSp_6$ given by $h\mapsto \diag(h,h,h)$. On the other hand, the group $\GSO_4$ has two Siegel parabolic subgroups which give us two embeddings from $\GL_2$ into $\GSO_4$. Combining these embeddings we get two embeddings from $\PGL_2$ into $G'$ which will be denoted by $\iota_1$ and $\iota_2$. We assume that the Half-Spin representation of $\Spin_4(\BC)$ appeared in $\rho_{s,+}$ (resp. $\rho_{s,-}$) corresponds to the Siegel Levi subgroup of $\GSO_4$ associated to $\iota_1$ (resp. $\iota_2$). We define
$$m_{geom}(\theta')=c_{\theta'}(1)+\sum_{T\in \CT_{ell}(\PGL_2)} \frac{1}{2}\int_{T(F)} D^{\GL_2}(t)^3 (c_{\theta'}(\iota_1(t))+c_{\theta'}(\iota_2(t)))dt,$$
$$m_{geom,D}(\theta')=\sum_{T\in \CT_{ell}(\PGL_2)} \frac{1}{2}\int_{T(F)} D^{\GL_2}(t)^3 (c_{\theta'}(\iota_1(t))-c_{\theta'}(\iota_2(t)))dt.$$

If $G'=G(\Sp_2\times \SO_4)\times \GSp_4/\GL_{1}^{diag}$, we have $s_2=1$ and the projection of $s_1\in \Spin_7(\BC)$ to $\SO_7(\BC)$ is conjugated to $\diag(I_3,-I_4)$. When we restrict the representation $\rho_X$ to $\hat{G}'=(\widehat{G/Z_{G,H}})_{s'}$, we can decompose it as $\rho_{s,+}\oplus \rho_{s,-}$ where $\rho_{s,+}$ (resp. $\rho_{s,-}$) is the tensor product of the Spin representation of $\Spin_3(\BC)$, a Half-Spin representation of $\Spin_4(\BC)$ and the Spin representation of $\Spin_5(\BC)$, and it is the $+1$ (resp. $-1$) eigenspace of $\rho_X(s')$.

Like in the previous case, we still have two embeddings from $\PGL_2$ into $G'$ which will be denoted by $\iota_1$ and $\iota_2$.  For $T_1,T_2\in \CT_{ell}(\GL_2)$, up to conjugation there are either one (when $T_1\neq T_2$) or two (when $T_1=T_2$) embeddings from $(T_1\times T_2)^0(F)$ into $\GSp_4(F)$. We fix one of such embeddings (the choice does not matter since $\theta'$ is stable). 
Meanwhile, we can embed $T_1$ into $\GSp_2$, and there are two ways to embed $T_2$ into $\GSO_4$ (again corresponding to the two Siegel parabolic subgroups). 
This gives us two embeddings from $(T_1\times T_2)^0(F)/\GL_{1}(F)^{\diag}$ into $G'(F)$ which will be denoted by $\nu_{T_1,T_2,i},\;1\leq i\leq 2$.
 We still assume that the Half-Spin representation of $\Spin_4(\BC)$ appeared in $\rho_{s,+}$ (resp. $\rho_{s,-}$) corresponds to the Siegel Levi subgroup of $\GSO_4$ associated to $\iota_1, \nu_{T_1,T_2,1}$ (resp. $\iota_2,\nu_{T_1,T_2,2}$). We define $m_{geom}(\theta')$ to be
\begin{eqnarray*}
&&c_{\theta'}(1)+\frac{1}{2}\sum_{T\in \CT_{ell}(\PGL_2)} \int_{T(F)} D^{\GL_2}(t)^2 (c_{\theta'}(\iota_1(t))+c_{\theta'}(\iota_2(t)))dt\\
&&+\sum_{T_1,T_2\in \CT_{ell}(\GL_2)}d(T_1,T_2)\int_{(T_1\times T_2)^0(F)/\GL_{1}(F)^{\diag}}\\
&& D^{\GL_2\times \GL_2}(t_1, t_2) (c_{\theta'}(\nu_{T_1,T_2,1}(t_1,t_2))+c_{\theta'}(\nu_{T_1,T_2,2}(t_1,t_2)))dt_1dt_2,
\end{eqnarray*}
For the quaternion side, we define $m_{geom,D}(\theta')$ to be
\begin{eqnarray*}
&&\frac{1}{2}\sum_{T\in \CT_{ell}(\PGL_2)} \int_{T(F)} D^{\GL_2}(t)^2 (c_{\theta'}(\iota_1(t))-c_{\theta'}(\iota_2(t)))dt\\
&&+\sum_{T_1,T_2\in \CT_{ell}(\GL_2)}d(T_1,T_2)\int_{(T_1\times T_2)^0(F)/\GL_{1}(F)^{\diag}}\\
&& D^{\GL_2\times \GL_2}(t_1, t_2) (c_{\theta'}(\nu_{T_1,T_2,1}(t_1,t_2))-c_{\theta'}(\nu_{T_1,T_2,2}(t_1,t_2)))dt_1dt_2
\end{eqnarray*}
where $d(T_1,T_2)=1$ if $T_1=T_2$ and $d(T_1,T_2)=\frac{1}{4}$ if $T_1\neq T_2$.

\begin{prop}\label{prop GSp(6)xGSp(4)}
Let $\theta$ (resp. $\theta_D$) be a quasi-character on $G/Z_{G,H}(F)$. Assume that $\theta$ (resp. $\theta_D$) is the endoscopic transfer of a stable quasi-character $\theta'$ of $G'(F)$. 
We have
$$m_{geom}(\theta)=m_{geom}(\theta'),\;m_{geom}(\theta_D)=m_{geom,D}(\theta').$$
\end{prop}

\begin{proof}
We will only prove the split case, and the quaternion case follows from a similar argument. Like in Proposition \ref{prop GSO(12)}, the only difference between the split case and the quaternion case is that there is an extra $-1$ in the transfer factor for certain conjugacy classes. 

The identity $c_\theta(1)=c_{\theta'}(1)$ is easy and we will skip the proof. Next we study the terms corresponding  to $T\in \CT_{ell}(H)$ in $m_{geom}(\theta)$. We would like to show that these terms are equal to zero. 
To do this, we only need to show that the transfer factor is non-trivial for elements in $T(F)$. We first describe the elliptic conjugacy classes in $H(F)$. We can view $H$ as a subgroup of $\GSp_6$. Then the elliptic conjugacy classes of $H(F)$ are just the elliptic conjugacy classes of $\GSp_6$ that have a nonempty intersection  with $H(F)$. To be specific, $\GSp_6$ has three types of elliptic conjugacy classes corresponding to
\begin{enumerate}
\item $(K,K_{\pm})$, $K_{\pm}/F$ is a cubic extension and $K/K_{\pm}$ is a quadratic extension. 
\item $(E,E_{\pm})\cup (E_1,F)$, $E/E_{\pm},\;E_1/F$ and $E_{\pm}/F$ are quadratic extensions.
\item $(E_1,F)\cup (E_2,F)\cup (E_3,F)$, $E_i/F$ is a quadratic extension for $1\leq i\leq 3$.
\end{enumerate}
The first type has no intersection with $H(F)$, so we only have Type (2) and Type (3). 

For a Type (2) (resp. Type (3)) conjugacy class, if there exists a conjugacy class in $G'=(G(\Sp_2\times \SO_4)\times \GSp_4)/\GL_{1}^{diag}$ corresponding to it, then $\eta_{E/E_{\pm}}$ is trivial on $F^{\times}$ (resp. $E_i=E_j$ for some $i\neq j$). In both cases the group 
$$\Pi_{i\in I} (ker(\tr_{F_i/F_{\pm i}})\cap F_{i}^{\times})/Im(N_{F_i/F_{\pm i}}) /\sim$$ 
is isomorphic to $\BZ/2\BZ$ or $(\BZ/2\BZ)^2$. Note that $(\BZ/2\BZ)^2$ only happens in Type (3) when $E_1=E_2=E_3$.

On the other hand, for a Type (2) (resp. Type (3)) conjugacy class, if there exists a conjugacy class in $G'=\GSO_6\times \GSp_4/\GL_{1}^{diag}$ corresponding to it, then $\eta_{E/E_{\pm}}|_{F^\times}=\eta_{E_1/F}$ (resp. $E_i\neq E_j$ for all $i\neq j$ and $E_1$ is contained in $E_2\otimes_F E_3$). In both cases the group 
$$\Pi_{i\in I} (ker(\tr_{F_i/F_{\pm i}})\cap F_{i}^{\times})/Im(N_{F_i/F_{\pm i}}) /\sim$$ 
is isomorphic to $\BZ/2\BZ$.

By our definition of the transfer factors in Section \ref{section transfer factor}, under the endoscopic relation between $\GSp_6$ and $\GSO_6$ (resp. $G(\Sp_2\times \SO_4)$), the transfer factors associated to the conjugacy classes of Type (2) and (3) are non-trivial. Moreover, they are equal to a constant times a non-trivial character on 
$$\Pi_{i\in I} (ker(\tr_{F_i/F_{\pm i}})\cap F_{i}^{\times})/Im(N_{F_i/F_{\pm i}}) /\sim.$$ 
This implies that when $G'=\GSO_6\times \GSp_4/\GL_{1}^{diag}$ or $G(\Sp_2\times \SO_4)\times \GSp_4/\GL_{1}^{diag}$, the term corresponding to $T\in \CT_{ell}(H)$ in $m_{geom}(\theta)$ is equal to 0.  

On the other hand, if $G'=\GSp_6\times \GSO_4/\GL_{1}^{diag}$, we need to study the projection of the above conjugacy classes to $\GSp_4$. The projection of Type (2) conjugacy classes to $\GSp_4$ corresponds to  
\begin{itemize}
\item[(2)'] $(E,E_{\pm})$, $E/E_{\pm}$ and $E_\pm/F$ are quadratic extensions.
\end{itemize}
The projection of Type (3) conjugacy classes to $\GSp_4$ corresponds to 
\begin{itemize}
\item[(3)'] $(E_1,F)\cup (E_2,F)$, $E_i/F$ is a quadratic extension for $1\leq i\leq 2$.
\end{itemize}
For Type (2)' (resp. (3)') conjugacy classes, there exist conjugacy classes in $\GSO_4$ corresponding to them if and only if $\eta_{E/E_{\pm}}$ is trivial on $F^{\times}$ (resp. $E_1=E_2$). If this is the case,  the group 
$$\Pi_{i\in I} (ker(\tr_{F_i/F_{\pm i}})\cap F_{i}^{\times})/Im(N_{F_i/F_{\pm i}}) /\sim$$ 
is isomorphic to $\BZ/2\BZ$. Moreover, the transfer factors are non-trivial, and  equal to a constant times the non-trivial character on $$\Pi_{i\in I} (ker(\tr_{F_i/F_{\pm i}})\cap F_{i}^{\times})/Im(N_{F_i/F_{\pm i}}) /\sim.$$ 
This shows that the term corresponding to $T\in \CT_{ell}(H)$ in $m_{geom}(\theta)$ is equal to 0. 

Then we need to study the terms correspond to $T,\;T^{2,0}$ and $(T_1\times T_2)^0$ in $m_{geom}(\theta)$ for $T,T_i\in \CT_{ell}(\GSp_2)$ with $T_1\neq T_2$. First we study the term corresponding to $T$. Like in the previous cases, $T$ corresponds to a quadratic extension $E_T$ of $F$ and we can view $t\in T(F)$ as an element of $E_{T}^{\times}$. For $t\in T(F)$, the regular germ 
$$D^H(\nu_{T}(t))c_\theta(\nu_{T}(t))$$ 
is equal to $\frac{D^{\GL_2}(t)^{-1/2}}{4}$ times the limits of $D^G(\cdot)^{1/2}\theta(\cdot)$ at
$$((E_T\oplus E_T,E_T,(\lambda t,\lambda^{-1} \bar{t}))\cup (E_T,F,t)) \;\times \;(E_T\oplus E_T,E_T,(\lambda t,\lambda^{-1} \bar{t}))$$
as $\lambda\rightarrow 1$. 
In this case the group 
$$\Pi_{i\in I} (ker(\tr_{F_i/F_{\pm i}})\cap F_{i}^{\times})/Im(N_{F_i/F_{\pm i}}) /\sim$$ 
is the trivial group. 

If $G'=\GSO_6\times \GSp_4/\GL_{1}^{diag}$, there is no conjugacy class in $G'$ corresponding to the above conjugacy classes of $G(F)$. 
As a result, the term corresponding to $T$ in $m_{geom}(\theta)$ is equal to 0. 

If $G'=\GSp_6\times \GSO_4/\GL_{1}^{diag}$ or $G(\Sp_2\times \SO_4)\times \GSp_4/\GL_{1}^{diag}$, the transfer factors are trivial since the quadratic character $\eta_{E_T\oplus E_T/E_T}$ is trivial. By the same argument as in the cases of the previous section, we know that 
$$D^H(\nu_{T}(t))c_\theta(\nu_{T}(t))=D^{\GL_2}(t)^k(c_{\theta'}(\iota_1(t))+c_{\theta'}(\iota_2(t)))$$
where $k=2$ if $G'=G(\Sp_2\times \SO_4)\times \GSp_4/\GL_{1}^{diag}$ and $k=3$ if $G'=\GSp_6\times \GSO_4/\GL_{1}^{diag}$. Hence the terms correspond to $T$ in $m_{geom}(\theta)$ and $m_{geom}(\theta')$ are equal to each other. 

For the term corresponding to $T^{2,0}$, the regular germ (here $t=(t_1,t_2)\in T^{2,0}(F)$)
$$D^H(\nu_{T^{2,0},i}(t))\theta(\nu_{T^{2,0},i}(t)), \;1\leq i\leq 2$$  
is equal to 
$$\frac{D^H(\nu_{T^{2,0},i}(t))D^G(\nu_{T^{2,0},i}(t))^{-1/2}}{2}$$ 
times the limit of $D^G(\cdot )^{1/2}\theta(\cdot)$ at
$$((E_T\oplus E_T,E_T,(\lambda t_2,\lambda^{-1} \bar{t}_{2}))\cup (E_T,F,t_1)) \;\times \;((E_T,F,t_1)\cup (E_T,F,t_2),c_i)$$
as $\lambda\rightarrow 1$. In this case the group 
$$\Pi_{i\in I} (ker(\tr_{F_i/F_{\pm i}})\cap F_{i}^{\times})/Im(N_{F_i/F_{\pm i}}) /\sim$$ 
is isomorphic to $\BZ/2\BZ$. More specifically, the $\GSp_6$-component of this group is trivial and the $\GSp_4$-component of this group is $\BZ/2\BZ$. We use $c_i,\;1\leq i\leq 2$ to denote the two elements in this group. 

If $G'=\GSO_6\times \GSp_4/\GL_{1}^{diag}$, there is no conjugacy class in $G'$ corresponding to the above conjugacy classes of $G(F)$. As a result, the term corresponding to $T^{2,0}$ in $m_{geom}(\theta)$ is equal to 0. 

If $G'=\GSp_6\times \GSO_4/\GL_{1}^{diag}$, the transfer factors are non-trivial and they are equal to a constant times the sign character of $$\Pi_{i\in I} (ker(\tr_{F_i/F_{\pm i}})\cap F_{i}^{\times})/Im(N_{F_i/F_{\pm i}}) /\sim.$$ 
This implies that the term corresponds to $T^{2,0}$ in $m_{geom}(\theta)$ is equal to 0. 

If $G'=G(\Sp_2\times \SO_4)\times \GSp_4/\GL_{1}^{diag}$, the transfer factors are trivial since the quadratic character $\eta_{E_T\oplus E_T/E_T}$ is trivial. By the same argument as in the cases of the previous section, we know that 
$$D^H(\nu_{T^{2,0},1}(t))c_\theta(\nu_{T^{2,0},1}(t))+D^H(\nu_{T^{2,0},2}(t))c_\theta(\nu_{T^{2,0},2}(t))$$
$$=2\cdot D^{\GL_2\times \GL_2}(t_1,t_2)(c_{\theta'}(\nu_{T,T,1}(t_1,t_2))+c_{\theta'}(\nu_{T,T,2}(t_1,t_2))).$$
This shows that the terms corresponding to $T^{2,0}$ in $m_{geom}(\theta)$ and $m_{geom}(\theta')$ are equal to each other. 

Finally, for the term corresponds to $(T_1\times T_2)^{0}$ ($T_1\neq T_2$) in $m_{geom}(\theta)$, the regular germ $D^H(\nu_{T_1,T_2}(t))\theta(\nu_{T_1,T_2}(t))$ (here $t=(t_1,t_2)\in (T_1\times T_2)^{0}(F)$) is equal to $$\frac{D^H(\nu_{T_1,T_2}(t))D^G(\nu_{T_1,T_2}(t))^{-1/2}}{2}$$ 
times the limit of $D^G(\cdot )^{1/2}\theta(\cdot)$ at
$$((E_{T_2}\oplus E_{T_2},E_{T_2},(\lambda t_2,\lambda^{-1} \bar{t}_{2}))\cup (E_{T_1},F,t_1)) \;\times \;((E_{T_1},F,t_1)\cup (E_{T_2},F,t_2))$$
as $\lambda\rightarrow 1$. In this case the group 
$$\Pi_{i\in I} (ker(\tr_{F_i/F_{\pm i}})\cap F_{i}^{\times})/Im(N_{F_i/F_{\pm i}}) /\sim$$ is trivial.  

If $G'=\GSO_6\times \GSp_4/\GL_{1}^{diag}$ or $\GSp_6\times \GSO_4/\GL_{1}^{diag}$, there is no conjugacy class in $G'$ corresponding to the above conjugacy classes of $G(F)$. As a result, the term corresponding to $(T_1\times T_2)^{0}$ in $m_{geom}(\theta)$ is equal to 0. 

If $G'=G(\Sp_2\times \SO_4)\times \GSp_4/\GL_{1}^{diag}$, the transfer factors are trivial since the quadratic character $\eta_{E_{T_2}\oplus E_{T_2}/E_{T_2}}$ is trivial. By the same argument as in the cases of the previous section, we know that 
\begin{eqnarray*}
&&D^H(\nu_{T_1,T_2}(t))c_\theta(\nu_{T_1,T_2}(t))\\
&=& D^{\GL_2\times \GL_2}(t_1,t_2)(c_{\theta'}(\nu_{T_1,T_2,1}(t_1,t_2))+c_{\theta'}(\nu_{T_1,T_2,2}(t_1,t_2))).
\end{eqnarray*}
This shows that the terms correspond to $(T_1\times T_2)^{0}$ in $m_{geom}(\theta)$ and $m_{geom}(\theta')$ are equal to each other. This finishes the proof of the proposition.
\end{proof}

\subsection{The model $(\GSp_4\times \GL_2\times \GL_2,(\GL_2\times \GL_2)^0)$} \label{sec:model-GSp4-GL2}
In this section we discuss the model $(\GSp_4\times \GL_2\times \GL_2,(\GL_2\times \GL_2)^0)$. This model is smaller than $(\GSp_6\times \GSp_4, (\GSp_4\times \GSp_2)^0)$ and we need to assume the weak conjecture holds for this model in order to prove Theorem \ref{main theorem} for $(\GSp_6\times \GSp_4, (\GSp_4\times \GSp_2)^0)$.

Let $G=\GSp_4\times \GL_2\times \GL_2$, we have an embedding from $(\GL_2\times \GL_2)^0$ into $\GSp_4$ (resp. $\GL_2\times \GL_2$) which induces a diagonal embedding from this group to $G$, we will use $H\subset G$ to denote the image. There is also a quaternion version $(G_D,H_D)$ of this model with $G_D=\GSp_2(D)\times \GL_1(D)\times \GL_1(D)$ and $H_D=(\GL_1(D)\times \GL_1(D))^0$. The models $(G,H)$ and $(G_D,H_D)$ are essentially the Gan--Gross--Prasad model for $\GSpin_5\times \GSpin_4$. The representation $\rho_X$ in this case is the tensor product of the standard representations of the two $\GL_2(\BC)$ copies with the Spin representation of $\GSpin_5(\BC)$.

We can define the character of the component group $\omega_{\phi,H}$ by the same formula as all the cases in Table \ref{fig:1}. This allows us to formulate the epsilon dichotomy conjecture for this model as in all the cases in Table \ref{fig:1}. We can also formulate the weak form of the conjecture as in Conjecture \ref{weak conjecture}.

\begin{conj}\label{weak conjecture GSp_4}
Let $\Pi_\phi$ be a tempered $L$-packet whose central character is trivial on $Z_{G,H}(F)$. The unique distinguished element in $\Pi_\phi$ for the model $(\GSp_4\times \GL_2\times \GL_2,(\GL_2\times \GL_2)^0)$ belongs to $\Pi_\phi(G)$ (resp. $\Pi_\phi(G_D)$) if and only if $\epsilon(\frac{1}{2},\Pi_\phi,\rho_X)=1$ (resp. $\epsilon(\frac{1}{2},\Pi_\phi,\rho_X)=-1$). 
\end{conj}

Lastly, we discuss the multiplicity formula of this model. For a quasi-character $\theta$ (resp. $\theta_D$) of $G(F)$ (resp. $G_D(F)$), we define the geometric multiplicity by the formula \eqref{GSpin multiplicity 1} (resp. \eqref{GSpin multiplicity 2}). The multiplicity formulas
$$m(\pi)=m_{geom}(\pi),\;m(\pi_D)=m_{geom}(\pi_D)$$
for tempered representations can be proved by the same argument as the orthogonal Gan--Gross--Prasad model case. Moreover, the multiplicity formulas imply that the summation of the multiplicities is equal to 1 over every tempered local $L$-packet.

\begin{rmk}
Our argument in this paper can also be applied to the model $(\GSp_4\times \GL_2\times \GL_2,(\GL_2\times \GL_2)^0)$. It proves the epsilon dichotomy conjecture when the packet is not discrete with $|\Pi_\phi(G)|=1$. On the other hand, when the packet is discrete with $|\Pi_\phi(G)|=1$, if we assume the packet has trivial central character, then the epsilon dichotomy conjecture follows from the epsilon dichotomy of the Gan-Gross-Prasad model $(\SO_5\times \SO_4,\SO_4)$ proved in \cite{Wal3}.
\end{rmk}

\subsection{The proof of Theorem \ref{main theorem} and \ref{thm weak conjecture smaller models} for $(\GSp_6\times \GSp_4, (\GSp_4\times \GSp_2)^0)$}\label{sec GSp 4}
In this subsection we will prove Theorem \ref{main theorem} for the model $(\GSp_6\times \GSp_4, (\GSp_4\times \GSp_2)^0)$. The argument is very similar to the four models in the previous section, we will only give a sketch of the proof. Let $\Pi_\phi=\Pi_\phi(G)\cup \Pi_\phi(G_D)$ be a tempered $L$-packet whose central character is trivial on $Z_{G,H}(F)$. We assume that $\Pi_\phi$ is not discrete with $|\Pi_\phi(G)|=1$.

The first step is still to prove that the unique distinguished element belongs to $\Pi_\phi(G)$ if and only if $\epsilon(\frac{1}{2},\Pi_\phi,\rho_X)=1$. There are two cases.

The first case is when the packet is induced from a maximal parabolic subgroup $M$ of $G$. If $M$ does not contain the Levi subgroup $L$ of $G$ that is isomorphic to $(\GL_2\times \GL_2)\times (\GL_2\times \GL_1)$, Proposition \ref{GSp(6)xGSp(4) parabolic induction} implies that the unique distinguished element belongs to $\Pi_\phi(G)$. It is also easy to see that $\epsilon(\frac{1}{2},\Pi_\phi,\rho_X)=1$ in this case.

If $M=\GSp_6\times (\GL_2\times \GL_1)$ (resp. $M=(\GL_2\times \GL_2)\times \GSp_4$), then the statement follows from Proposition \ref{GSp(6)xGSp(4) parabolic induction} and Conjecture \ref{weak conjecture} for the model $(\GSp_6\times \GL_2,\GL_2\ltimes U)$ (resp. Conjecture \ref{weak conjecture GSp_4} for the model $(\GSp_4\times \GL_2\times \GL_2,(\GL_2\times \GL_2)^0)$). Note that both models are smaller than $(G,H)$. This also proves Theorem \ref{thm weak conjecture smaller models} for the model $(\GSp_6\times \GSp_4, (\GSp_4\times \GSp_2)^0)$.

The second case is when the packet $\Pi_\phi(G)$ is discrete. By our assumption, we must have $|\Pi_\phi(G)|>1$. Hence there exists a proper elliptic extended endoscopic triple $(G',s',{}^L\eta)$ of $G$ such that $\phi$ factors through ${}^L\eta$ and $s'\in Z_\phi$. We may also assume that $s'=(s_1,1)$ or $(1,s_2)$. If the order of $s'$ is equal to 4, then $G'=\GSO_6\times \GSp_4/\GL_{1}^{\diag}$ and Proposition \ref{prop GSp(6)xGSp(4)} implies that 
$$\sum_{\pi\in \Pi_\phi(G)} \tr(\chi_\pi(s'))m(\pi)=c_{\theta_{\Pi_\phi(G')}}(1)=1,$$
$$\sum_{\pi_D\in \Pi_\phi(G_D)} \tr(\chi_{\pi_D}(s'))m(\pi_D)=0,$$
i.e. the unique distinguished element belongs to $\Pi_\phi(G)$. In this case, by our discussion in Section \ref{sec epsilon factor}, we also know that $\epsilon(\frac{1}{2},\Pi_\phi,\rho_X)=1$.

If $G'=\GSp_6\times \GSO_4/\GL_{1}^{\diag}$ (resp. $G(\Sp_2\times \SO_4)\times \GSp_4/\GL_{1}^{\diag}$), in Section \ref{sec epsilon factor}, we have decomposed $\rho_X\circ\phi$ into $\rho_{s',\phi,+}\oplus \rho_{s',\phi,-}$. By Proposition \ref{prop GSp(6)xGSp(4)} and the same argument as in the previous section, we have
$$\sum_{\pi\in \Pi_\phi(G)} \tr(\chi_\pi(s'))m(\pi)=\frac{\epsilon(\frac{1}{2},\rho_{s',\phi,-})+\epsilon(\frac{1}{2},\rho_{s',\phi,+})}{2},$$
$$\sum_{\pi_D\in \Pi_\phi(G_D)} \tr(\chi_{\pi_D}(s'))m(\pi_D)=\frac{\epsilon(\frac{1}{2},\rho_{s',\phi,-})-\epsilon(\frac{1}{2},\rho_{s',\phi,+})}{2}.$$
Here we need to use Conjecture \ref{weak conjecture} for the model $(\GSp_6\times \GL_2,\GL_2\ltimes U)$ (resp. Conjecture \ref{weak conjecture GSp_4} for the model $(\GSp_4\times \GL_2\times \GL_2,(\GL_2\times \GL_2)^0)$). As a result, we know that the unique distinguished element belongs to $\Pi_\phi(G)$ if and only if $\epsilon(\frac{1}{2},\rho_{s',\phi,+})=\epsilon(\frac{1}{2},\rho_{s',\phi,-})$ which is equivalent to $\epsilon(\frac{1}{2},\Pi_\phi,\rho_X)=1$.

Now we can prove the theorem. Let $\omega_\phi\in \hat{S}_\phi$ corresponds to the unique distinguished element in the packet. By Remark \ref{distinguished is character} we know that $\omega_\phi$ is a character and we view it as a character of $Z_\phi$. For $s\in S_\phi$, by Lemma \ref{lem extended endoscopic triple}, there exists an elliptic extended endoscopic triple $(G',s',{}^L\eta)$ of $G/Z_{G,H}$ such that $s'\in sZ_{\phi}^{\circ}$ and $\phi$ factors through ${}^L\eta$. We need to show that $\omega_\phi(s')=\omega_{\phi,H}(s)$. The above discussion implies that $\omega_\phi(s')=\omega_{\phi,H}(s)$ if $s'$ belongs to the center of the dual group.

For general $s'=(s_1,s_2)$ with $s_1\in \Spin_7(\BC)$ and $s_2\in \Spin_5(\BC)$, by our definition of $\omega_{\phi,H}$ we know that 
$$\omega_{\phi,H}(s_1,s_2)=\omega_{\phi,H}(s_1,1)\omega_{\phi,H}(1,s_2).$$
Hence it is enough to consider the case when $s'=(s_1,1)$ or $s'=(1,s_2)$.

If the order of $s'$ is equal to 4, by the discussion above we know that the unique distinguished element belongs to $\Pi_\phi(G)$. By the definition of $\omega_{\phi,H}$ we know that $\omega_{\phi,H}(s)=1$. This implies that
$$\omega_\phi(s')=\tr(\omega_\phi(s'))=\sum_{\pi\in \Pi_\phi(G)} \tr(\chi_\pi(s'))m(\pi)=1=\omega_{\phi,H}(s).$$

If the order of $s'$ is equal to 2, by our discussion above, we have 
$$\sum_{\pi\in \Pi_\phi(G)} \tr(\chi_\pi(s'))m(\pi)=\frac{\epsilon(\frac{1}{2},\rho_{s',\phi,+})+\epsilon(\frac{1}{2},\rho_{s',\phi,-})}{2},$$
$$\sum_{\pi_D\in \Pi_\phi(G_D)} \tr(\chi_{\pi_D}(s'))m(\pi_D)=\frac{-\epsilon(\frac{1}{2},\rho_{s',\phi,+})+\epsilon(\frac{1}{2},\rho_{s',\phi,-})}{2}.$$
By the definition in Section \ref{sec epsilon factor}, we have
$$\omega_{\phi,H}(s)=\epsilon(\frac{1}{2},\rho_{s',\phi,-}).$$

We have two cases. If the unique distinguished element belongs to $\Pi_\phi(G)$, we have
$$\epsilon(\frac{1}{2},\Pi_\phi,\rho_X)=1,\;\epsilon(\frac{1}{2},\rho_{s',\phi,+})=\epsilon(\frac{1}{2},\rho_{s',\phi,-}).$$
This implies that
$$\omega_\phi(s')=\tr(\omega_\phi(s'))=\sum_{\pi\in \Pi_\phi(G)} \tr(\chi_\pi(s'))m(\pi)$$
$$=\frac{\epsilon(\frac{1}{2},\rho_{s',\phi,+})+\epsilon(\frac{1}{2},\rho_{s',\phi,-})}{2}=\epsilon(\frac{1}{2},\rho_{s',\phi,-}).$$

If the unique distinguished element belongs to $\Pi_\phi(G_D)$, we have
$$\epsilon(\frac{1}{2},\Pi_\phi,\rho_X)=-1,\;\epsilon(\frac{1}{2},\rho_{s',\phi,+})=-\epsilon(\frac{1}{2},\rho_{s',\phi,-}).$$
This implies that
$$\omega_\phi(s')=\tr(\omega_\phi(s'))=\sum_{\pi_D\in \Pi_\phi(G_D)} \tr(\chi_{\pi_D}(s'))m(\pi_D)$$
$$=\frac{-\epsilon(\frac{1}{2},\rho_{s',\phi,+})+\epsilon(\frac{1}{2},\rho_{s',\phi,-})}{2}=\epsilon(\frac{1}{2},\rho_{s',\phi,-}).$$

This finishes the proof of Theorem \ref{main theorem} for the model $(\GSp_6\times \GSp_4,G(\Sp_4\times \Sp_2))$.

\section{The model $(E_7,\PGL_2\ltimes U)$}\label{sec E7}

\subsection{The model and the multiplicity formula}
Let $G=E_7$ be the split adjoint reductive group of Type $E_7$, and let $P=LU$ be the parabolic subgroup of $G$ of Type $A_1\times A_1\times A_1$ defined in Section 7 of \cite{WZ2}. Let $\xi:U(F)\rightarrow \BC^\times$ be the generic character defined in loc. cit. and let $H_0\subset  L$ be the stabilizer of $\xi$ which is isomorphic to $\PGL_2$. Let $H=H_0\ltimes U$ and we extend the character $\xi$ to $H(F)$ by making it trivial on $H_0(F)$. We can also define the quaternion version $(G_D,H_D=H_{0,D}\ltimes U_D,\xi_D)$ where $G_D$ is the unique pure inner form of $G$ ($G_D$ has split rank 4) and $H_{0,D}(F)\simeq \PGL_1(D)$. We refer the reader to Section 7 of \cite{WZ2} for more details of this model. Let $\pi$ (resp. $\pi_D$) be an irreducible representation of $G(F)$ (resp. $G_D(F)$), we define the multiplicities
$$m(\pi)=\dim(\Hom_{H(F)}(\pi,\xi)),\;m(\pi_D)=\dim(\Hom_{H_D(F)}(\pi_D,\xi_D)).$$

For the multiplicity formula, let $\theta$ (resp. $\theta_D$) be a quasi-character of $G(F)$ (resp. $G_D(F)$). Define the geometric multiplicities 
$$m_{geom}(\theta)=c_\theta(1)+\sum_{T\in \CT_{ell}(H_0)} \frac{1}{2}\int_{T(F)/Z_{G,H}(F)}D^H(t) c_{\theta}(t)dt,$$
$$m_{geom}(\theta_D)=\sum_{T_D\in \CT_{ell}(H_{0,D})}\frac{1}{2} \int_{T_D(F)/Z_{G_D,H_D}(F)} D^{H_D}(t) c_{\theta_D}(t)dt.$$
For the rest of this section we will assume that the multiplicity formulas
$$m(\pi)=m_{geom}(\theta_\pi),\;m(\pi_D)=m_{geom}(\theta_{\pi_D})$$
hold for all tempered representations $\pi$ (resp. $\pi_D$) of $G(F)$ (resp. $G_D(F)$). 

To end this subsection, we will discuss the behavior of the geometric multiplicities under parabolic induction. 
Let $M$ be a proper Levi subgroup of $G$ and $\theta^M$ be a quasi-character on $M(F)$. If $M$ does not contain the Levi subgroup $L$ up to conjugation, define $m_{geom}(\theta^M)=c_{\theta^M}(1)$. Otherwise, $M$ corresponds to a proper Levi subgroup $M_D$ of $G_D$. Moreover, up to conjugation we may assume that $L\subset M$ and $L_D\subset M_D$. Let $\theta_{D}^{M_D}$ be a quasi-character on $M_D(F)$. Define
\begin{eqnarray*}
m_{geom}(\theta^M)&=&c_{\theta^M}(1)+\sum_{T\in \CT_{ell}(H_0)}\frac{1}{2} \int_{T(F)/Z_{G,H}(F)} \\
&&  D^M(t)^{1/2}(t)D^{H_0}(t)^{-1/2} c_{\theta^M}(t)dt,\\
m_{geom}(\theta_{D}^{M_D})&=&\sum_{T_D\in \CT_{ell}(H_{0,D})} \frac{1}{2} \int_{T_D(F)/Z_{G_D,H_D}(F)}\\
&&D^{M_D}(t)^{1/2}D^{H_{0,D}}(t)^{-1/2} c_{\theta_{D}^{M_D}}(t)dt.
\end{eqnarray*}
The following proposition is a direct consequence of Proposition \ref{germ parabolic induction} (one just need to use the fact that $D^H(t)=D^G(t)^{1/2}D^{H_0}(t)^{-1/2}$ for $t\in H_0(F)$). 

\begin{prop}\label{E7 parabolic induction}
Let $\theta$ (resp. $\theta_D$) be a quasi-character on $G(F)$ (resp. $G_D(F)$). Assume that $\theta$ (resp. $\theta_D$) is the parabolic induction of a quasi-character $\theta^M$ (resp. $\theta_{D}^{M_D}$) of a proper Levi subgroup $M$ of $G$ (resp. $M_D$ of $G_D$). We have
$$m_{geom}(\theta)=m_{geom}(\theta^M),\;m_{geom}(\theta_D)=m_{geom}(\theta_{D}^{M_D}).$$
\end{prop}

\subsection{The smaller models}
In this subsection we will discuss the models that are smaller than the model $(E_7,\PGL_2\ltimes U)$. There are three smaller models. The first one is $(\GSpin_{10}\times \GSpin_3,\GSpin_3\ltimes U)$ which is an analogy of the Gan--Gross--Prasad model $(\SO_{10}\times \SO_3,\SO_3\ltimes U)$. To be specific, let $Q=MN$ be the parabolic subgroup of $\GSpin_{10}$ associated to the simple roots $\frac{e_4\pm e_5}{2}$. We can define a generic character $\xi_N$ of $N(F)$ similar to the Gan--Gross--Prasad model case and its stabilizer in $M(F)$ is isomorphic to $\GSpin_3(F)=\GL_2(F)$. This defines the model $(\GSpin_{10}\times \GSpin_3,\GSpin_3\ltimes U)$.

To define the other two models, we need to use the group $\GHSpin_{4n}=\GSpin_{4n}/\{1,z\}$ where $z$ is an order 2 element in the center of $\GSpin_{4n}$ that does not belong to the connected component of the center (there are two such elements differed by the outer automorphism). Note that the map $\GSpin_{4n}(F)\rightarrow \GHSpin_{4n}(F)$ is not surjective. An example would be $\GHSpin_{4}\simeq \GL_2\times \PGL_2$. 

The center of the group $\GHSpin_{4n}$ is $\GL_1$, it has a unique Half-Spin representation, and it is equipped with a similitude character $l:\GHSpin_{4n}\rightarrow \GL_1$. We use $\GHSpin_{4n}^{\vee}$ to denote the dual group of $\GHSpin_{4n}$ and it is also equipped with a similitude character $l:\GHSpin_{4n}^{\vee}\rightarrow \GL_1(\BC)$ whose kernel is $\Spin_{4n}(\BC)$. Moreover, the group $\GHSpin_{4n}^{\vee}$ has two Half-Spin representations, one of them has determinant 1 and the other one has a nontrivial determinant. We use $\HSpin_{4n}^{+}$ (resp. $\HSpin_{4n}^{-}$) to denote the Half-Spin representation with determinant 1 (resp. nontrivial determinant.)

The two remaining smaller models of the model $(E_7,\PGL_2\ltimes U)$ are related to the group $\GHSpin_{12}$. Another way to describe the group $\GHSpin_{12}$ is that it is the Levi subgroup of the group $E_7$ of Type $D_6$. Similarly, its dual group $\GHSpin_{12}^{\vee}$ is the Levi subgroup of the group $E_{7,sc}(\BC)$ of Type $D_6$ and we have $\GHSpin_{12}^{\vee}\simeq \Spin_{12}(\BC)\times \GL_1(\BC)/\{1,(z,-1)\}$. Under this isomorphism, the $\HSpin_{12}^{+}$ representation is just a Half-Spin representation of $\Spin_{12}$ and the $\HSpin_{12}^{-}$ representation is a Half-Spin representation of $\Spin_{12}$ tensor with the standard representation of $\GL_1$.

One of the reduced model is an analogy of the model $(\GSO_{12},\GL_2\ltimes U)$ for the group $\GHSpin_{12}$. Let $(G,H=H_0\ltimes U,\xi)=(E_7,\PGL_2\ltimes U,\xi)$ be the model defined in the previous subsection. Recall that we also have the Levi subgroup $P=LU$. Let $Q=MN$ be the Levi subgroup of $G=E_7$ with $P\subset Q$ and $M\simeq \GHSpin_{12}$ and let $M_H=M\cap H=H_0\ltimes (U\cap M)$. The first smaller model is just $(M,M_H,\xi|_{M_H})$. We will denote this model by $(\GHSpin_{12},\PGL_2\ltimes U)$.

\begin{rmk}
With the notation above, $U\cap M$ is the unipotent subgroup of the parabolic subgroup $P\cap M=L\ltimes (U\cap M)$ of $M$. The character $\xi|_{U\cap M}$ is a generic character whose centralizer in $L(F)$ is $H_0(F)\times Z_M(F)\simeq \PGL_2(F)\times \GL_1(F)$.
\end{rmk}

The other one is an analogy of the Gan--Gross--Prasad model $(\SO_{12}\times \SO_3,\SO_3\ltimes U)$. To be specific, let $Q=MN$ be the parabolic subgroup of $\GHSpin_{12}$ associated to the simple roots $\frac{e_5\pm e_6}{2}$. We can define a generic character $\xi_N$ of $N(F)$ similar to the Gan--Gross--Prasad model case whose stabilizer in $M(F)$ is isomorphic to $\GSpin_3(F)=\GL_2(F)$. This defines the model $(\GHSpin_{12}\times \GSpin_3,\GSpin_3\ltimes N,\xi_N)$. We will denote this model by $(\GHSpin_{12}\times \GSpin_3,\GSpin_3\ltimes U)$.

Let $(G,H=H_0\ltimes U,\xi)$ be one of the smaller models above. We can also define the quaternion version of the model in a similar way. We will denote it by $(G_D,H_D=H_{0,D}\ltimes U_D,\xi_D)$. Let $\pi$ (resp. $\pi_D$) be an irreducible representation of $G(F)$ (resp. $G_D(F)$) whose central character is trivial on $Z_{G,H}(F)=Z_G(F)\cap H(F)$ (resp. $Z_{G_D,H_D}(F)=Z_{G_D}(F)\cap H_D(F)$), we define the multiplicities
$$m(\pi)=\dim(\Hom_{H(F)}(\pi,\xi)),\;m(\pi_D)=\dim(\Hom_{H_D(F)}(\pi_D,\xi_D)).$$

For the multiplicity formula, let $\theta$ (resp. $\theta_D$) be a quasi-character of $G(F)$ (resp. $G_D(F)$). Define the geometric multiplicities 

$$m_{geom}(\theta)=c_\theta(1)+\sum_{T\in \CT_{ell}(H_0)} \frac{1}{2}\int_{T(F)/Z_{G,H}(F)}D^H(t) c_{\theta}(t)dt,$$
$$m_{geom}(\theta_D)=\sum_{T_D\in \CT_{ell}(H_{0,D})}\frac{1}{2}\int_{T_D(F)/Z_{G_D,H_D}(F)} D^{H_D}(t) c_{\theta_D}(t)dt.$$
We will assume that the multiplicity formulas
$$m(\pi)=m_{geom}(\theta_\pi),\;m(\pi_D)=m_{geom}(\theta_{\pi_D})$$
hold for all tempered representations. 

\begin{rmk}
The multiplicity formula for the models $(\GHSpin_{12}\times \GSpin_3,\GSpin_3\ltimes U)$ and $(\GSpin_{10}\times \GSpin_3,\GSpin_3\ltimes U)$ can be proved by a similar argument as the Gan--Gross--Prasad model case. When $F$ is $p$-adic, the multiplicity formula for the model $(\GHSpin_{12},\PGL_2\ltimes U)$ can be proved by a similar argument as the model $(\GSO_{12},\GL_2\ltimes U)$.
\end{rmk}

Like in all the previous cases, combining the   multiplicity formula and the local Langlands correspondence, we know that each tempered $L$-packet contains a unique distinguished element, and the unique distinguished element corresponds to a character of the component group.

To formulate the weak conjecture for the smaller models, we need to define the representation $\rho_X$ of the dual group. If the model is $(\GSpin_{10}\times \GSpin_3,\GSpin_3\ltimes U)$, we let $\rho_X$ be the 30-dimensional tensor product $L$-function of $\GSO_{10}(\BC)\times \GSp_2(\BC)$. If the model is 
$(\GHSpin_{12},\PGL_2\ltimes U)$, we let $\rho_X$ be the representation $\HSpin_{12}^{+}$ of $\GHSpin_{12}^{\vee}$. If the model is $(\GHSpin_{12}\times \GL_2,\GL_2\ltimes U)$, let $\rho_X$ be the tensor product of the 12-dimensional standard representation of $\GHSpin_{12}^{\vee}$ with the 2-dimensional standard representation of $\GL_2(\BC)$. Let $\Pi_\phi=\Pi_\phi(G)\cup \Pi_\phi(G_D)$ be a tempered $L$-packet whose central character is trivial on $Z_{G,H}(F)$. We can formulate the weak conjecture in this case.

\begin{conj}\label{weak conjecture for E7 smaller model}
The unique distinguished element in $\Pi_\phi$ belongs to $\Pi_\phi(G)$ (resp. $\Pi_\phi(G_D)$) if and only if 
$$\epsilon(\frac{1}{2},\Pi_\phi,\rho_X)=1,\; \text{(resp.}\; \epsilon(\frac{1}{2},\Pi_\phi,\rho_X)=-1).$$ 
\end{conj}

\begin{rmk}
We can also formulate the epsilon dichotomy conjecture for these smaller models.
\end{rmk}

\subsection{The endoscopic relation}
In this subsection we will study the behavior of the geometric multiplicity under endoscopy. Let $G=E_7$ and $G_D$ be its pure inner form. Let $\theta$ (resp. $\theta_D$) be a quasi-character of $G(F)$ (resp. $G_D(F)$). Recall that we have defined the geometric multiplicities 
$$m_{geom}(\theta)=c_\theta(1)+\sum_{T\in \CT_{ell}(H_0)} \frac{1}{2}\int_{T(F)/Z_{G,H}(F)}D^H(t) c_{\theta}(t)dt,$$
$$m_{geom}(\theta_D)=\sum_{T_D\in \CT_{ell}(H_{0,D})}\frac{1}{2} \int_{T_D(F)/Z_{G_D,H_D}(F)} D^{H_D}(t) c_{\theta_D}(t)dt.$$

Let $(G',s',{}^L\eta)$ be a proper elliptic extended endoscopic triple of $G$, and let $\theta'$ be a stable quasi-character of $G'(F)$. Assume that $\theta$ (resp. $\theta_D$) is the endoscopic transfer of $\theta'$. To define $m_{geom}(\theta')$ and $m_{geom,D}(\theta')$, we have 4 situations. Note that like in the previous cases, we always choose ${}^L\eta$ to be the natural embedding from ${}^LG'$ into ${}^LG$.

If $\hat{G}'=\SL_8(\BC)/\BZ_2$, we let 
$$m_{geom}(\theta')=c_{\theta'}(1),\;m_{geom,D}(\theta')=0.$$ 
If 
\begin{eqnarray*}
\hat{G}'&=&\SL_6(\BC)\times \SL_3(\BC)/\BZ_3=\{(g_1,g_2)\in \GL_6(\BC)\times \GL_3(\BC) \mid \\
&&\det(g_1)=\det(g_2)^4\}/\{(a^2I_6,aI_3) \mid a\in \GL_1(\BC)\},
\end{eqnarray*}
we have
$$G'=\{(g_1,g_2)\in \GL_6\times \GL_3 \mid  \det(g_1)^2=\det(g_2)\}/\{(aI_6,a^4I_3) \mid a\in \GL_1\}.$$
We can embed $\PGL_2$ into $G'$ via the map $h\mapsto (\diag(h,h,h)\times \det(h)^2I_3)$ and we will denote this embedding by $\nu$. We define
$$m_{geom}(\theta')=c_{\theta'}(1)+\sum_{T\in \CT_{ell}(\PGL_2)}\frac{1}{2} \int_{T(F)} D^{\PGL_2}(t)^4 c_{\theta'}(\nu(t))  dt,$$
$$m_{geom,D}(\theta')=\varepsilon(s')\sum_{T\in \CT_{ell}(\PGL_2)}\frac{1}{2} \int_{T(F)} D^{\PGL_2}(t)^4 c_{\theta'}(\nu(t))  dt$$
where $\varepsilon(s')$ is equal to 1 if the order of $s'$ is 3 and it is equal to $-1$ if the order of $s'$ is 6.

If 
\begin{eqnarray*}
\hat{G}'&=&\SL_4(\BC)\times \SL_4(\BC)\times \SL_2(\BC)/\BZ_4 \\
&=&\{(g_1,g_2,g_3)\in \GL_4(\BC)\times \GL_4(\BC)\times \GL_2(\BC) \mid \det(g_1)=\det(g_2)\\
&&=\det(g_3)^{-1}\}/\{(aI_4,aI_4,a^{-2}I_2) \mid a\in \GL_1(\BC)\},    
\end{eqnarray*}
we have
$$G'=\{(g_1,g_2,g_3)\in \GL_4\times \GL_4\times \GL_2 \mid  \det(g_1)\det(g_2)$$
$$=\det(g_3)^2\}/\{(aI_4,bI_4,abI_2) \mid a,b\in \GL_1\}.$$
In this case, $s'$ is equal to $(I_4,\pm iI_4,\pm I_2)$.
We have two embeddings $\nu_1,\nu_2$ from $\PGL_2$ into $G'$ given by
$$\nu_1(h)=(\diag(h,h),I_4,h),\;\nu_2(h)=(I_4,\diag(h,h),h).$$
We define
$$m_{geom}(\theta')=c_{\theta'}(1)+\sum_{i=1}^{2}\sum_{T\in \CT_{ell}(\PGL_2)}\frac{1}{2} \int_{T(F)} D^{\PGL_2}(t)^2 c_{\theta'}(\nu_i(t))  dt,$$
$$m_{geom,D}(\theta')=\varepsilon(s')\sum_{i=1}^{2}(-1)^i\sum_{T\in \CT_{ell}(\PGL_2)}\frac{1}{2}\int_{T(F)} D^{\PGL_2}(t)^2 c_{\theta'}(\nu_i(t))  dt$$
where $\varepsilon(s')$ is equal to $-1$ if $s'=(I_4,\pm iI_4,I_2)$ and it is equal to $1$ if $s'=(I_4,\pm iI_4,-I_2)$.

The last case is when 
\begin{eqnarray*}
\hat{G'}&=&\Spin_{12}(\BC)\times \SL_2(\BC)/\BZ_2\\
&=&\{(g_1,a,g_2) \mid (g_1,a)\in \Spin_{12}(\BC)\times \GL_1(\BC)/\{1,(z,-1)\},\\
&&g_2\in \GL_2(\BC),\;\det(g_2)=a^{-2}\}/\{(1,a,aI_2) \mid a\in \GL_1(\BC)\}\\
&=&\{(g_1,g_2)\in \GHSpin_{12}^{\vee}\times \GL_2(\BC) \mid \\
&&l(g_1)=\det(g_2)^{-1}\}/\GL_1(\BC)^{anti-diag}.
\end{eqnarray*}
There are two choices of $s'$, the $-1$ eigenspace of $\rho_X(s')$ (here $\rho_X$ is the 56-dimensional representation of $E_{7,sc}(\BC)$) is either 24 dimensional or 32 dimensional depends on the choice of $s'$. In this case, we have
$$G'=\{(g_1,g_2)\in \GHSpin_{12}\times \GL_2 \mid  l(g_1)=\det(g_2)\}/\GL_{1}^{\diag}.$$
By our discussion of the smaller model $(\GHSpin_{12},\PGL_2\ltimes U)$ in the previous subsection, we have an embedding from $\PGL_2$ into $\GHSpin_{12}$ which induces an embedding from $\PGL_2$ into $G$ by making it trivial on the $\GL_2$-component. We denote this embedding by $\nu_1$. By our discussion of the smaller model $(\GHSpin_{12}\times \GL_2,\GL_2\ltimes U)$, we have a diagonal embedding from $\GL_2$ into $\GHSpin_{12}\times \GL_2$ which induces an embedding from $\PGL_2$ into $G$. We denote this embedding by $\nu_2$. We define
\begin{eqnarray*}
m_{geom}(\theta')&=&c_{\theta'}(1)+\sum_{i=1}^{2}\sum_{T\in \CT_{ell}(\PGL_2)}\frac{1}{2}\\
&&\cdot\int_{T(F)} D^{\PGL_2}(t)^{-1/2}D^{G'}(\nu_i(t))^{1/2} c_{\theta'}(\nu_i(t))  dt,
\end{eqnarray*}
\begin{eqnarray*}
m_{geom,D}(\theta')&=&\varepsilon(s')\sum_{i=1}^{2}(-1)^i\sum_{T\in \CT_{ell}(\PGL_2)}\frac{1}{2} \\
&&\cdot \int_{T(F)} D^{\PGL_2}(t)^{-1/2}D^{G'}(\nu_i(t))^{1/2} c_{\theta'}(\nu_i(t))  dt
\end{eqnarray*}
where $\varepsilon(s')$ is equal to $-1$ (resp. 1) if the $-1$ eigenspace of $\rho_X(s')$ is 24-dimensional (resp. 32-dimensional).

\begin{prop}\label{prop E7}
Let $\theta$ (resp. $\theta_D$) be a quasi-character on $G(F)$ (resp. $G_D(F)$). Assume that $\theta$ (resp. $\theta_D$) is the endoscopic transfer of a stable quasi-character $\theta'$ of $G'(F)$ . We have
$$m_{geom}(\theta)=m_{geom}(\theta'),\;m_{geom}(\theta_D)=m_{geom,D}(\theta').$$
\end{prop}

\begin{proof}
We will only prove the case when $\hat{G'}=\Spin_{12}(\BC)\times \SL_2(\BC)/\BZ_2$. The rest case follows from a similar argument. Like in all the previous cases, the only difference between the split case and the quaternion case is the extra sign in the transfer factor. Hence we will only consider the split case.

The proof of the equation $c_\theta(1)=c_{\theta'}(1)$ is easy and we will skip it here. We fix a quadratic extension $E/F$ and let $T_E\in \CT_{ell}(\PGL_2)=\CT_{ell}(H_0)$ correspond to $E$ (for simplicity we identify $H_0$ with $\PGL_2$). We just need to show that the term corresponds to $T_E$ in $m_{geom}(\theta)$ is equal to the term corresponds to $T_E$ in $m_{geom}(\theta')$.  

Let $T_G$ be the centralizer of $T_E$ in the Levi subgroup $L$ of $G$, which is a maximal torus of $G$. On the other hand, let $L_1$ be the Levi subgroup of $G'$ which is of Type $A_1\times A_1\times A_1$ on the $\GHSpin_{12}$-copy and is a maximal torus on the $\GL_2$-copy such that it contains $\nu_1(T_E)$. Let $L_2$ be the Levi subgroup of $G'$ which is of Type $D_2$ on the $\GHSpin_{12}$-copy and is equal to $\GL_2$ on the $\GL_2$-copy such that it contains $\nu_2(T_E)$. Let $T_{G',i}$ be the centralizer of $\nu_i(T_E)$ in $L_i$, which is a maximal torus of $G'$. 

Let $W=W(G,T_G)$ and $W_i=W(G',T_{G',i})$ be the Weyl groups. We have $$|W|=9216,\;|W_1|=768,\;|W_2|=1536$$ 
and $W_i$ can be naturally identified as a subgroup of $W$ for $i=1,2$. Note that $W$ (resp. $W_1$, $W_2$) is of Type 
$$F_4\times (A_1)^3\; \text{(resp.} \;C_3\times (A_1)^4, \;D_4\times (A_1)^3).$$ 
The $W$-action stabilizes $T_E$ and its action on $T_E$ factors through the Weyl group $W(\PGL_2,T_E)\simeq \BZ/2\BZ$. It is easy to see that there are natural isomorphisms $f_i:T_G(F)\simeq T_i(F)$ whose restriction to $T_E(F)$ are the identity map (here by abusing of notation we identify $T_E$ with $\nu_i(T_E)$) and satisfy the following condition:
\begin{itemize}
\item for $\gamma\in T_G(F)\cap G_{reg}(F)$, there are exactly $$18=12+6=\frac{9216}{768}+\frac{9216}{1536}$$ 
conjugacy classes $\gamma_{G'}$ of $G'(F)$ with $\Delta(\gamma,\gamma_{G'})\neq 0$. Each of them is represented by an element $f_i(w\gamma w^{-1})$ for $1\leq i\leq 2$ and $w\in W/W_i$.
\end{itemize}

Next we show that the transfer factor $\Delta(\gamma,\gamma_{G'})$ is always equal to 1 for any $\gamma\in T_G(F)\cap G_{reg}(F)$ and $\gamma_{G'}=f_i(w\gamma w^{-1})$. We follow the notation in Section 3 of \cite{LS}. It is easy to see that in this case $\textbf{s}_{T_G}'\in \pi_0(\hat{T}_{G,ad}^{\Gamma})$ is the identity component. This implies that the terms $\Delta_{I}(\gamma,\gamma_{G'})$ and $\Delta_{III_{1}}(\gamma,\gamma_{G'})$ are equal to 1. Also for any root $\alpha$ of $T_G$ outside $G'$, we have $F_{\pm \alpha}=F_\alpha$ and hence we can choose the $\chi$-data $\chi_\alpha$ to be the trivial character. This implies that $\Delta_{II}(\gamma,\gamma_{G'})=1$. Lastly, it is easy to see that $\textbf{a}\in H^1(W_F,\hat{T}_G)$ is the trivial cocycle (note that for any regular semisimple element $t\in T_G(F)$, the stable conjugacy class of $t$ only contains one rational conjugacy class). This implies that $\Delta_{III_2}(\gamma,\gamma_{G'})=1$. This proves that the transfer factor $\Delta(\gamma,\gamma_{G'})$ is always equal to 1 for any $\gamma\in T_G(F)\cap G_{reg}(F)$ and $\gamma_{G'}=f_i(w\gamma w^{-1})$.

For $t\in T_E(F)\cap H_{0,reg}(F)$, we have 
$$D^G(t)^{1/2}c_\theta(t)=\frac{1}{1152}\lim_{t'\in T_G(F)\cap G_{reg}(F)\rightarrow t} D^G(t')^{1/2}\theta(t')$$
where $1152$ is the cardinality of the Weyl group of $G_t(F)$ (which is of Type $F_4$). Similarly, for $t\in T_E(F)\cap \PGL_{2,reg}(F)$, we have 
$$D^{G'}(\nu_1(t))^{1/2}c_{\theta'}(\nu_1(t))=\frac{1}{96}\lim_{t'\in T_{G',1}(F)\cap G_{reg}'(F)\rightarrow \nu_1(t)} D^{G'}(t')^{1/2}\theta(t'),$$
$$D^{G'}(\nu_2(t))^{1/2}c_{\theta'}(\nu_2(t))=\frac{1}{192}\lim_{t'\in T_{G',2}(F)\cap G_{reg}'(F)\rightarrow \nu_2(t)} D^{G'}(t')^{1/2}\theta(t').$$
Here $96$ (resp. $192$) is the cardinality of the Weyl group of $(G')_{\nu_1(t)}(F)$ (resp. $(G')_{\nu_2(t)}(F)$), which is of Type $C_3\times A_1$ (resp. $D_4$).

Combining the above discussion, we know that  $D^H(t)c_\theta(t)$ is equal to 
$$D^{\PGL_2}(t)^{-1/2}(D^{G'}(\nu_1(t))^{1/2}c_{\theta'}(\nu_1(t))+D^{G'}(\nu_2(t))^{1/2}c_{\theta'}(\nu_2(t)))$$
for all $t\in T_E(F)\cap H_{0,reg}(F)=T_E(F)\cap \PGL_{2,reg}(F)$. Here we have used the identity $D^H(t)=D^{\PGL_2}(t)^{-1/2}D^G(t)^{1/2}$. Hence the term corresponds to $T_E$ in $m_{geom}(\theta)$ is equal to the term corresponds to $T_E$ in $m_{geom}(\theta')$. This proves the proposition. 
\end{proof}

\subsection{The main result and the proof}
In this subsection we are going to state and prove our main results for the model $(E_7,\PGL_2\ltimes U)$. Let $G=E_7$, $\phi:W_F\rightarrow \hat{G}$ be a tempered Langlands parameter, and $\Pi_\phi=\Pi_\phi(G)\cup \Pi_\phi(G_D)$ be the associated tempered $L$-packet. Like in all the other cases, we assume that the local Langlands correspondence holds for $G$.

\begin{thm}\label{main theorem for E7}
Assume that Conjecture \ref{weak conjecture for E7 smaller model} holds.  If the packet $\Pi_\phi$ is not discrete with $|\Pi_\phi(G)|=1$, then Conjecture \ref{main conj} holds for packet $\Pi_\phi$.
\end{thm}

\begin{cor}
Conjecture \ref{main conj} holds when $F=\BR$.
\end{cor}

\begin{thm}\label{thm weak conj for E7}
Assume that the Conjecture \ref{weak conjecture} holds for the model $(E_7,\PGL_2\ltimes U)$. Then Conjecture \ref{weak conjecture for E7 smaller model} holds.
\end{thm}

\begin{cor}
Conjecture \ref{weak conjecture} is equivalent to Conjecture \ref{main conj} for the model $(E_7,\PGL_2\ltimes U)$.
\end{cor}

By using Propositions \ref{E7 parabolic induction} and \ref{prop E7}, the proof of Theorem \ref{main theorem for E7} and \ref{thm weak conj for E7} is almost the same as all the previous cases. We will only give a sketch of the proof. 

The first step is to prove Conjecture \ref{weak conjecture} when $\Pi_\phi$ is not discrete with $|\Pi_\phi(G)|=1$. When $\Pi_\phi$ is not discrete, it is induced from a maximal parabolic subgroup $M$ of $G$. In this case, if $M$ does not contain $L$ up to conjugation, then Proposition \ref{E7 parabolic induction} implies that the unique distinguished element belongs to $\Pi_\phi(G)$. Also in this case it is easy to see that the epsilon factor $\epsilon(\frac{1}{2},\Pi_\phi,\rho_X)$ is equal to 1.

If $M$ is of Type $D_6$, then then Conjecture \ref{weak conjecture} follows from Proposition \ref{E7 parabolic induction}, Conjecture \ref{weak conjecture for E7 smaller model} for the model $(\GHSpin_{12},\PGL_2\ltimes U)$ and the multiplicity formula for the model $(\GHSpin_{12},\PGL_2\ltimes U)$. 

\begin{rmk}\label{E7 imply D6 GL(6)}
As in the previous cases, the above discussion also implies that Conjecture \ref{weak conjecture} for the model $(E_7,\PGL_2\ltimes U)$ would imply Conjecture \ref{weak conjecture for E7 smaller model} for the model $(\GHSpin_{12},\PGL_2\ltimes U)$. Similarly, Conjecture \ref{weak conjecture for E7 smaller model} for the model $(\GHSpin_{12},\PGL_2\ltimes U)$ would imply Conjecture \ref{conj GL(6) general} (we just need to consider the maximal Levi subgroup of $\GHSpin_{12}$ that is isomorphic to $\GL_6\times \GL_1$. By Remark \ref{GL(6) implies GL(4)}, it would also imply Conjecture \ref{conj GL(4)xGL(2)}.
\end{rmk}

If $M$ is of Type $D_5\times A_1$, Conjecture \ref{weak conjecture} follows from Proposition \ref{E7 parabolic induction}, Conjecture \ref{weak conjecture for E7 smaller model} for the model $(\GHSpin_{10}\times \GSpin_3,\GSpin_3\ltimes U)$ and the multiplicity formula for the model $(\GHSpin_{10}\times \GSpin_3,\GSpin_3\ltimes U)$. Note that in this case, $\hat{M}=\Spin_{10}(\BC)\times \SL_2(\BC)\times \GL_1(\BC)/(\BZ/4\BZ)$. We have a projection map 
\begin{eqnarray*}
&&\hat{M}=\Spin_{10}(\BC)\times \SL_2(\BC)\times \GL_1(\BC)/(\BZ/4\BZ) \\
&\rightarrow& \Spin_{10}(\BC)\times  \SL_2(\BC)/(\BZ/4\BZ)=\SO_{10}(\BC)\times  \SL_2(\BC)/(\BZ/2\BZ)\\
&&=(\GSO_{10}(\BC)\times \GL_2(\BC))^0/\GL_1(\BC).
\end{eqnarray*}
Combining with Theorem 8.1 of \cite{La}, each Langlands parameter of $M$ induces a Langlands parameter of $\GSpin_{10}\times \GSpin_3$, this allows us to apply Conjecture \ref{weak conjecture for E7 smaller model} for the smaller model $(\GSpin_{10}\times \GSpin_3,\GSpin_3\ltimes U)$. 

\begin{rmk}\label{E7 imply D5}
The above discussion also implies that Conjecture \ref{weak conjecture} for the model $(E_7,\PGL_2\ltimes U)$ would imply Conjecture \ref{weak conjecture for E7 smaller model} for the model $(\GSpin_{10}\times \GSpin_3,\GSpin_3\ltimes U)$. Note that by the above description of $\hat{M}$ and Theorem 8.1 of \cite{La}, a tempered L-packet of $\GSpin_{10}\times \GSpin_3$ whose central character is trivial on the diagonal $\GL_1$ would induce a $L$-packet of $M$.
\end{rmk}

If $M$ is of Type $A_5\times A_1$, Conjecture \ref{weak conjecture} follows from Proposition \ref{E7 parabolic induction}, Conjecture \ref{conj GL(6) general} (see Remark \ref{E7 imply D6 GL(6)}), and the multiplicity formula for the model $(\GL_6,\GL_2\ltimes U)$. Note that in this case, $\hat{M}=\SL_2(\BC)\times \SL_6(\BC)\times \GL_1(\BC)/(\BZ/6\BZ)$. We have a projection map 
$$\hat{M}=\SL_2(\BC)\times \SL_6(\BC)\times \GL_1(\BC)/(\BZ/6\BZ)$$
$$\rightarrow \SL_6(\BC)\times \GL_1(\BC)/(\BZ/6\BZ)=\GL_6(\BC).$$
Hence each Langlands parameter of $M$ induces a Langlands parameter of $\GL_6$, which allows us to apply Conjecture \ref{conj GL(6) general} for the model $(\GL_6,\GL_2\ltimes U)$.

If $M$ is of Type $A_3\times A_2\times A_1$, then Conjecture \ref{weak conjecture} follows from Proposition \ref{E7 parabolic induction}, Conjecture \ref{conj GL(4)xGL(2)} (see Remark \ref{E7 imply D6 GL(6)}), and the multiplicity formula for the model $(\GL_4\times \GL_2,\GL_2\times \GL_2)$.  Note that in this case, $\hat{M}=\SL_2(\BC)\times \SL_3(\BC)\times \SL_4(\BC)\times \GL_1(\BC)/(\BZ/12\BZ)$. We have a projection map 
$$\hat{M}=\SL_2(\BC)\times \SL_3(\BC)\times \SL_4(\BC)\times \GL_1(\BC)/(\BZ/12\BZ)$$
$$\rightarrow \SL_2(\BC)\times \SL_4(\BC)/(\BZ/4\BZ)=\GL_4(\BC)\times \GL_2(\BC)/\{(aI_4,a^2I_2) \mid a\in \BC^{\times}\}.$$
Combining with Theorem 8.1 of \cite{La}, each Langlands parameter of $M$ induces a Langlands parameter of $\GL_4\times \GL_2$, this allows us to apply Conjecture \ref{conj GL(4)xGL(2)} for the model $(\GL_4\times \GL_2,\GL_2\times \GL_2)$.

If $\Pi_\phi(G)$ is discrete, since $|\Pi_\phi(G)|>1$, there exists a proper elliptic extended endoscopic triple $(G',s',{}^L\eta)$ of $G$ such that $\phi$ factors through ${}^L\eta$ and $s'\in Z_\phi$. We can view $\phi$ as a Langlands parameter of $G'$. If $\hat{G}'=\SL_8(\BC)/\BZ_2$, Proposition \ref{prop E7} implies that
$$\sum_{\pi\in \Pi_\phi(G)}\tr(\chi_\pi(s'))m(\pi)=1,\;\sum_{\pi\in \Pi_\phi(G_D)}\tr(\chi_\pi(s'))m(\pi)=0.$$
In this case, by the discussion in Section \ref{sec epsilon factor} we also know that 
$$\epsilon(\frac{1}{2},\Pi_\phi,\rho_X)=1.$$ 
This proves Conjecture \ref{weak conjecture}.

If 
\begin{eqnarray*}
\hat{G}'&=&\SL_6(\BC)\times \SL_3(\BC)/\BZ_3=\{(g_1,g_2)\in \GL_6(\BC)\times \GL_3(\BC) \mid \\
&&\det(g_1)=\det(g_2)^4\}/\{(a^2I_6,aI_3) \mid a\in \GL_1(\BC)\}, 
\end{eqnarray*}
by Theorem 8.1 of \cite{La},  we can lift a Langlands parameter of $G'$ to a Langlands parameter of $\GL_6\times \GL_3$. Then by Proposition \ref{prop E7}, Conjecture \ref{conj GL(6) general} and the multiplicity formula of the model $(\GL_6,\GL_2\ltimes U)$, we have (recall that $\varepsilon(s')=1$ if the order of $s'$ is 3 and it is equal to $-1$ if the order of $s'$ is 6)
$$\sum_{\pi\in \Pi_\phi(G)}\tr(\chi_\pi(s'))m(\pi)=\frac{1+\epsilon(\frac{1}{2},\rho_{s',\phi})}{2},$$
$$\sum_{\pi\in \Pi_\phi(G_D)}\tr(\chi_\pi(s'))m(\pi)=-\varepsilon(s')\cdot\frac{\epsilon(\frac{1}{2},\rho_{s',\phi})-1}{2},$$
where $\rho_{s',\phi}$ is defined in Section \ref{sec epsilon factor}. 
In particular, we know that the unique element belongs to the packet $\Pi_\phi(G)$ if and only if $\epsilon(\frac{1}{2},\rho_{s',\phi})=1$. By our discussion in Section \ref{sec epsilon factor}, we have $\epsilon(\frac{1}{2},\rho_{s',\phi})=\epsilon(\frac{1}{2},\Pi_\phi,\rho_X)$. This proves Conjecture \ref{weak conjecture}.

If 
\begin{eqnarray*}
\hat{G}'&=&\SL_4(\BC)\times \SL_4(\BC)\times \SL_2(\BC)/\BZ_4\\
&=&\{(g_1,g_2,g_3)\in \GL_4(\BC)\times \GL_4(\BC)\times \GL_2(\BC) \mid \det(g_1)=\det(g_2)\\
&&=\det(g_3)\}/\{(aI_4,aI_4,a^2I_2) \mid a\in \GL_1(\BC)\},  
\end{eqnarray*}
by Theorem 8.1 of \cite{La},  we can lift a Langlands parameter of $G'$ to a Langlands parameter of $\GL_4\times \GL_4\times \GL_2$. Then by Proposition \ref{prop E7}, Conjecture \ref{conj GL(4)xGL(2)} and the multiplicity formula of the model $(\GL_4\times \GL_2,\GL_2\times \GL_2)$, we have (recall that $\varepsilon(s')=-1$ if $s'=(I_4,\pm iI_4,I_2)$ and $\varepsilon(s')=1$ if $s'=(I_4,\pm iI_4,-I_2)$)
$$\sum_{\pi\in \Pi_\phi(G)}\tr(\chi_\pi(s'))m(\pi)=\frac{\epsilon(\frac{1}{2},\rho_{s',\phi,1}\oplus \rho_{s',\phi,2})}{2},$$
$$\sum_{\pi\in \Pi_\phi(G_D)}\tr(\chi_\pi(s'))m(\pi)=\varepsilon(s')\cdot\frac{\epsilon(\frac{1}{2},\rho_{s',\phi,1})-\epsilon(\frac{1}{2},\rho_{s',\phi,2})}{2}$$
where $\rho_{s',\phi,i}$ is defined in Section \ref{sec epsilon factor}. In particular, we know that the unique element belongs to the packet $\Pi_\phi(G)$ if and only if 
$$\epsilon(\frac{1}{2},\rho_{s',\phi,1}\oplus \rho_{s',\phi,2})=1.$$ 
By our discussion in Section \ref{sec epsilon factor}, we have $$\epsilon(\frac{1}{2},\rho_{s',\phi,1}\oplus \rho_{s',\phi,2})=\epsilon(\frac{1}{2},\Pi_\phi,\rho_X).$$ 
This proves Conjecture \ref{weak conjecture}.

If 
\begin{eqnarray*}
\hat{G}'&=&\Spin_{12}(\BC)\times \SL_2(\BC)/\BZ_2=\{(g_1,g_2)\in \GHSpin_{12}^{\vee}(\BC)\times \GL_2(\BC) \mid  \\
&&l(g_1)\det(g_2)=1\}/\GL_1(\BC)^{anti-diag},
\end{eqnarray*}
by Theorem 8.1 of \cite{La},  we can lift a Langlands parameter of $G'$ to a Langlands parameter of $\GHSpin_{12}\times \GL_2$. Then by Proposition \ref{prop E7}, Conjecture \ref{weak conjecture for E7 smaller model} and the multiplicity formula of the models $(\GHSpin_{12},\PGL_2\ltimes U)$ and $(\GHSpin_{12}\times \GSpin_3,\GSpin_3\ltimes U)$, we have 
$$\sum_{\pi\in \Pi_\phi(G)}\tr(\chi_\pi(s'))m(\pi)=\frac{\epsilon(\frac{1}{2},\rho_{s',\phi,+}\oplus \rho_{s',\phi,-})}{2},$$
$$\sum_{\pi\in \Pi_\phi(G_D)}\tr(\chi_\pi(s'))m(\pi)=\frac{\epsilon(\frac{1}{2},\rho_{s',\phi,-})-\epsilon(\frac{1}{2},\rho_{s',\phi,+})}{2}$$
where $\rho_{s',\phi,+}$ and $\rho_{s',\phi,-}$ are defined in Section \ref{sec epsilon factor}. In particular, we know that the unique element belongs to the packet $\Pi_\phi(G)$ if and only if $\epsilon(\frac{1}{2},\Pi_\phi,\rho_X)=\epsilon(\frac{1}{2},\rho_{s',\phi,+}\oplus \rho_{s',\phi,-})=1$. This finishes the proof of Conjecture \ref{weak conjecture} when $\Pi_\phi$ is not discrete with $|\Pi_\phi(G)|=1$.

Now we are ready to prove Theorem \ref{main theorem for E7}. Let $\omega_\phi\in \hat{S}_\phi$ correspond  to the unique distinguished element in the packet. By Remark \ref{distinguished is character} we know that $\omega_\phi$ is a character and we view it as a character of $Z_\phi$. For $s\in S_\phi$, by Lemma \ref{lem extended endoscopic triple}, there exists an elliptic extended endoscopic triple $(G',s',{}^L\eta)$ of $G/Z_{G,H}$ such that $s'\in sZ_{\phi}^{\circ}$ and $\phi$ factors through ${}^L\eta$. We need to show that $\omega_\phi(s')=\omega_{\phi,H}(s)$. The above discussion implies that $\omega_\phi(s')=\omega_{\phi,H}(s)$ if $s'$ belongs to the center of the dual group.

If $s'$ does not belong to the center of the dual group, there are four cases. If $\hat{G}'=\SL_8(\BC)/\BZ_2$, the above discussion implies that the unique distinguished element belongs to $\Pi_\phi(G)$ and 
$$\sum_{\pi\in \Pi_\phi(G)} \tr(\chi_\pi(s'))m(\pi)=1.$$
By the definition of $\omega_{\phi,H}$ we know that $\omega_{\phi,H}(s)=1$. This implies that
$$\omega_\phi(s')=\tr(\omega_\phi(s'))=\sum_{\pi\in \Pi_\phi(G)} \tr(\chi_\pi(s'))m(\pi)=1=\omega_{\phi,H}(s).$$

If $\hat{G}'=\SL_6(\BC)\times \SL_3(\BC)/\BZ_3$, by our discussion above, we have 
$$\sum_{\pi\in \Pi_\phi(G)}\tr(\chi_\pi(s'))m(\pi)=\frac{1+\epsilon(\frac{1}{2},\rho_{s',\phi})}{2},$$
$$\sum_{\pi\in \Pi_\phi(G_D)}\tr(\chi_\pi(s'))m(\pi)=-\varepsilon(s')\cdot\frac{\epsilon(\frac{1}{2},\rho_{s',\phi})-1}{2}.$$
By the definition in Section \ref{sec epsilon factor}, we have $\omega_{\phi,H}(s)=\epsilon(\frac{1}{2},\rho_{s',\phi})$ if the order of $s'$ is 6 and $\omega_{\phi,H}(s)=1$ if the order of $s'$ is 3. We have two cases. If the unique distinguished element belongs to $\Pi_\phi(G)$, we have $\epsilon(\frac{1}{2},\rho_{s',\phi})=1$. This implies that
$$\omega_\phi(s')=\tr(\omega_\phi(s'))=\sum_{\pi\in \Pi_\phi(G)} \tr(\chi_\pi(s'))m(\pi)$$
$$=\frac{1+\epsilon(\frac{1}{2},\rho_{s',\phi})}{2}=\epsilon(\frac{1}{2},\rho_{s',\phi})=1=\omega_{\phi,H}(s).$$
If the unique distinguished element belongs to $\Pi_\phi(G_D)$, we have $\epsilon(\frac{1}{2},\rho_{s',\phi})=-1.$ This implies that
$$\omega_\phi(s')=\tr(\omega_\phi(s'))=\sum_{\pi_D\in \Pi_\phi(G_D)} \tr(\chi_{\pi_D}(s'))m(\pi_D)$$
$$=-\varepsilon(s')\epsilon(\frac{1}{2},\rho_{s',\phi})=\omega_{\phi,H}(s).$$

If $\hat{G}'=\SL_4(\BC)\times \SL_4(\BC)\times \SL_2(\BC)/\BZ_4$, by our discussion above, we have 
$$\sum_{\pi\in \Pi_\phi(G)}\tr(\chi_\pi(s'))m(\pi)=\frac{\epsilon(\frac{1}{2},\rho_{s',\phi,1}\oplus \rho_{s',\phi,2})}{2},$$
$$\sum_{\pi\in \Pi_\phi(G_D)}\tr(\chi_\pi(s'))m(\pi)=\varepsilon(s')\cdot\frac{\epsilon(\frac{1}{2},\rho_{s',\phi,1})-\epsilon(\frac{1}{2},\rho_{s',\phi,2})}{2}.$$
By the definition in Section \ref{sec epsilon factor}, we have $\omega_{\phi,H}(s)=\epsilon(\frac{1}{2},\rho_{s',\phi,2})$ if the $s'=(I_4,\pm iI_4,I_2)$ and $\omega_{\phi,H}(s)=\epsilon(\frac{1}{2},\rho_{s',\phi,1})$ if the $s'=(I_4,\pm iI_4,-I_2)$. We have two cases. If the unique distinguished element belongs to $\Pi_\phi(G)$, we have $\epsilon(\frac{1}{2},\rho_{s',\phi,1})=\epsilon(\frac{1}{2},\rho_{s',\phi,2})=\omega_{\phi,H}(s)$. This implies that
$$\omega_\phi(s')=\tr(\omega_\phi(s'))=\sum_{\pi\in \Pi_\phi(G)} \tr(\chi_\pi(s'))m(\pi)$$
$$=\frac{\epsilon(\frac{1}{2},\rho_{s',\phi,1}\oplus \rho_{s',\phi,2})}{2}=\omega_{\phi,H}(s).$$
If the unique distinguished element belongs to $\Pi_\phi(G_D)$, we have $\epsilon(\frac{1}{2},\rho_{s',\phi,1})=-\epsilon(\frac{1}{2},\rho_{s',\phi,2})$. 
This implies that
$$\omega_\phi(s')=\tr(\omega_\phi(s'))=\sum_{\pi_D\in \Pi_\phi(G_D)} \tr(\chi_{\pi_D}(s'))m(\pi_D)$$
$$=\varepsilon(s')\cdot\frac{\epsilon(\frac{1}{2},\rho_{s',\phi,1})-\epsilon(\frac{1}{2},\rho_{s',\phi,2})}{2}=\omega_{\phi,H}(s).$$

If $\hat{G}'=\Spin_{12}(\BC)\times \SL_2(\BC)/\BZ_2$, by our discussion above, we have 
$$\sum_{\pi\in \Pi_\phi(G)}\tr(\chi_\pi(s'))m(\pi)=\frac{\epsilon(\frac{1}{2},\rho_{s',\phi,+}\oplus \rho_{s',\phi,-})}{2},$$
$$\sum_{\pi\in \Pi_\phi(G_D)}\tr(\chi_\pi(s'))m(\pi)=\frac{\epsilon(\frac{1}{2},\rho_{s',\phi,-})-\epsilon(\frac{1}{2},\rho_{s',\phi,+})}{2}.$$
By the definition in Section \ref{sec epsilon factor}, we have $\omega_{\phi,H}(s)=\epsilon(\frac{1}{2},\rho_{s',\phi,-})$. We have two cases. If the unique distinguished element belongs to $\Pi_\phi(G)$, we have $\epsilon(\frac{1}{2},\rho_{s',\phi,+})\epsilon(\frac{1}{2},\rho_{s',\phi,-})=1$, i.e. $\epsilon(\frac{1}{2},\rho_{s',\phi,+})=\epsilon(\frac{1}{2},\rho_{s',\phi,-})$. This implies that
$$\omega_\phi(s')=\tr(\omega_\phi(s'))=\sum_{\pi\in \Pi_\phi(G)} \tr(\chi_\pi(s'))m(\pi)$$
$$=\frac{\epsilon(\frac{1}{2},\rho_{s',\phi,+}\oplus \rho_{s',\phi,-})}{2}=\epsilon(\frac{1}{2},\rho_{s',\phi,-})=\omega_{\phi,H}(s).$$
If the unique distinguished element belongs to $\Pi_\phi(G_D)$, we have $$\epsilon(\frac{1}{2},\rho_{s',\phi,+})\epsilon(\frac{1}{2},\rho_{s',\phi,-})=-1\Rightarrow \epsilon(\frac{1}{2},\rho_{s',\phi,+})=-\epsilon(\frac{1}{2},\rho_{s',\phi,-}).$$ 
This implies that
$$\omega_\phi(s')=\tr(\omega_\phi(s'))=\sum_{\pi_D\in \Pi_\phi(G_D)} \tr(\chi_{\pi_D}(s'))m(\pi_D)$$
$$=\frac{\epsilon(\frac{1}{2},\rho_{s',\phi,-})-\epsilon(\frac{1}{2},\rho_{s',\phi,+})}{2}=\epsilon(\frac{1}{2},\rho_{s',\phi,-})=\omega_{\phi,H}(s).$$
This finishes the proof of Theorem \ref{main theorem for E7}.

Lastly, we prove Theorem \ref{thm weak conj for E7}. Assume that the Conjecture \ref{weak conjecture} holds for the model $(E_7,\PGL_2\ltimes U)$. We need to prove Conjecture \ref{weak conjecture for E7 smaller model}. By Remark \ref{E7 imply D6 GL(6)} and Remark \ref{E7 imply D5} we know that Conjecture \ref{weak conjecture for E7 smaller model} holds for the models $(\GHSpin_{12},\PGL_2\ltimes U)$ and $(\GSpin_{10}\times \GSpin_3,\GSpin_3\ltimes U)$. It remains to prove it for the model $(\GHSpin_{12}\times \GL_2,\GL_2\ltimes U)$. We just need to use the endoscopic relation in Proposition \ref{prop E7} for the case when $G'$ is of Type $D_6\times A_1$ together with Conjecture \ref{weak conjecture for E7 smaller model} for the model $(\GHSpin_{12},\PGL_2\ltimes U)$. The argument is the same as the proof of Theorem \ref{thm weak conjecture smaller models} for the model $(\GU_6,\GU_2\ltimes U)$ in Section \ref{sec GU(6)} and we will skip it here. This completes the proof of Theorem \ref{thm weak conj for E7}.

\end{document}